%% file: QCD_TIT_Aug17_arXiv.tex
\documentclass{article} 

\pdfoutput=1  

\usepackage{amsmath,amsfonts}
\usepackage{algorithmic}
\usepackage{algorithm}
\usepackage{array}
\usepackage[caption=false,font=normalsize,labelfont=sf,textfont=sf]{subfig}
\usepackage{textcomp}
\usepackage{stfloats}
\usepackage{url}
\usepackage{verbatim}

\usepackage{graphicx}
  
\usepackage{cite}

\usepackage{geometry}
\usepackage{amssymb}
\usepackage{bm,upgreek}

\usepackage[dvipsnames]{xcolor}

\usepackage[pdfa,hidelinks]{hyperref}
\usepackage{cleveref}

\usepackage{titlesec}
	 \usepackage{accents}

 \usepackage{wrapfig}

\newtheorem{theorem}{Theorem}[section]
\newtheorem{lemma}[theorem]{Lemma}
\newtheorem{proposition}[theorem]{Proposition}
\newtheorem{corollary}[theorem]{Corollary}

\newlength{\noteWidth}
\long\def\notes#1{\ifinner
	{\tiny #1}
\else
\marginpar{\parbox[t]{\noteWidth}{\raggedright\tiny #1}}
\fi}
\def\notes#1{}

\Crefname{corollary}{Corollary}{Corollaries}
\Crefname{eqnarray}{eq.}{eqs.}
\Crefname{equation}{eq.}{eqs.}
\Crefname{figure}{Fig.}{Figs.}
\Crefname{tabular}{Tab.}{Tabs.}
\Crefname{table}{Tab.}{Tabs.}
\Crefname{lemma}{Lemma}{Lemmas}
\Crefname{theorem}{Thm.}{Thms.}
\Crefname{definition}{Definition}{Definitions}
\Crefname{section}{Section}{Sections}
\Crefname{proposition}{Prop.}{Propositions}

\input{bookmacrosRL_IEEE.tex}

\def\PF{\wham{Proof:}}

\def\Ass#1{A.\ref{#1}} 
\def\thone{\theta^\wr}
\def\thone{\theta^\ddagger}

\def\thzz{\theta^\varsigma}

\def\thmix{\theta^\wr}
\def\Fmix{F^\wr}

\newcommand{\clThm}{%
  \ifvmode\noindent\fi
  \unskip
  \nobreak\hfill\penalty50
  \hskip1em\hbox{}\nobreak\hfill
  \hbox{\footnotesize$\square$}%
  \par
  \gobblepars
}
 
\def\Pdens{\upvarpi}

\def\kmax{k_{\textsf{max}}}

\def\Fmax{F_{\textsf{max}}}

 \def\filt{\textsf{T}}  

\def\Tmin{T_{\textsf{min}}}

\def\bst{\upomega}

\def\cg{{\check{g}}}

\def\cX{{\check{X}}}
\def\cY{{\check{Y}}}

\def\whamb{\wham{$\bullet$} }

\def\Wham#1{\wham{\large#1}}

\DeclareMathAccent{\widecheck}{0}{mathx}{"71}

\def\xH{x^{\text{\sf\tiny(H)}}}
\def\sH{s^{\text{\sf\tiny(H)}}}

\def\whamit#1{\smallbreak\pagebreak[3]%
\noindent\textit{#1}\ \ \gobblepars}

\def\wham#1{\smallbreak\pagebreak[3]%
\noindent\textbf{#1}\ \ \gobblepars}

\newcommand{\overbar}[1]{\mkern 1.5mu\overline{\mkern-1.5mu#1\mkern-1.5mu}\mkern 1.5mu}

\def\tilG{\widetilde{G}}

\def\pFA{p_{\textsf{\tiny FA}}}

\def\del{\updelta} 

\def\PFeig{\upzeta}		%

\def\Obs{Y}
\def\bfObs{\bm{Y}}

\def\strY#1{\Obs^{(#1)}}

\def\stroX#1{\accentset{\circ}{X}^{(#1)}}

\def\preObs{X^0}
\def\postObs{X^1}
\def\bfpreObs{\bfmX^0}
\def\bfpostObs{\bfmX^1}

\def\condDist{\Uppi}

\def\InfoState{\clX}
\def\cm{\widecheck{m}}
\def\thexp{\upupsilon}

\def\twmarg{\widecheck{\uppi}}
\def\twtrajP{\widecheck{\upmu}}

\def\marg{\uppi} 
 
\def\cmarg{\widecheck{\uppi}}

\def\trajP{\upmu}

\def\metamarg{\upgamma}

\def\ometamarg{\check{\upgamma}}

\def\tchange{\uptau_{\text{\scriptsize\sf a}}}
\def\tstop{\uptau_{\text{\scriptsize\sf s}}}

\def\expa{\varrho_{\text{\scriptsize\sf a}}}

\def\thresh{%
\mathchoice
	{\text{\small\rm H}}%
	{\text{\small\rm H}}%
	{\text{\scriptsize\rm H}}%
	{\text{\tiny\rm H}}}

\def\barthresh{\overbar{\text{\small\rm H}}}

\def\threshSmall{%
\mathchoice
	{\text{\small\rm h}}%
	{\text{\small\rm h}}%
	{\text{\scriptsize\rm h}}%
	{\text{\tiny\rm h}}
}

 \def\threshSmallStar{\threshSmall^{\!*}}

\def\fineTag{\blacktriangle}

\def\coarseTag{\infty}

\def\barthreshFineB{\barthresh^\circ_{\fineTag}}
\def\barthreshFineStar{\barthresh^*_{\fineTag}}

\def\barJFineB{\barJ^\circ_{\!\fineTag}}
\def\barJFineStar{\barJ^*_{\!\fineTag}}

\def\xkFine{x_\kappa^{\fineTag}}

\def\barJStarCoarse{\barJ_{\coarseTag}^{*}}

\def\aFine{\upalpha}  
\def\bFine{\upbeta}  

\def\Fstar{F^*}

\def\surL{\breve{L}}
\def\surg{\breve{g}}

\def\MDD{{\sf MDD}}
\def\MDE{{\sf MDE}}
\def\MDDfine{{\sf MDD}_\fineTag}
\def\MDEfine{{\sf MDE}_\fineTag}

   	\def\lorFA{{\sf FAR}}
\def\lorMDD{{\sf WADD}}

\title{\LARGE 
Quickest Change Detection Using Mismatched {CUSUM}
}
\author{Austin Cooper and Sean Meyn%
\thanks{AC and SM are with the  Department of Electrical and Computer Engineering,  University of Florida, Gainesville,  FL, USA;  austin.cooper,meyn @ufl.edu.   Financial support from   NSF awards  CCF-2306023 and DMS-2427265.}}

\begin{document}
	
\maketitle

   \begin{abstract}
Quickest change detection concerns estimation of an unknown change time
\(\tchange\) from a sequence of partial observations
\(\{Y_k:k\ge 0\}\).  We consider stopping rules of CUSUM form,
\[
        \InfoState_{n+1}
        =
        \max\{0,\InfoState_n+F(Y_{n+1})\},
        \quad
        \tstop=\min\{n\ge 0:\InfoState_n\ge \thresh\},
\]
where the function \(F\) and threshold \(\thresh\) are design
parameters.

The observations and change time are modeled jointly through a hidden Markov
model, and \( F\) is selected from a prescribed function class \(\clG\) to
minimize the weighted criterion
\[
        \Expect\bigl[
            (\tstop-\tchange)_+
            +
            \kappa(\tstop-\tchange)_-
        \bigr].
\]
 When \(\clG\) is a linear function class, the optimizer
  \(F^*\) is characterized by a convex program, whose dual yields
extensions of classical likelihood-ratio constructions.
This conclusion is based on analysis that is asymptotic in the regime \(\kappa\to\infty\).      
We show that the hidden Markov model admits an asymptotically equivalent conditionally independent approximation of the type commonly used in the
quickest change detection literature.  We then develop the design and
asymptotic theory for a substantially broader class of conditionally
independent models, so that the resulting conclusions are not tied to the
particular POMDP reduction. 

 Combining renewal theory and large deviations
for reflected random walks, we obtain for each $F\in\clG$ asymptotically accurate approximations
of the optimal threshold and average cost, with error vanishing as
\(\kappa\to\infty\).   It is found in numerical
experiments that the resulting approximations are accurate for
moderate values of \(\kappa\).

\end{abstract}

 %
%


\def\oshoot{\clV}	
\def\bigstate{\Upomega} 
\def\Noise{B}
\def\bfNoise{\bfmath{B}}

\def\oX{\accentset{\circ}{X}} 

\def\oPhi{\accentset{\circ}{\Phi}} 
\def\bfoPhi{\accentset{\circ}{\bfPhi}}

  \def\oP{\accentset{\circ}{P}} 
  \def\oQ{\accentset{\circ}{Q}} 
  \def\oY{\accentset{\circ}{Y}} 
  \def\oq{\accentset{\circ}{q}} 
  \def\op{\accentset{\circ}{p}} 
  
  \def\ogamma{\accentset{\circ}{\gamma}} 

  \def\oZ{\accentset{\circ}{Z}}

\def\owmetamarg{\accentset{\circ}{\metamarg}}

\def\ometamarg{\accentset{\circ}{\metamarg}}   

\def\ocondDist{\accentset{\circ}{\condDist}}

\section{Introduction}
\label{s:intro}
	
The  QCD problem  considers  the  real time monitoring of an observation process $\bfObs\eqdef \{ \Obs_k :  k\ge 0\}$ evolving on a set $\ystate$, taken in this paper to be a subset of Euclidean space.    At an unknown  ``change time'' denoted by $\tchange$, the statistics of the observations  undergo a change due to underlying  anomalous behavior.  Examples include the onset of a heart attack,  a fault in a mechanical device,  or a security breach in a computer system or physical environment.

The goal  is to construct a ``stopping time",  denoted $\tstop$ adapted to the observations: on denoting $\Obs_0^k = (\Obs_0; \cdots; \Obs_k)$,  for each $k$ we may write $\ind\{\tstop \le k\} =  \fee_k(\Obs_0^k)$ for some Borel-measurable mapping $\fee_k\colon \ystate^{k+1} \to \{0,1\}$. We refer to any such sequence $\fee = \{ \fee_k : k \ge 0 \}$ as the \textit{policy}.   

The present paper follows a large literature on algorithm design approached through the construction of a real-valued stochastic process $\{\InfoState_n\} $  adapted to the observations, sharing a role similar to the information state of POMDP theory.
	The stopping rule is  of the threshold form,
\begin{equation}
\tstop = \min\{ n\ge 0 :   \InfoState_n\ge \thresh \}\, ,
\label{e:threshold}
\end{equation}
	with $\thresh>0$.    
	Two examples are the tests of Shiryaev–Roberts, and Page's CUSUM.    The latter is the focus of this paper, in which    $\{\InfoState_n \}$ is defined as a reflected random walk,  
\begin{equation}
\InfoState_{n+1} =    \max\{ 0,  \InfoState_n +  F_{n+1}  \} 
\label{e:CUSUM}
\end{equation}
initialized with $\InfoState_0=0$,   where $\{ F_{n+1} : n\ge 0\}$  is a stochastic process adapted to the observations. 

  \begin{subequations}
  
Any approach to policy design  must balance  two costs:    1.~\textit{Delay}, expressed $(\tstop -\tchange)_+\eqdef \max(0,\tstop -\tchange)$,  and 2.~\textit{false alarm}, meaning that $\tstop -\tchange<0$.     A standard performance metric for a policy $\fee$ in a Bayesian setting is the sum
\begin{equation}
		J_0(\fee)  = \MDD+\kappa\pFA  
\label{e:MDD+kappa_pFA}
\end{equation} 
in which $\kappa >0$,  $ \pFA=  \Prob\{\tstop < \tchange \} $ denotes the  \textit{probability of false alarm},
 and the \textit{mean detection delay} is 
\begin{equation}
\MDD  =   \Expect[ (\tstop -\tchange)_+ ]  \,, 
\label{e:MDDpFA}
\end{equation}
	with   $x_+=\max(x,0)$, $x_-=\max(-x,0)$ for $x\in\Re$. 
	
  \end{subequations}

The following representation for the QCD model is adopted throughout this work:
\begin{equation}
\Obs_k = \preObs_k \ind_{k< \tchange } +  \postObs_k \ind_{k \ge \tchange } \,, \qquad k \ge 0.
\label{e:QCDmodel}
\end{equation} 
Much of the  literature imposes the following assumption:

\wham{Conditionally independent (C.I.) model:}  The observation equation \eqref{e:QCDmodel}
in which  $\bfpreObs$, $\bfpostObs$ and $\tchange$ are statistical independent.

 This leads to elegant theory, but unfortunately  excludes many applications of interest.   A simple example makes this clear: engine failure in an automobile.  The driver may hear noise,  followed by a gradual loss in power,  and an hour later the engine freezes.   If $\tchange$ is taken to be the time at which the engine fails,  then there is a great deal of information regarding this change time based on the pre-change behavior.

This example not only illustrates lack of independence in a representative application, but also a flaw in the performance criterion \eqref{e:MDD+kappa_pFA}:  if the policy commonly induces a false alarm, but the \textit{detection eagerness} $  (\tstop -\tchange)_- $ is small,  then the driver of the vehicle will be pleased---this policy may be more valuable than one that always achieves $\tstop -\tchange = 0$.

For these reasons, in this paper we modify the objective and allow statistical dependence between $\tchange$ and $\bfmX^i$ 
 in the observation model \eqref{e:QCDmodel}:

\wham{1.}   The probability of false alarm $\pFA$  is replaced by the
\textit{mean detection eagerness} $\MDE=\Expect[ (\tstop -\tchange)_- ]$:   
\begin{equation}
		J(\fee)  = \MDD+  \kappa\MDE =   \Expect\big[ (\tstop -\tchange)_+  +  \kappa (\tstop -\tchange)_-  \big]
\label{e:MDD+kappaMDE}
\end{equation} 

\wham{2.}    The observation process and change time are modeled as a hidden Markov model (HMM) with hidden state $\bfPhi$ evolving on a state space $\bigstate$.  The model is constructed so that the performance criterion \eqref{e:MDD+kappaMDE} may be expressed in a form standard in optimal control theory:
\begin{equation}
			J(\Phi_0, U_0^\infty ) =   \Expect\Big[ \sum_{k=0}^{\tstop -1}  c_\circ(\Phi_k  )     + c _\bullet(\Phi_{\tstop}  )  \Big]
\label{e:ObjPOMDPQCD}
\end{equation} 
An   optimal policy can be expressed as state feedback based on an \textit{information state};  this is most commonly prescribed as the sequences of conditional distribution of $\Phi_k$ given observations $Y_0^k$.     Computation of an optimal policy is intractable outside of very special cases such as Shiryaev's model,  for which the information state reduces to a scalar stochastic process  \cite{shi77}.

\begin{figure*}[h]
	\centering
	\includegraphics[width=0.95\hsize]{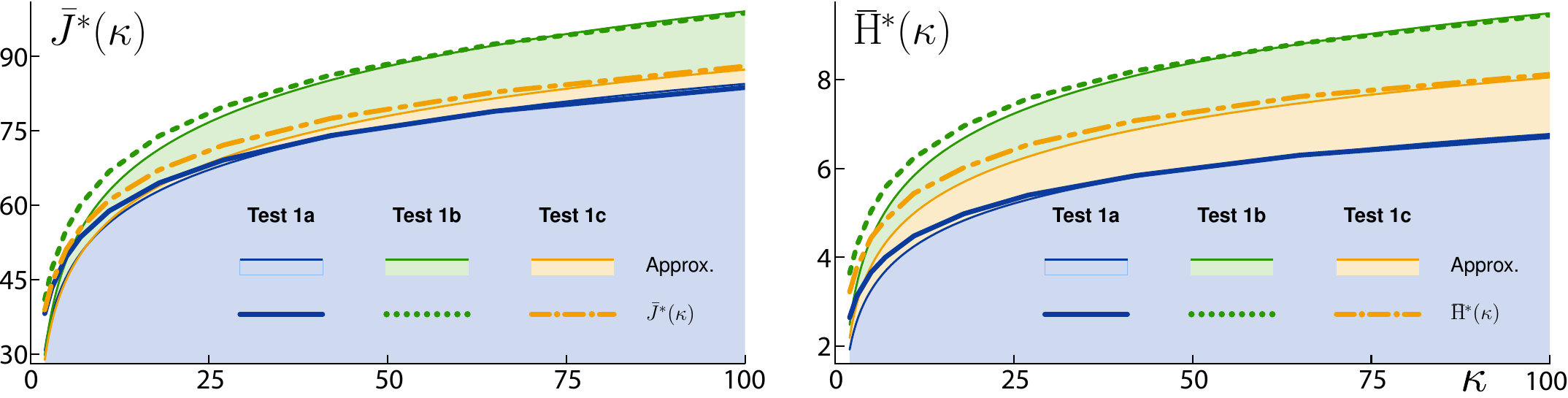}  
	\caption{CUSUM* performance and thresholds for Model~1 compared to the approximations in     \eqref{e:barFineBody}}
		\label{f:model1_approx}
\end{figure*}

Denote by $\barJ(\thresh, \kappa)  $  the value of the expectation \eqref{e:MDD+kappaMDE} using CUSUM with threshold $\thresh>0$,   and    
\begin{equation}
\barJ^* (\kappa)  =\min_{\thresh} \barJ( \thresh, \kappa)  \,,
\qquad
\barthresh^*(\kappa)  =\argmin_{\thresh} \barJ(\thresh, \kappa)  
\label{e:Best_H_J}
\end{equation}  
The CUSUM algorithm using the optimal threshold $\barthresh^*(\kappa)  $ is denoted CUSUM*.

\wham{Contributions}   The paper focuses on approximations that are asymptotic in the penalty term $\kappa$ appearing in \eqref{e:MDD+kappaMDE}.
The contributions are of three varieties:

\wham{1.}  \textit{Strong approximation by a conditionally independent model}.    
The POMDP model is introduced in 
\Cref{s:POMDP}.  It is shown in \Cref{t:POMDPhazard} the observation process can be approximated by the conditionally independent model in the large-$\kappa$ regime.

\wham{2.}    \textit{Performance approximation for CUSUM*}.  
The following is obtained in \Cref{t:barcApprox} under
   mild assumptions:  
\begin{equation}
	\barthresh^*(  \kappa)    =  \barthreshFineStar(\kappa)  +   o( 1)	 \,,
	\qquad
	\barJ^*(  \kappa)    =  \barJFineStar ( \kappa) 
				+  o( 1)  \,,
  \label{eq:refined_bdds} 
\end{equation}
\begin{subequations}%
 where the terms $o(1)$ vanish as $\kappa \to\infty$,   
 \begin{align}
\barthreshFineStar(\kappa)
			&=
\Bigl[\log\kappa+\frac{1}{2}\log\log\kappa+\bFine\Bigr]   \frac{1}{\thexp_+}
\label{e:barHFineBody}
\\
\barJFineStar(\kappa)
			&=
\Bigl[\log\kappa+\frac{1}{2}\log\log\kappa+ \aFine\Bigr] 
\frac{1}{m_1\thexp_+}  \,, 
\label{e:barJFineBody}
\end{align}
 $\thexp_+ >0$,  $   \aFine ,\bFine \in\Re$ are identified in the proposition, and $m_1 = \lim_{n\to\infty}\Expect[F_n \mid \tchange \le n ]$ is assumed positive.
   \label{e:barFineBody}
    \end{subequations} 
    
    \smallskip

\Cref{f:model1_approx} provides an illustration of the approximations in \eqref{eq:refined_bdds} for a particular model, showing a tight approximation for $\kappa\ge 25$  for each of the three variants of CUSUM*   considered---details  may be found in \Cref{s:num_model1}.

\wham{3.}  
 \textit{Optimality and convex duality.}  
 \Cref{s:thetaCUSUM}  concerns optimization of the increment process $\{ F_{n+1} : n\ge 0\}$ appearing in \eqref{e:CUSUM}.   
 We take $F_n = F(Y_n)$ in \eqref{e:CUSUM}, and  seek to  
minimize  over $F\in\clG$ the cost approximation
\begin{align}
\barJStarCoarse( \kappa ; F)  &=\Big[  \log \kappa    +\frac{1}{2}\log\log\kappa\Bigr]  
\frac{1}{ \marg_1(F) }\frac{1}{\thexp_+(F)}  
\label{e:barJ_infty_F}
\end{align}
where  $\thexp_+(F)>0$ is the solution to $\Lambda_0( \thexp F )  = \expa$.

It is shown in   \Cref{t:calcOpt} that the first order condition for optimality reduces to a moment matching problem.  

Minimization of \eqref{e:barJ_infty_F}
 may be expressed as a convex program, for which  strong duality is established in \Cref{t:DualProps} along with an information-theoretic representation of the dual function.


\notes{Tweak for arxiv}
 
\smallbreak

The conclusions in \Cref{s:thetaCUSUM} both unify and extend existing theory.    
Consider the C.I.\ Markov model in which $\bfmX^i$ is a time-homogeneous Markov chain for each $i$,
and let $L_\infty\colon \ystate^2\to\Re$ denote the log likelihood of the respective transition kernels (recalled in   \Cref{s:trad}). 
The following conclusions are well known  using $F_{n+1} = L_\infty (Y_n, Y_{n+1})$:
\textbf{(a)} 
 CUSUM* is approximately optimal for large $\kappa$,
subject to the performance criterion \eqref{e:MDD+kappa_pFA}.    
\textbf{(b)} 
Sample path large deviations theory tells us that conditioned on a false alarm, the
statistics of $ \{ Y_n  : n \le \tstop\}$ are approximated by those of the   Markov chain    $\bfmX^1$.   

\smallskip
Parallel conclusions are obtained in  \Cref{s:thetaCUSUM}:

\wham{(a')}
It is shown in  \Cref{t:bestF} that the choice   
$F_{n+1} = L_\infty (Y_n, Y_{n+1})  + \expa$ is approximately optimal with respect to the   performance criterion  \eqref{e:MDD+kappaMDE}, 
 with $\expa$  the \textit{survival exponent},  
 \begin{equation}
 \expa \eqdef -
\lim_{n\to\infty}\frac{1}{n}\log\Prob\{\tchange\ge n\} 
\label{e:hazardAss}
\end{equation}

\wham{(b')}
Suppose that $F_{n+1} = F^* (Y_n, Y_{n+1}) $ where $F^*$ is optimal over the linear function class $\clG$.     
The most likely behavior before false alarm may be represented by  a \textit{twisted law} denoted $\twtrajP_0$:   It defines a stationary Markov chain,  
whose bivariate marginals coincide with those of $\bfmX^1$ on this    function class:  
\begin{equation}
\Expect [ G( \cX^0_{n-1},\cX^0_{n}) ] = 
\Expect [ G( X^1_{n-1},X^1_{n}) ] \,,\quad G\in\clG\,,
\label{e:mmIntro}
\end{equation}
where the expectation on the left hand side uses $\twtrajP_0$.   
This is an instance of the moment-matching property  established in  \Cref{t:calcOpt}.  It  recovers the classical conclusion \textbf{(b)} when the function class $\clG$  is \textit{dense} (see \Cref{t:completeOpt}).

\wham{Literature}

The extended abstract \cite{coomey24d} contains preliminary results, including a crude version of  \eqref{eq:refined_bdds}:   $\barthresh(  \kappa) = \log(\kappa) / \thexp_+ + o( \log(\kappa) )$,  $\barJ(  \kappa)  = \log(\kappa) /(m_1 \thexp_+) + o( \log(\kappa) )$.

Our interest in the objective function \eqref{e:MDD+kappaMDE}  is well aligned with the large literature on optimal decision making for preventative maintenance  in which a goal is timely detection before failure  \cite{iseReview05,ise06,dejsca20,mcc65}.   There are other commonalities:   a POMDP setting is typical and   the existence of the limit \eqref{e:hazardAss} is implied by the hazard rate conditions found in this literature.    

In the Bayesian setting an optimal policy minimizing the objective \eqref{e:ObjPOMDPQCD} may be characterized using POMDP theory, and this approach is successful in special cases such as  Shiryaev's model  \cite{shi77} and generalizations in \cite{yak94,lai95}.    See   \cite{kri11,zhoenl13} for history and numerical techniques to solve this optimization problem  in more general settings.

The complexity of POMDP theory may be avoided by turning to a large-$\kappa$ asymptotic setting.  Much of the literature addresses an adversarial setting:    A test is said to be optimal under Lorden's criterion if it minimizes  the \textit{worst-case average detection delay}  $\lorMDD  $, subject to a lower bound on the  \textit{worst case false alarm rate}   
 $\lorFA \le 1/\kappa$, with  
\begin{equation}
\begin{aligned}
 \lorMDD  &= 
   \sup_n  \Expect[  (\tstop - n)_+ \mid \tchange = n]    
\\
\lorFA &=  \frac{1}{ \Expect[  \tstop \mid  \tchange = \infty]} 
\end{aligned}
\label{e:Lorden}
\end{equation}

CUSUM was introduced in \cite{pag54} following a series of works on sequential hypothesis testing, beginning with \cite{wal45}.   It is known that CUSUM and some other sequential tests are approximately optimal, as $\kappa\to\infty$,  for both the Lorden's criterion and the criterion \eqref{e:MDD+kappa_pFA};   the theory also provides estimates of the optimal threshold as a function of $\kappa$.     Early work can be found in \cite{mou86,rit90,mou98},   based on   \cite{lor71}.   Approximate optimality for the closely related   Shiryaev-Roberts statistic may be found in  \cite{mou98,mostar09,poltar10}.

To obtain approximate optimality the functions $\{F_k :  k\ge 1\}$ are chosen so that  their partial sums approximate a sequence of log likelihood ratios.  Consider the C.I.\   model \eqref{e:QCDmodel} in which $\ystate$ is discrete, and for $n\ge0$ let $p_i^n$ denote the pmf for the random string $(X_0^i,\dots, X_n^i)$, $i=0,1$.    In key results of \cite{mou86,lai95,fuh03} and many others, it is assumed that for each $n$,
\begin{equation}
        \log
        \frac{p_1^n (Y_0,\ldots,Y_n)}
             {p_0^n(Y_0,\ldots,Y_n)}
        =       \log \frac{p^0_1 (Y_0 ) }{p^0_0 (Y_0 ) } +        \sum_{k=1}^{n}  F_k
\label{e:LLRn}
\end{equation}
Approximations of $\lorMDD$ or related performance metrics are expressed in terms of  the relative entropy rate $	\clK(\trajP_1  \|  \trajP_0) $ between the two stochastic processes $\bfmX^0$,   $\bfmX^1$:    \begin{equation}
	\clK(\trajP_1  \|  \trajP_0)  \eqdef \lim_{n\to\infty} \frac{1}{n}   D( p_1^n  \| p_0^n )  
	= \lim_{n\to\infty} \frac{1}{n}   \sum_{k=1}^n F_k  
\label{e:FuhAssumption}
\end{equation}
where the first equality is by definition, and the second 
by construction where the almost sure convergence is subject to $\tchange =0$ in \eqref{e:QCDmodel}.

For example, in the C.I.\ i.i.d.\ model we obtain \eqref{e:LLRn}
  using $F_k =L(\Obs_k) $ with $L$ the LLR, and the Law of Large Numbers gives $\clK(\trajP_1  \|  \trajP_0)  =   D( p_1^0  \| p_0^0 )$.  More recent work establishes similar results for the  C.I.\  model in which  $\bfmX^i$, $i=0,1$ are Markov  or hidden Markov models (HMMs)  \cite{fuh03,fuhtar18},
  as well as theory and algorithms requiring minimal a~priori statistical information  \cite{liatarvee22,zhasunherzou23b}.
  See these papers and \cite{veeban14,xiezouxievee21,liavee22,tarbarfuhxin26} for surveys.    

\notes{  zhasunherzou23 is a conference paper on the C.I. Markov model.   And note that I dropped this fuh04 which considers Shiryaev–Roberts for HMMs}

In this paper we do not assume that  $\{F_k :  k\ge 1\}$ satisfy \eqref{e:LLRn}, as it requires excessive statistical information.   Moreover, the optimal choice is not the LLR even in the i.i.d.\ setting (recall discussion surrounding \eqref{e:hazardAss}).
However, the relative entropy rate is central to analysis.

The approximation of $\MDD$ in the proof of  \Cref{t:barcApprox} is based on approximation of the   \textit{overshoot}:\begin{equation}
 \oshoot_\infty \eqdef \lim_{\thresh \to\infty}  \oshoot_{\thresh} \,,
 \qquad 
  \oshoot_{\thresh} \eqdef \Expect_1[ \tstop^{\thresh} ] - \thresh/m_1 \,, 
 \label{e:oshootDef}
\end{equation} 
where  $\tstop^{\thresh} $ denotes \eqref{e:threshold} using threshold $\thresh$,
and the subscript in the expectation indicates $\tchange=0$.    
 
Approximation of $\MDE = \Expect[ (\tstop -\tchange)_-  ]$   requires new techniques since we allow for freedom in the choice of  $\{F_k\}$.   The analysis is based on sample path large deviations theory for reflected random walks,  e.g.~\cite{sie86,ganoco02,ganocowis04,dufmet05b,bazblarhezwa24}.    Relative entropy rates appear throughout the large deviations theory literature.   In particular,     on  the event that $\tstop < \tchange$, when the threshold is large we can approximate the distribution of $\{ X_k^0 : 0\le k \le \tstop \}$ by the \textit{twisted law}  $\twtrajP_0$ defined in \Cref{s:extend},
and the  error exponent appearing in the probability of this rare event is proportional to $\clK(\twtrajP_0 \| \trajP_0) $.

Recall from \Cref{t:Jentropy},
\begin{equation}
\barJFineStar(\kappa)
=
\frac{ \log\kappa+\frac{1}{2}\log\log\kappa+ \aFine  
}{ \expa
				+
	\clK(\trajP_1\|\trajP_0)
	-
	\clK(\trajP_1\|\twtrajP_0).  }  
 \label{e:Jentropy}
\end{equation}
Hence   the approximation $\barJFineStar(\kappa)$ is minimized if   the CUSUM recursion is constructed so that $\clK(\trajP_1\|\twtrajP_0) =0$,
which is achieved if $\{F_k\}$ are chosen so that \eqref{e:FuhAssumption} holds.    Maintaining only the most significant terms we obtain in this case   $\barJFineStar(\kappa) \sim \log\kappa / [\expa + \clK(\trajP_1\|\trajP_0)]$.  In the limiting case   $\expa =0$,
 this coincides with   Lorden's asymptotic for $\lorMDD $ in the C.I.\ i.i.d.\ model  \cite[Theorem 1]{lor71}.

The law $\twtrajP_0$ appeared previously in the moment-matching expression
\eqref{e:mmIntro}.
Similar conclusions are obtained in the theory of  complexity constrained approaches to universal hypothesis testing in \cite{huamey10,huaunnmeyveesur11}.

\wham{Organization}   
Details of the POMDP model are provided  in  \Cref{s:POMDP} followed by our first main result justifying the approximation of this model by a C.I.\ HMM.   \Cref{s:extend} focuses on the C.I.\ model subject to milder assumptions on $\bfmX^0,\bfmX^1$.    
Justification of  \eqref{eq:refined_bdds} is established for these models as well as the POMDP model introduced in   \Cref{s:POMDP}.  
\Cref{s:thetaCUSUM} concerns the selection of $F$ within a function class to minimize the dominant term in the approximation $ \barJFineStar ( \kappa) $.
Numerical results are surveyed in \Cref{s:numres}, followed by conclusions and future research in \Cref{s:conc}. 	

Proofs of the technical results 
and details of the numerical experiments
 are contained in the Appendix.
  
\section{POMDP Reduction to a C.I.\ Model}
\label{s:POMDP}

In this section we provide details of the POMDP model for the QCD problem, and establish  
a useful approximation by a conditionally independent (C.I.) model in the large-$\kappa$ analysis of CUSUM that is the focus of this paper.  
 
\subsection{Notation and basic assumptions}

The following are assumed throughout the paper:
 
\whamb The   limit \eqref{e:hazardAss} exists
 with $\expa\in(0,\infty)$.

\whamb
For some $\del\ge 1$ and a measurable function $F\colon\ystate^\del\to\Re$ we have $F_n = F(\strY{\del}_n)$ with $\strY{\del}_{n} = (Y_{(n-\del+1)_+},\dots,Y_{n})$
in the CUSUM recursion \eqref{e:CUSUM}.

\whamb
A real-valued random variable \(X\) is called   \emph{strongly non-lattice}   if 
\begin{equation}
\limsup_{|\theta|\to\infty}
\left|
\Expect\!\left[e^{j\theta X}\right]
\right|
<1 
\label{e:strong_non-lattice}
\end{equation}
This assumption is imposed on $X=F_n$  to obtain approximations of   $\MDD$  for a C.I.\ HMM  in \cite{fuhlai01}.

\whamb
The following limits exist:
\begin{equation}
\begin{aligned}
m_0 &= \lim_{n\to\infty}   \Expect[ F_n \mid \tchange >n]   
\\
m_1 &= \lim_{n\to\infty}   \Expect[ F_n \mid \tchange \le  n]    
\end{aligned}
 \label{e:DriftAssumptionsGen}
\end{equation}
satisfying    $m_0<0$ and $m_1>0$.

 \begin{subequations}

\whamb 
For $\del\ge 1$ and  a measurable function $G\colon\ystate^{\del} \to\Re$ 
the (conditional)
log \textit{cumulant generating functions} (CGFs) are defined by
\begin{align}
\Lambda_0(G) &= \lim_{n\to\infty} \frac{1}{n} \log \Expect\Big[ \exp\Big(  \sum_{k=0}^{n-1} G(\strY{\del}_k) \Big)
				\big|  \tchange >   n\Big]
\label{e:logCGF_POMDP0}
\\
\Lambda_1(G) &= \lim_{n\to\infty} \frac{1}{n} \log \Expect\Big[ \exp\Big(  \sum_{k=0}^{n-1} G(\strY{\del}_k) \Big)
				\big|  \tchange \le   n\Big]
\label{e:logCGF_POMDP1}
\end{align}%
It is assumed that $\Lambda_i(G)$ exists  for each $i$, with  $G =  \thexp F$ and a suitably large range of $\thexp$.
\end{subequations}

\begin{subequations}
For the C.I.\ model  these admit the alternative expressions,    
\begin{align}
\Lambda_0(G) &= \lim_{n\to\infty} \frac{1}{n} \log \Expect_0\Big[ \exp\Big(  \sum_{k=0}^{n-1} G(\strY{\del}_k) \Big)
				 \Big]
\label{e:logCGF0}
\\
\Lambda_1(G) &= \lim_{n\to\infty} \frac{1}{n} \log \Expect_1\Big[ \exp\Big(  \sum_{k=0}^{n-1} G(\strY{\del}_k) \Big)
				 \Big]
\label{e:logCGF1}
\end{align}%
where the subscripts in the expectation correspond to $   \tchange =\infty$ in \eqref{e:logCGF0} and $   \tchange = 0$ in \eqref{e:logCGF1}. 
 
\label{e:logCGF}%
\end{subequations}%

\wham{Shiryaev's model}
This  is the  C.I.\   i.i.d.\ model in which   the change time     has  a geometric distribution. It follows that \eqref{e:hazardAss} holds with $\expa=-\log\Prob\{\tchange\ge n\}/n$ for each $n\ge1$.  	 Moreover, we obtain a POMDP model with $\bfPhi =  ( \bfmZ, \bfmB)$  in which $\bfmB =
( \bfpreObs;  \bfpostObs)$  is i.i.d., and $Z_k = \ind \{  \tchange \le  k \} $ is Markovian.   
Recalling \cite{shi77} and the tutorial \cite{veeban14}, 
an optimal policy  has the following form, for some   threshold $\threshSmallStar>0$:
\begin{equation}
U_k^* =  \ind\{  p_k \ge \threshSmallStar\}  \,, \quad p_k     =  \Prob\{  \tchange \le  k \mid \Obs_0^k \}  
\label{e:SyiryaevPolicy}
\end{equation}
The same structure holds for the C.I.\ Markov model with geometric change time:
 \eqref{e:SyiryaevPolicy} holds with the constant threshold $\threshSmallStar$ replaced by $\threshSmallStar(Y_k)$  for a function $\threshSmallStar\colon\ystate\to\Re_+$
\cite{yak94}.

\subsection{Beyond Shiryaev's model}
\label{s:POMDPa}

The application of POMDP theory to QCD goes back to Shiryaev and remains an active research area.
The model considered in this section is described as follows:

\whamb  The observations take the form $\Obs_k = h(\Phi_k)$,  $k\ge 0$,  where $\bfPhi$ is a   time-homogeneous Markov chain evolving on a state space $\bigstate$,  and $h\colon\bigstate\to\ystate$ (measurable in an appropriate sense). 
\whamb   There is a decomposition $\bigstate = \bigstate_0 \cup \bigstate_1$,   for which $\bigstate_1$ is \textit{absorbing}:   $\Phi_k\in \bigstate_1$ for all $k\ge 0$ if $\Phi_0\in \bigstate_1$.     The change time is defined by $\tchange = \min\{k\ge 0 :  \Phi_k \in \bigstate_1 \}$;  consequently,  $\tchange =0$ when $\Phi_0 \in\bigstate_1$.  

Consequently, the distribution of $\tchange $ is ``phase-type''  when $\bigstate_0 $ is finite.   These  distributions are used in
areas of applied probability including 
 risk theory and 
 queueing theory \cite{blanie17}.

\whamb The input is denoted $U_k \in  \ustate = \{0,1\}$,  with  $\tstop$ defined as the first value of $k$ such that $U_k=1$.

\whamb  
Cost functions  adopted so that \eqref{e:ObjPOMDPQCD}
is consistent with \eqref{e:MDD+kappaMDE}:
\begin{equation}
\begin{aligned}
				c_\circ(\bst)   &=  \ind\{ \bst \in \bigstate_1\}    
	\\
				c_\bullet(\bst)  &= \kappa  \Expect[ \tchange \mid \Phi_0=\bst]   
\end{aligned}
\label{e:costPOMDP}
\end{equation} 
Our goal is  to minimize \eqref{e:ObjPOMDPQCD} or equivalently  \eqref{e:MDD+kappaMDE} over all inputs adapted to the observations.

This appears to be a special case of a POMDP since the observations are a \textit{deterministic} function of the Markov chain.  In most applications we have $\Phi_k = (Z_k, \Noise_k) $ in which $\bfmZ$ is a Markov chain and $\bfNoise$ is i.i.d., which captures the standard POMDP setting \cite{kri16}.

It is known that, under general conditions, the optimal solution is obtained as \textit{information-state feedback}  $U_k = \fee^*(\condDist_k)$ in which $\{\condDist_k \}$ is the conditional distribution of $\Phi_k$ given observations $Y_0^k$,  and $\fee^*$ takes values in $\{0, 1\}$. 
The policy \eqref{e:SyiryaevPolicy} obtained for 
Shiryaev's C.I.\ i.i.d.\ model is an elegant illustration, since   $ p_k =  \Expect[ Z_k   \mid \Obs_0^k ] = \condDist_k(1)$,
so that $\fee^*(\condDist_k) = \ind\{   \condDist_k(1)  \ge \threshSmallStar\} $.

		\smallskip

The main result of this section asserts that  the POMDP model very closely approximates a C.I.\ model, subject to Assumption~\Ass{s:POMDP} that follows.  Moreover, 
  the limit \eqref{e:hazardAss} exists with   $\expa>0$.

\subsection{Metastability and approximations}
\label{s:meta}
   
The following will be imposed on the POMDP model:

\wham{\Ass{s:POMDP}.}  \textit{Assumptions   for the  POMDP model}.

\whamrm{(i)}    $\bigstate$ is finite and 
		the set $\bigstate_1$ is absorbing:   $P(\bst,\bst') = 0$ for each $\bst\in\bigstate_1$ and $\bst'\in\bigstate_0$.

\whamrm{(ii)}    The Markov chain $\bfPhi$ is uni-chain and aperiodic:  there is  $\bst^\bullet \in \bigstate_1$  and  $n^\bullet\ge 1$ such that   $P^n(\bst , \bst^\bullet)>0$  for all $n \ge n^\bullet$ and $\bst\in\bigstate$.

\whamrm{(iii)}    
There is $\bst^\circ \in \bigstate_0$  and $n^\circ\ge 1$  such that 
$\Prob\{  \Phi_n = \bst^\circ    \}   >0  $ for $ n\ge n^\circ$ and  each $\Phi_0 \in\bigstate_0$.

\smallskip 

Since $\bigstate_1$ is absorbing it follows under \Ass{s:POMDP}~(i) and (ii)  that $\Prob\{  \Phi_n = \bst    \}   \to 0$ as  $n\to\infty$ for $\bst\in \bigstate_0$, but this does not rule out  \Ass{s:POMDP}~(iii).

 The assumption that $\bigstate$ is finite is far from necessary, but will simplify the proofs that follow.   The remaining assumptions are crucial to obtain the desired approximation by a C.I.\ model in which  $\{ \preObs_k , \postObs_k \}$ in    \eqref{e:QCDmodel} take the form
\begin{equation}
\preObs_k = h (\oPhi_k^\infty)   \,,  \qquad   
\postObs_k = h (\Phi_k^\infty)  \,,  \qquad  k\ge 0\,,   
\label{e:ExplicitApprox}
\end{equation}
where    $\{ \oPhi_k^\infty\} $, $\{ \Phi_k^\infty\} $ are stationary Markov chains.

 The statement of the proposition that follows requires notation from the theory of non-negative matrices:  Consider the non-negative matrix with indices in $\bigstate_0$ defined by   $\haP(\bst,\bst')= P(\bst,\bst') $ for $\bst,\bst' \in \bigstate_0$.    The assumptions imposed allow the application of Perron-Frobenius theory:  $\haP$ has  an eigenvector  $\xi$ with   entries strictly positive on $\bigstate_0$  and eigenvalue $\PFeig>0$.    The  \textit{twisted} (or \textit{tilted})  
   transition matrix   is  defined for $\bst, \bst' \in \bigstate_0$ by 
\begin{equation}
\oP(\bst,\bst') = \frac{1}{\PFeig} \frac{\xi(\bst')}{\xi(\bst)}  \haP(\bst,\bst') 
\label{e:twisted}
\end{equation}

\begin{subequations}

\begin{proposition}[Approximation by a C.I.\ model]
\label[proposition]{t:POMDPhazard}
Consider the POMDP model     subject to Assumption \Ass{s:POMDP}.   Then,  $\bfPhi$ has a unique invariant pmf  $   \metamarg$ supported on $\bigstate_1  $, and the following hold for any $m\ge 1$ and states $\{\bst_1,\dots, \bst_m \}$.
In each part the sequence $\{ \epsy_n\}$ converges to zero geometrically quickly.   It may differ in each appearance, and 
 may depend on these states, as well as $\Phi_0$.
 
	\whamrm{(i)}  					
Letting $\bfPhi^\infty$ denote a stationary version of $\bfPhi$,  
\begin{equation}
 \begin{aligned}
 \! \! \! 
 \! 
& \Prob  \{  \Phi_{n-k+1} = \bst_k \,, 1\le k\le m  \, , \,    \tchange \le   n  \}     
		\\
 \! \! \! 
 \! 
 		&  \ =    [1 +   \epsy_n  ]  \ 
		\Prob\{  \Phi^\infty_{n-k+1} = \bst_k \, , 1\le k\le m \}
			\Prob\{      \tchange \le   n  \}   
 \end{aligned}
\label{e:POMDPergodic_i}
\end{equation}
The conditional log CGF $  \Lambda_1(G) $ defined in \eqref{e:logCGF_POMDP1} exists for any function $G$  and coincides with \eqref{e:logCGF1}.

\whamrm{(ii)} 
The transition matrix \eqref{e:twisted}
defines a Markov chain $\bfoPhi$ on $\bigstate_0$ that  is uni-chain and aperiodic, with unique invariant pmf $
\ometamarg$.  
Let  $\bfoPhi^\infty$ denote a stationary version of  this Markov chain,  and suppose    $	\Prob\{  \oPhi^\infty_{n-k+1} = \bst_k \, , 1\le k\le m \} >0$.   
Then,  for any   $\ell\ge 1$,    
\begin{equation}
 \begin{aligned}
  \! \! \! 
 \! 
		\Prob\{ & \Phi_{n-\ell -k+1} = \bst_k \,, 1\le k\le m   |   \tchange >   n  \}     
		\\
 \! \! \! 
 \! 
		&\  =    [1 +   \epsy_n + \delta_\ell ]  \ 
		\Prob\{  \oPhi^\infty_{n-k+1} = \bst_k \, , 1\le k\le m \}
 \end{aligned}
\label{e:POMDPergodic_ii}
\end{equation}
where   $\delta_\ell\to 0 $ geometrically quickly as $\ell\to\infty$.

The conditional log CGF $  \Lambda_0(G) $ defined in \eqref{e:logCGF_POMDP0} exists for any function $G$, and is equivalently expressed as the standard log CGF for $\oX_n =  h (\oPhi_n) $:
 \begin{equation}
  \Lambda_0(G) 
= \lim_{n\to\infty} \frac{1}{n} \log \Expect\Big[ \exp\Big(  \sum_{k=0}^{n-1} G( \stroX{\del}_k   ) \Big)  \Big]
\label{e:logCGF0twist}
\end{equation}
with arbitrary $\oPhi_0 \in \bigstate_0$,  
where $\stroX{\del}_k =
( \oX_{(k-\del+1)_+},\dots, \oX_{k})$     for $k\ge \del$.
 
\whamrm{(iii)}  The limit  \eqref{e:hazardAss} holds with $   \expa = - \log(\PFeig)$, where $\PFeig>0$ is the eigenvalue used to define the transition matrix $\oP$ in \eqref{e:twisted}.   Moreover,  there is a constant $b_0$ such that  
\begin{equation}
				\Prob\{    \tchange >  n  \}     =     [  b_0 +  \epsy_n   ] \exp( - \expa n )  
\label{e:POMDPapprox_tchange}	
\end{equation}
\clThm			
\end{proposition}
			
\label{e:POMDPapproxInd}
\end{subequations}

%
%
%

\section{Performance for C.I.\ Models}
\label{s:extend}

Having justified the C.I. approximation in \Cref{s:POMDP}, we now state conditions directly on a general C.I. model under which the    approximations    \eqref{e:barFineBody} hold for CUSUM.

Recall from the discussion surrounding \eqref{e:Best_H_J} that   $\barJ(\thresh, \kappa)  $  denotes    \eqref{e:MDD+kappaMDE} using CUSUM with threshold $\thresh>0$,   and  $
\barthresh^*(\kappa)  =\argmin_{\thresh} \barJ(\thresh, \kappa)  $.   
And CUSUM*   denotes the algorithm using the optimal threshold $\barthresh^*(\kappa)  $.

The cost $\barJ( \thresh, \kappa)  $ consists of two components which require separate analysis:  
\whamb
Approximation of $\MDD$ requires the large threshold approximation:  
with $\oshoot_\infty$ defined in  \eqref{e:oshootDef},
\begin{equation}
\MDD =    \thresh/m_1 + \oshoot_\infty  + o(1)  
\label{e:MDD_renewal}
\end{equation}

\whamb  A simple proof of the following representation may be found in \Cref{t:MDE1}:
\begin{equation}
\MDE   =   
        \sum_{n=0}^{\infty}
        \Prob\{\tstop\le n,\ \tchange>n\}.
\label{e:MDEstep0}
\end{equation} 
The terms in the sum  may be approximated through application of the theory of rare events for reflected random walks (e.g.,   \cite{ganoco02,ganocowis04}).

\subsection{General C.I.\ model and assumptions}

Recall that the C.I.\ model refers to   \eqref{e:QCDmodel}
in which  $\bfpreObs$, $\bfpostObs$ and $\tchange$ are statistically independent.   This is assumed throughout this section, and in addition   the observation space $\ystate$ is Euclidean space.

  Then,   approximation of $\MDD$ is   based on    $\{    \postObs_k \}$ while approximation of $\MDE$ is   based on   $\{  \preObs_k \}$.   Analysis reduces to consideration of the partial sums
 \[
        S_n^i=\sum_{k=1}^{n} F^i_k,\qquad  n\ge 1 \,, \ i=0,1\, 
\]
where $\{F_k^i = F( X^i_k)\}$ is a stationary sequence,  and   $F\colon\ystate\to\Re$ is Borel measurable.


Approximation of  $\MDE $  is based on large deviations theory. 
We  perform a temporal and spatial scaling that is standard in this literature: 
for a given threshold $\thresh>0$    consider the continuous functions defined for $k\ge 1$ by
\begin{equation}	
 \xH_{t} =  \thresh^{-1} \InfoState_k \,,  \ \   \sH_{t} =  \thresh^{-1} S_k^0 \,, \quad \textit{	 for $t = k/\thresh$, }
\label{e:xHxS}
\end{equation} 
with $ \xH_0=\sH_0=0$, 
 and defined for $t>0$  
 by piecewise linear interpolation.
 	
When $\tstop$ is defined using   threshold $\thresh$, then  $\tstop\le  T\thresh $  if and only if $    \xH_t    \ge 1$ for some $t\le T$.  
On denoting     $s_n = n/\thresh$, the identity  \eqref{e:MDEstep0} becomes 
 \[
\MDE 
		=    \sum_{n=0}^{ \infty }    \Prob\{  \max_{0\le t\le s_n}  \xH_t \ge 1 \,,  \tchange  >  n  \} 
 \]
 and hence for the C.I.\ model, 
\begin{equation}
\MDE 
		=    \sum_{n=0}^{ \infty }    \Prob_0\{  \max_{0\le t\le s_n}  \xH_t \ge 1 \}    \Prob\{ \tchange  >  n  \} 
\label{e:MDEstep1}
\end{equation} 
where the subscript in $ \Prob_0$  indicates that $Y_k = X_k^0$ for all $k$ in this probability.   

LDP theory provides estimates of  
\[
\frac{1}{\thresh} \log
   \Prob_0\{  \max_{0\le t\le s_n}  \xH_t \ge 1 \}  \,, \quad   \frac{1}{\thresh} \log \Prob\{  \max_{0\le t\le s_n}  \sH_t \ge 1 \} \,, 
\]
which coincide in the large $\thresh$ limit.    These estimates   require the cumulative log moment generating function  and its dual:  for $ \thexp\,,\, m\in\Re$,
\begin{equation}
\Upupsilon_0(\thexp) \eqdef \Lambda_0(\thexp F)   \,, 
\quad
			I_0(m) \eqdef \sup_\thexp [  \thexp m - \Upupsilon_0(\thexp   )  ]   
\label{e:I0}
\end{equation}
where $\Lambda_0$ is defined in \eqref{e:logCGF0}.    
 We then define,
\begin{equation}
  e_0(t)  = \begin{cases}
				t I_0(1/t)       &   t <  1/\cm_0  
	\\
				I_0(\cm_0)/ \cm_0  &   t \ge  1/\cm_0  
\end{cases} 
\label{e:rate1}
\end{equation}
Let $ \Fmax^0 $ denote the essential supremum $[\marg_0]$  of the random variable $F(X_k^0)$ (independent of $k$),   and
$\Tmin \eqdef 1/( \Fmax^0)$.
It is shown in \Cref{t:e0infinite} that  $e_0(t) =  \infty$ for   $0<t < \Tmin $.

\wham{\Ass{s:extend}.}  \textit{Assumptions   for the C.I.\   model}.

\whamrm{(i)}  
$\Fmax^0> \expa > 0$,
$m_0<0$,  and 
$m_1>0$.

\whamrm{(ii)} 
The limit following 
exists:
\begin{equation}
 \oshoot_\infty \eqdef \lim_{\thresh \to\infty} \big [ \Expect[\tau^1_{\thresh} ] - \thresh/m_1 \big]
\label{e:oshootA4}
\end{equation} 
where $
        \tau^1_{\thresh}=\inf\{n\ge 1:S^1_n \ge \thresh\}$.

\whamrm{(iii)} 
$\Upupsilon_0 (\thexp)$ is defined for each $\thexp$ to define a convex extended valued function.    
There are solutions   $ \thexp_+ > \thexp_0>0$ to 
\begin{equation}
\Upupsilon_0(\thexp_0) = 0 \,, \ \ 
\Upupsilon_0(\thexp_+) = \expa
\label{e:thexp+}
\end{equation}
and $\Upupsilon_0$ is twice   continuously differentiable on an open interval $\clI$ containing $[0,\thexp_+]$.
Moreover, 
  $\cm_+ > \cm_0 >0$, where 
\begin{equation}
\cm_0 = \Upupsilon_0'(\thexp_0)\,, \ \ \cm_+ = \Upupsilon_0'(\thexp_+). 
\label{e:cm0cm+}
\end{equation}

\whamrm{(iv)} 
The partial-sum process   $\{S_n^0\}$ satisfies the following sample-path large deviations principle on each compact time
interval:   
For every $T>0$, the law of $\{\sH_t:0\leq t\leq T\}$ satisfies 
an LDP with speed $\thresh$
and good rate function
\[
          I_{[0,T]}(\phi)
        =
        \begin{cases}
        \displaystyle
        \int_0^T I_0\bigl( \ddt \phi(t)\bigr)\,dt,
        &
        \phi  \in  \text{AC}[0,T] \,, \  \phi(0)=0 
        \\[2ex]
        +\infty,
        & \text{otherwise}.
        \end{cases}
\]
where $    \phi  \in \text{AC}[0,T]$ indicates that $\phi$ is absolutely continuous on this interval.
 \qed
 
 \smallskip

Observe that \eqref{e:oshootA4} is equivalent to   \eqref{e:oshootDef} for the C.I.\ model.   
It is shown in  \cite{fuhlai01} that this holds for a class of  C.I.\ HMMs.   Among the assumptions is the strong non-lattice condition imposed on the common distribution of $\{F_k^1 : k\ge 1 \}$.

 
Conditions for   (iii) and (iv) are known if $\bfmX^0$ is  Markovian  or  a  hidden Markov model   \cite{konmey03a,konmey05a,ganoco02}.

 \subsection{Performance approximations}
 \label{s:perCUSUM}

\begin{proposition}[Performance and threshold approximations]
\label[proposition]{t:barcApprox}
Subject to Assumptions \Ass{s:extend} (i)--(iv)  the approximations \eqref{eq:refined_bdds} hold with 
\begin{equation}
 \begin{aligned} 
 \bFine = \log\bigl(m_1\sqrt{2\pi\gamma^2\thexp_+}\bigr) \,,
 \quad
\gamma^2=\Upupsilon_0''(\thexp_+)/\cm_+^3
 \end{aligned}
\label{e:RefinedPars}
\end{equation}
and $ \aFine
= 1 + \bFine+ \oshoot_\infty m_1 \thexp_+  $.
\clThm
\end{proposition}

These conclusions hold for the POMDP model under assumptions compatible with \Ass{s:extend}.
Recall the pmf $\ometamarg$ is defined in \Cref{t:POMDPhazard}.

\begin{corollary}
\label[corollary]{t:barcApproxPOMDP}
Consider the POMDP model     subject to  
Assumption~\Ass{s:POMDP} and the   following variants of parts (i) and (ii) of Assumption~\Ass{s:extend}:

 \whamrm{(i)} 
$     \max\{F(h(x)):\ometamarg(x)>0\}>\expa$, and with
\[
 m_0= \sum_x F(h( x)) \ometamarg (x) \qquad
	m_1= \sum_x F(h( x)) \metamarg (x) 
\] 
we have $m_0<0$ and $m_1>0$.

\whamrm{(ii)}    For any initial condition $\bst \in\bigstate$, the following limit exists, possibly depending on $\bst$:  
\[
 \oshoot_\infty \eqdef \lim_{\thresh \to\infty} \big [ \Expect[\tstop(\thresh)  ] - \thresh/m_1 \big]
\]
where $\tstop(\thresh) $ is defined in \eqref{e:threshold}.

Then, the conclusions of \Cref{t:barcApprox} hold,
in which $\Lambda_0$ is defined in   \eqref{e:logCGF0twist},  and  
$\Upupsilon_0(\thexp) $, $  I_0(m)$ in 
\eqref{e:I0}
are defined with respect to this log~CGF to obtain the parameters \eqref{e:RefinedPars}.
\clThm
\end{corollary}

\Cref{t:barcApprox} is established in \Cref{s:BigProof}
of the Appendix.

The proof of the corollary is omitted as it is identical to the proof of \Cref{t:barcApprox}.
However,  we point out here the key steps of the proof in the POMDP setting:  
The expression \eqref{e:MDEstep0} motivates consideration of the limit,
 \begin{equation}
G(t) = - \lim_{\thresh\to\infty}    \tfrac{1}{\thresh} \log         \Prob\{\tstop\le \thresh t ,\ \tchange>  \thresh t \} \,, \quad t >0 \,, 
\label{e:MDEstep0.1}
\end{equation}  
with $\tstop$ within the probability   defined using the threshold $\thresh$.   The limit  may be expressed,  
\begin{equation}
 G(t) =    e_0(t)     +  \expa t     
\label{e:rate1G}
\end{equation}
where  $e_0(t)   $ is defined in \eqref{e:rate1}, recalling that for the POMDP model $I_0$ is the convex dual of $\Lambda_0$   defined in   \eqref{e:logCGF0twist}.

Subject to  \eqref{e:hazardAss}, approximation of the probability on the right hand side of \eqref{e:MDEstep0.1} reduces to approximating the conditional probability,
\begin{equation}
    \Prob\{\tstop\le \thresh t \mid \tchange>  \thresh t \}  =  
    \Prob\Big\{ \max_{0 \le r\le  t} \xH_r \ge 1 \mid \tchange>  \thresh t \Big \} 
\label{e:cond_stop_rare}
\end{equation}
The sum \eqref{e:MDEstep0} is then approximated by an integral of $\exp(-\thresh G(t))$ over a small time-interval, whose width vanishes as $\thresh\to\infty$.   
This justifies a second-order Taylor series approximation of $G$ near its minimizer $s^*$.   It is shown in  \Cref{t:Gprops} that $s^* = 1/\cm_+$   and  $G''\, (s^*)  = 1/ \gamma^2 $.

These arguments provide an approximation $\MDE(\thresh) $ of $\MDE(\thresh)$ for large positive thresholds.

  It is perhaps surprising that logarithmic asymptotics such as
\eqref{e:hazardAss} ultimately lead to the refined approximations in  \Cref{t:barcApprox} with additive $o(1)$ error.  The key observation in the proof   
  is that the
weighted eagerness term  $ \kappa \MDE(	\barthreshFineStar(\kappa) ) $  is bounded as a function of $\kappa$.
   Consequently,
\begin{equation} 
 \kappa \MDE(	 \barthreshFineStar(\kappa)  )  [1+o(1)]
				 = \kappa \MDE(	\barthreshFineStar(\kappa)  )  + o(1)
 \label{e:whyTight} 
\end{equation} 
This converts the relative error in the approximation of
$\MDE$ into the additive $o(1)$ error appearing in    \eqref{eq:refined_bdds}.

 \subsection{Twisted laws and information-theoretic interpretation}
 \label{s:InfoInterpretations}
 
 Here we explain how the performance bound can be expressed in terms of relative entropy rates to
obtain  \eqref{e:Jentropy}  and  the connection with Lorden's asymptotic for $\lorMDD $ discussed in the Introduction.
This requires an explanation of the twisted law $\twtrajP_0$ appearing in  \eqref{e:Jentropy}.

 For each $n$ and  $\thexp$
denote by $\twtrajP_0^{(n,\thexp)}$ the probability measure on $\clB(\ystate^{n+1})$  defined for $Z\in \clB(\ystate^{n+1})$ by
\begin{equation}
\begin{aligned}
\twtrajP_0^{(n,\thexp)}\{Z\} =\beta_{(n,\thexp)}^{-1}\Expect\big[\exp(\thexp S_n^0)\,\ind_Z(X_0^0,\dots,X_n^0)\big] \,, 
\end{aligned}
\label{e:twtrajP+}
\end{equation}
  $\beta_{(n,\thexp)}=\Expect[\exp(\thexp S_n^0)]$,
whose  last marginal is denoted  
\[
\twmarg_0^{(n,\thexp)} \{ A \}  =  \beta_{(n,\thexp)}^{-1}     \Expect\big[ \exp( \thexp S_n^0)   \ind_A(  X_n^0 ) \big]  \,,\quad A\in\by
\]
For the special case $\thexp = 0$ these probability measures are denoted $\trajP_0^{(n)}  $,   $\marg_0^{(n)} $ respectively;  	stationarity implies that $\marg_0^{(n)} = \marg_0$ for each $n$.    We also reserve the notation $\trajP_1^{(n,\thexp)}  $ when $S_n^1$ is used in \eqref{e:twtrajP+}.

   It will be assumed that the following construction of the    \textit{twisted process} $\{ X^\thexp_k : k\ge 0 \}$ is well defined for $\thexp \in \clI$ (recall that the interval $\clI$ is introduced in \Ass{s:extend}).   It 
is the stationary stochastic process with law $\twtrajP_0^{\thexp}$,  whose
 finite dimensional marginals are defined as follows:   for $m\ge0$ and $ Z\in \clB(\ystate^{m})$,
\begin{equation}
\begin{aligned}
\Prob &\{ (X^\thexp_{k+1}, \dots, X^\thexp_{k+m} ) \in Z \} 
\\
&\eqdef  
	\lim_{n\to\infty}  \beta_{(n,\thexp)}^{-1}\Expect\big[\exp(\thexp S_n^0)\,\ind_Z(X_{n-m+1}^0,\dots,X_n^0)\big] 
 \, .
 \end{aligned}
	\label{e:twtrajP}
\end{equation}
The marginal $	\twmarg_0^\thexp\{ \varble  \} = \Prob\{ X^\thexp_k \in \varble \} $    coincides with the limit, 
\begin{equation}
\twmarg_0^\thexp\{ A \} 
= \lim_{n\to\infty} \twmarg_0^{(n,\thexp)} \{ A \}   \,,\quad A\in\by
	\label{e:twmarg0}
\end{equation}
We simplify notation to  $\twtrajP_0$  and $\twmarg_0$  when $\thexp =1$.

\wham{\Ass{s:extend} (v)}  \textit{Twisted process   for the C.I.\   model}.
 
\whamb
For each $\thexp \in \clI$ the    \textit{twisted law} $\twtrajP_0^{\thexp}$ defining the twisted process exists, through the construction above.

 \whamb
\textit{Relative entropy rate:}    The following   exists as a finite limit for $\thexp \in \clI$:
\begin{equation}
	\clK(\twtrajP_0^{\thexp} \|  \trajP_0)  = \lim_{n\to\infty} \frac{1}{n}   D( \twtrajP_0^{(n,\thexp)}  \|  \trajP_0^{(n)}   ) 
	\label{e:RelEntRate}
\end{equation}
where $D$ denotes relative entropy.        It is also assumed that $\clK(\trajP_1  \|  \trajP_0) <\infty$ where   
\begin{equation}
	\clK(\trajP_1  \|  \trajP_0)  = \lim_{n\to\infty} \frac{1}{n}   D( \trajP_1^{(n)}  \|  \trajP_0^{(n)}   ) 
	\label{e:RelEntRateModel}
\end{equation}

 \whamb
\textit{Entropic calculus:}  For $\thexp \in \clI$,
\begin{equation}
\frac{d}{d \thexp }  \Upupsilon_0(\thexp)    =    \twmarg_0^{\thexp} (F)  
	\label{e:twistedMean}
\end{equation}
Consequently,  letting $  F^0_{\textsf{min} } ,   F^0_{\textsf{max} } $ denote the respective   essential-infimum, essential-supremum   of   $F(X_k^0)$ w.r.t.\ $[\marg_0]$,
\begin{equation}
  F^0_{\textsf{min} } \le 
   \frac{d}{d \thexp } 
\Upupsilon_0(\thexp)    \le  F^0_{\textsf{max} } 
\label{e:derBdds}
\end{equation}

Conditions for   \Ass{s:extend}~(v) are easily established in the i.i.d., Markovian and even for hidden Markov models   \cite{konmey03a,konmey05a}.  
In each of these three model classes, the class of the law $\twtrajP_0$  is unchanged.  In particular, 
for the C.I.\ HMM model described in \Cref{s:trad},
under mild assumptions we find that   $\twtrajP_0$ is itself defined by an HMM  in which the transition law for $\bfPhi^0$ is modified and $h$ is unchanged---an example is provided in \Cref{s:Model2Relax}.


\begin{proposition}[Information-theoretic approximation]
\label[proposition]{t:Jentropy}
Under \Ass{s:extend} (v) along with 
 the assumptions of either \Cref{t:barcApprox} or \Cref{t:barcApproxPOMDP} the representation    \eqref{e:Jentropy}
holds
and
the approximation \eqref{e:barJFineBody} admits the lower bound,
\begin{equation}
\barJFineStar(\kappa)
			\ge 
\Bigl[\log\kappa+\frac{1}{2}\log\log\kappa+ \aFine\Bigr] 
\frac{1}{[\expa
				+
	\clK(\trajP_1\|\trajP_0)]}
\label{e:barJFineBodyBdd}
\end{equation}
which is achieved if and only if $\clK(\trajP_1\| \twtrajP_0^{\thexp_+})  = 0$. 		
\clThm
\end{proposition}
  
  The proof is immediate from the following:

 \begin{lemma}
\label[lemma]{t:m1_ent_F}
If $F$ satisfies \Ass{s:extend} then for $\thexp \in \clI$, 
\begin{equation}
	 \marg_1(F) \thexp
				=
	\Upupsilon_0(\thexp)
				+
	\clK(\trajP_1\|\trajP_0)
				-
	\clK(\trajP_1\| \twtrajP_0^{\thexp})  
	\label{e:m1_ent_F}
\end{equation}
Consequently,  with   $\thexp_+$ defined in \eqref{e:thexp+}  we have 
\begin{equation}
 \begin{aligned} 
m_1\thexp_+
				&=
\expa
				+
	\clK(\trajP_1\|\trajP_0)
				-
	\clK(\trajP_1\|\twtrajP_0^{\thexp_+}) 
	\\
	& \le \expa
				+
	\clK(\trajP_1\|\trajP_0) 
 \end{aligned} 
	\label{e:m1_ent_F+}
\end{equation}

\end{lemma}

\PF    
From \eqref{e:twtrajP+} we have  for any  $\thexp \in \clI$ and $n\ge 1$,
and sequence $(y_0,\dots,y_{n-1})$,
\[
\log\Bigl(\frac{d\twtrajP_0^{(n,\thexp)}}{d\trajP_0^{(n)} }\Bigr)
		=
		\thexp s_n-\log(\beta_{(n,\thexp)}),
\quad
		s_n=\sum_{k=0}^{n-1}F(y_k).
\]
Taking the expectation of each side with respect to $\trajP_1^{(n)} $ and dividing by $n$ gives
\[
		 \tfrac{1}{n} \big[ 		D(\trajP_1^{(n)}\|\trajP_0^{(n)} ) 
		 -
		 D(\trajP_1^{(n)}\|\twtrajP_0^{(n,\thexp)})  \big]
				=
		    \thexp  \marg_1(F)
		- 		 \tfrac{1}{n} \log(\beta_{(n,\thexp)}).
\]
Letting $n\to\infty$ we obtain  $\frac1n\log(\beta_{(n,\thexp)})\to \Upupsilon_0(\thexp )$.
    \qed

\section{Optimizing CUSUM}
\label{s:thetaCUSUM}

Here we consider the selection of $\del$ and  $F\colon\ystate^\del\to\Re$ in CUSUM.     To simplify exposition we restrict to the conditionally independent model \eqref{e:QCDmodel} in which $\bfpreObs$ and $\bfpostObs$ are mutually independent stationary stochastic processes, and independent of the change time $\tchange$.  We denote by $\marg_0$ and $\marg_1$ the corresponding stationary marginals on $\by$.  

For much of the theory surveyed in this section it is assumed that $\del=1$ so that  $F_n = F(\Obs_n)$, for a  function $F\colon\ystate\to\Re$.     This is without loss of generality since the specification of the observation process is one choice in algorithm design.  In particular, for Markovian models we require $\del =2$, but we may apply the theory for $\del =1$ by treating  $ \{ \strY{2}_{n} = (\Obs_{n-1}, \Obs_n ) : n\ge 0 \}$ as the observation process.

\subsection{Assumptions and notation}
\label{s:optGenAss}

Our goal is to minimize over some function class $\clG$ the approximate CUSUM cost \eqref{e:barJ_infty_F}.
We do not have theory to optimize  the finer approximation $\barJFineStar$ defined in  \eqref{e:barJFineBody}, which explains why the term $\aFine$ has been dropped in \eqref{e:barJ_infty_F}.

The minimization of $
\barJStarCoarse( \kappa ; F)  $ is subject to the constraint that  $\marg_0(F) \le 0$.    
Under mild conditions we find that this inequality is strict for an optimizer, and also   $\marg_1(F) > 0$.

\wham{\Ass{s:thetaCUSUM}.}  \textit{Assumptions   for the linear function class}.

Each $F\in \clG$ satisfies the following:   
  
\whamrm{(i)} 
  The mean $ \marg_i(F)$  exists for each $i$,

  \whamrm{(ii)} 
  If  $\Fmax^0> 0$ and 
  $ \marg_0(F) \le  0$    then parts (iii)--(v) of  Assumption~\Ass{s:extend} hold for this $F$;
  we set $\thexp_0 = 0$ if  $ \marg_0(F) =  0$.    In particular, 
  there is a unique value $\thexp_+ = \thexp_+(F)>0$   satisfying $\Lambda_0(\thexp_+ F) = \expa$.

  \whamrm{(iii)} 
    If  $ \marg_1(F) > 0$  then   \eqref{e:oshootA4} holds for some $ \oshoot_\infty \in\Re$.

  \whamrm{(iv)} 
The function $F\equiv 1$ lies in $\clG$, and there exists $\Fmix\in\clG$ satisfying $ \marg_0(\Fmix) <0$ and $ \marg_1(\Fmix)  >0  $. 
\qed

For $F\in\clG$ let $\twtrajP_0^F$ denote the law 
$\twtrajP_0^{\thexp}$ defined in \eqref{e:twtrajP} with $\thexp=1$.
Its marginal $\twmarg_0^F$  is defined using \eqref{e:twmarg0}, also with $\thexp =1$.

Whenever the assumptions of \Cref{t:m1_ent_F} hold we have,
 \begin{equation}
  \Lambda_0(F)-\marg_1(F) =  \clK(\trajP_1\|\twtrajP_0^F) - \clK(\trajP_1\|\trajP_0) 
\label{e:GammaK}
\end{equation}
which is part of the derivation of   \Ass{s:thetaCUSUM}:
\begin{equation}
\barJStarCoarse( \kappa ; F)  =
\frac{ \log\kappa+\frac{1}{2}\log\log\kappa   
}{ \expa
				+
	\clK(\trajP_1\|\trajP_0)
	 -
	\clK(\trajP_1\|\twtrajP_0^F) }  
 \label{e:JentropyF}
\end{equation}
Hence our objective in this section is to minimize $	\clK(\trajP_1\|\twtrajP_0^F) $.  

We adopt the following notation, provided the minimum is attained:  
\begin{align} 
F^\circ & \in \argmin\{ 
\clK(\trajP_1\| \twtrajP_0^F)  : F\in\clG \}
\label{e:Fcirc}
\\
\expa^0 & = \clK(\trajP_0\|\twtrajP_0^F)  \big|_{F = F^\circ}
\label{e:expa}
\end{align}

\subsection{Convexity and duality}

Minimization of $\barJStarCoarse$ may be posed as a convex program:

\begin{proposition}[Convex primal formulation]
\label[proposition]{t:CoarseConvex} 
Minimization of \eqref{e:barJ_infty_F} may be cast as the convex program under Assumption~\Ass{s:thetaCUSUM}:
\begin{equation}
 \begin{aligned}
\max \ \ &    \marg_1(F) 
\\
\textrm{s.t.} \ \ 
& \Lambda_0(F) \le \expa \, ,   \  \marg_0(F) \le 0  \, ,   \   F\in\clG 
\end{aligned}
 \label{e:jStarConvexRelaxation}
 \end{equation}
If a solution $F^*$ exists, then $ \Lambda_0(F^*) = \expa$ and    $\marg_1(F^* ) >0$. 
\clThm
  \end{proposition}

Under Assumption~\Ass{s:thetaCUSUM}~(iv) the function $G = r\Fmix \in\clG$ is strictly feasible for sufficiently small $r>0$.   That is, $\Lambda_0(G)  < \expa $ and $\marg_0(G)<0$, so that   Slater's condition holds.

The first step in the proof is scale invariance of the objective:
 \begin{lemma}
\label[lemma]{t:Jscalings}
The objective \eqref{e:barJ_infty_F} is invariant under positive scalings:
 $\barJStarCoarse( \kappa ; r F) = \barJStarCoarse( \kappa ; F) $ for $r>0$.
 \clThm
\end{lemma}

Consequently, we may assume without loss of generality that $ \thexp_+(F) =1$,
 which is equivalently expressed as the (non-convex) equality constraint 
  $ \Lambda_0(F^*) = \expa$.   Justification of the convex relaxation $ \Lambda_0(F) \le \expa $ and other details may be found in the Appendix.

We now turn to the information-theoretic representation of  $\barJStarCoarse( \kappa ; F) $ shown in  \eqref{e:JentropyF}.   
While the function $F^\circ \in\clG$ is assumed to be a minimizer,   it may not solve \eqref{e:jStarConvexRelaxation} since  the law 
$\twtrajP_0^F$ does not determine the function $F$.   
It is shown next that an optimal solution is obtained by modification through an additive constant.

\begin{proposition}[Strict sign constraint satisfaction]
\label[proposition]{t:jStarConvexRelaxationENT}
Suppose that   Assumption~\Ass{s:thetaCUSUM}   holds and that  the argmin in \eqref{e:Fcirc} is achieved to define 
   $F^\circ  \in\clG$.   It may be assumed without loss of generality that $\Lambda_0(F^\circ) = 0$.    Then, 

\whamrm{(i)}
$\marg_1(F^\circ) >0$,  $\marg_0(F^\circ) = -\expa^0 < 0  $. 

\whamrm{(ii)}
For any  $\expa \in [0, \expa^0]$  the function $F^* = F^\circ + \expa$  is a solution to \eqref{e:jStarConvexRelaxation} satisfying
   $ \Lambda_0(F^*) = \expa$,
   $\marg_1(F^*) > 0  $, 
 and  $\marg_0(F^*)  = \expa - \expa^0 \le 0  $.      \notes{!!! Error here -- no reason for equality.
\\
Also, go over proofs to see if I should use this:  
$ \expa^0 \eqdef \clK(\trajP_0\|\twtrajP_0^{F^\circ}) $.}
 \clThm
    \end{proposition}

We   say that  $\clG$  is  \textit{dense} if the following holds for any two probability measures $\upmu,\uppi$: if 
the equality $\upmu(F) =  \uppi(F) $ holds for all $F\in\clG$ then  $\upmu = \uppi$.

\begin{proposition}[Optimality over a dense function class]  
\label[proposition]{t:completeOpt} 
Suppose that   Assumption~\Ass{s:thetaCUSUM}   holds and  that   $\clG$ is   dense.     
Then $ \twmarg_0^{F} = \marg_1$ when $F=F^\circ$, 
and hence the optimizer    $F^* = F^\circ + \expa$  satisfies $ \twmarg_0^{F^*} = \marg_1$  for any  $\expa \in [0, \expa^0]$.   
\qed
\end{proposition}

The conclusion
$\clK(\trajP_1\|\twtrajP_0^F)  = 0 $ may fail to hold when $F=F^\circ$ even  when the function class is dense.
However, the proposition implies an analogous conclusion for the marginals,   $D( \marg_1 \|  \twmarg_0^{F^*} ) =0$.

\subsection{Traditional settings}
\label{s:trad}

The solution to the convex program in \Cref{t:CoarseConvex}  admits an explicit solution in the   most common settings.   
Each is a C.I.\ model with $\ystate$   Euclidean space for which we can characterize the best choice of $F$.

\wham{1. Conditionally independent i.i.d.\ model:}   
This is the classical setting in which $\bfpreObs$ and $\bfpostObs$ are i.i.d.\ (and mutually independent).
	Consequently, $\{\InfoState_n\}$ evolves as the workload in a GI/G/1 queue, with a change in load parameter at
	time $\tchange$ (e.g.\ \cite{ana89,ganoco02}).

Assuming mutually absolutely continuous marginals, the LLR is denoted $  L=\log(d\marg_1/d\marg_0)$.  
 The function $F=L$  satisfies the sign conventions,
 \begin{equation}
\marg_0(L)=-D(\marg_0\|\marg_1)<0,\quad
\marg_1(L)= D(\marg_1\|\marg_0)>0, 
\label{e:CI_Lsigns}
\end{equation}
where $D$ denotes relative entropy.

\wham{2. Conditionally independent  Markov model:}       the C.I.\ model with $\bfmX^i$  a time-homogeneous Markov chain  for each $i$.

%
%

 Much of the theory in the QCD  literature is subject to the following assumption on the respective transition kernels:    
 for a probability measure $\Pdens$ on $ \by$ and measurable functions
	$g_i\colon\ystate\times\ystate\to(0,\infty)$,
\begin{equation}
	P_i(x,dx')=g_i(x,x')\,\Pdens(dx'),\qquad i\in\{0,1\}.  
\label{e:macP}
\end{equation}
The best choice of $F$ requires $\del = 2$: 
rather than the LLR, consider the function on $\ystate^2$ defined by 
\begin{equation}
		L_\infty(x,x')\eqdef \log\Bigl(\frac{g_1(x,x')}{g_0(x,x')}\Bigr)\,,  \quad x,x'\in \ystate 
\label{e:LLR_Markov}
\end{equation}
Using   $F=L_\infty$
the sign conventions hold:
\begin{equation}
\begin{aligned}
m_0 &=    \Expect[F( X^0_{n+1},X^0_{n+2}) ]=-K(P_0\|P_1)<0
\\
m_1 & =  \Expect[F( X^1_{n+1},X^1_{n+2}) ]=K(P_1\|P_0)>0 \,,
\end{aligned}
\label{e:MarkovA1}
\end{equation}
	where $K$ denotes the Donsker--Varadhan relative-entropy rate for stationary Markov chains.  In particular,
\begin{equation}
		K(P_1\|P_0)=\dint L_\infty(x,x')\,\marg_1(dx)\,P_1(x,dx').
\label{e:DVrate}
\end{equation}
Observe that for $F=L_\infty$ we have  $\Lambda_0(F) =0$   and  $\clK(\twtrajP_0^\thexp  \|  \trajP_0) =   K( P_1 \| P_0) $ when $\thexp = 1$ (recall \eqref{e:RelEntRate}).

In the C.I.\ i.i.d.\ model the choice $F=L$ is approximately optimal  under the performance criterion \eqref{e:MDD+kappa_pFA}
(as well as others),  and in the C.I.\ Markov model analogous conclusions hold for $F=L_\infty$ 
 \cite{xiezouxievee21}.   For the performance criterion considered in this paper, the optimal function is modified by an additive constant:

\begin{proposition}
\label[proposition]{t:bestF} 
Consider the C.I.\ Markov model.   Choose $\del =2$ and the function class $\clG$ consisting  of all measurable functions satisfying Assumption~\Ass{s:thetaCUSUM}.   
	 \notes{careful here since $\del=1$ in those assumptions}
  Assume moreover that  $0\le \expa \le K(P_0\|P_1) $.   Then,
an optimizer   of    \eqref{e:jStarConvexRelaxation}  is   $\Fstar=L_\infty+\expa$, which satisfies $m_0 =  \expa  -K(P_0\|P_1)\le 
0$ and $\Lambda_0(\Fstar) =  \expa$.
\end{proposition}

\PF 
To begin, recall that on invoking  Assumption~\Ass{s:thetaCUSUM} with $\del =2$ we regard    $ \{ \strY{2}_{n} = (\Obs_{n-1}, \Obs_n ) : n\ge 0 \}$ as the observation process.

The assertion $\Lambda_0(\Fstar) =  \expa$ follows from  $\Lambda_0(L_\infty) = 0$.   
The remaining conclusions follow from  \Cref{t:Jentropy} since $\twtrajP_0^F=\trajP_1$ 
and hence $\clK(\trajP_1\|\twtrajP_0^F) =0 $ for $F=\Fstar$.   
\qed

\wham{3.  Conditionally independent HMM.} 

More recent in the literature is theory for the C.I.\  HMM.  Consider the model \eqref{e:QCDmodel} in which 
$X_k^i = h(\Phi_k^i)$ 
for a measurable function
 $h\colon\bigstate\to\ystate$,  
 and  $\{ \Phi_k^i \}$ is Markovian with transition kernel $P_i$.  

 Given  $F\colon\ystate\to\Re$ we may write  $  F(\preObs_{n+1})  =   G(\Phi_{n+1})$ using $G= F\circ h$.     
Hence $\Lambda_0( \thexp F) = \log(\lambda)$ with $\lambda>0$  the PF eigenvalue for  the positive kernel   $\haP (x, dx') =  \exp( \thexp G(\bst )) P_0(\bst, d\bst') $,
and the \textit{twisted law} is again Markovian:  
\begin{equation}
\cP(\bst, d\bst') =   
 \frac{1}{\lambda} \frac{\xi(\bst')}{\xi(\bst)}    \exp( \thexp G(\bst) ) P_0(\bst, d\bst') 
\label{e:twistedF_HMM}
\end{equation}

To understand the optimal choice of $\{ F_{n+1} \}$  for CUSUM,  for simplicity we restrict to a finite state space.   
For any value of $\del\ge 1$ we may optimize over the set of all functions  $F\colon \ystate^\del \to\Re$.   However, theory from \cite{fuh03} suggests that in general we require $\del = \infty$.  

Denote for $i=0,1$ the two random pmfs,
\begin{equation}
p^i_{n+1\mid n}(y  ) = \Prob\{Y_{n+1} = y   \mid Y_0^n \} \,, \quad y\in\ystate\,,
\label{e:CPY}
\end{equation} 
where    the conditional probability is obtained based on the state transition matrix for $\Phi_k^i$ and the function $h$.

  Combining Theorems~5--7 of \cite{fuh03} establishes a form of approximate optimality of CUSUM using
 \begin{equation}
F_{n+1} = \log \ell _{n+1\mid n} \,, \qquad \ell _{n+1\mid n} \eqdef     \frac{p^1_{n+1\mid n}(Y_{n+1}   ) }{  p^0_{n+1\mid n}(Y_{n+1}     )  }
\label{e:HMMopt}
\end{equation}
Note   that $\{ F_k \}$ defines the increments of a super-martingale when $\tchange = \infty$:  by Jensen's inequality,
\begin{equation}
\Expect[F_{n+1} \mid Y_0^n ]   \le \log   \Expect[  \ell _{n+1\mid n} \mid Y_0^n ]     =0
\label{e:CI_HMM_MD}
\end{equation}

The choice \eqref{e:HMMopt} might be  motivated by viewing the entire history $\clY_n^i = ( h (\Phi_0^i), \dots, h (\Phi_n^i) ) $ as a Markov chain, for each $i=0,1$, and then appealing to \Cref{t:bestF}.     
In fact, as in   the Markovian setting, theory in \cite{fuh03} begins with analysis of the   LLR associated with $ \clY_n^0,  \clY_n^1$, 
which may be expressed as \eqref{e:LLRn},  
 with $p^i_0$ the pmf for $h (\Phi_0^i)$.    
 
If  $D(  p_0^1 \|   p_0^0) $ and $D(  p_0^0 \|   p_0^1) $ are each finite, it follows from  \eqref{e:LLRn} that  the sign conventions are typically satisfied:  analogous to \eqref{e:MarkovA1},
 \[
  m_i  = \lim_{n\to\infty}\frac{1}{n}  \Expect_i\Big[     \sum_{k=1}^n  F_k  \Big] =
 \begin{cases}
  - 	\clK(\trajP_0  \|  \trajP_1)    \le 0  &  i=0
   \\
 \phantom{-} \clK(\trajP_1  \|  \trajP_0)  \ge 0  & i=1
\end{cases}
\]
where $i=0$ corresponds to $\tchange = \infty$ and $i=1$ corresponds to $\tchange = 0$.

By construction we have $\clK(\trajP_1\|\twtrajP_0^F) =0$ for this choice of $\{F_n \}$.   
Based on \Cref{t:Jentropy} it is likely that \Cref{t:bestF} can be extended to   the C.I.\ HMM and POMDP models.

\subsection{Approximation over a finite-dimensional function class}
\label{s:optFiniteDim}

Throughout the remainder of this section we consider a  linear  function class  that is finite dimensional, of dimension $d$.  
Let $\psi\colon\ystate\to\Re^d$ denote a basis, so that $\clG = \{ F_\theta = \theta^\transpose  \psi : \theta\in\Re^d \}$.  
We continue to assume that $\del = 1$ so that $F_\theta \colon\ystate\to\Re$.

\wham{Notation}
Subject to Assumption~\Ass{s:thetaCUSUM} we adopt the following notational conventions:

\whamb 
$m_i^\theta \eqdef \marg_i(F_\theta)$ for $i=0, 1$ and $\theta\in\Re^d$.

 \whamb There is  $\thone \in\Re^d$ such that ${\thone}^\transpose \psi\equiv 1$,  and   $\thmix \in\Re^d$ such that   
$ \marg_0(F_{\thmix} ) < 0$ and $ \marg_1(F_{\thmix} ) >0$.

\smallskip

The following notation is adopted when $F_\theta$ satisfies    
$m_0^\theta \le 0$ and $\Prob \{ F_\theta(\preObs_k) > \expa \} >0$.

 \whamb 
Let   $\thexp_+^\theta > \thexp_0^\theta  \ge 0$   denote the solutions to
\begin{equation}
\Lambda_0(\thexp_0^\theta F_\theta) = 0 \,, \quad 
\Lambda_0(\thexp_+^\theta F_\theta) = \expa
\label{e:thexp+theta}
\end{equation}
with $\thexp_0^\theta =0$ only when $m_0^\theta=0$.

 \whamb 
Let $\twmarg_0^{\thexp ,\theta}$ denote the probability measure \eqref{e:twmarg0} obtained using $\thexp $ with function $F_\theta$.      This is expressed  $  \twmarg_0^{\theta}$  for $\thexp =1 $ and 
 $  \twmarg_0^{ +,\theta}$ when using  $\thexp_+^\theta$.
 Denote $\cm_0^{+,\theta} =  \twmarg_0^{+,\theta} (F_\theta)$.   
		\qed

The notation \eqref{e:barJ_infty_F} is also modified in this subsection:
\begin{align}
\barJStarCoarse( \kappa ; \theta)  &=
\frac{1}{m_1^\theta }\frac{1}{\thexp_+^\theta}  \Big[  \log \kappa    +\frac{1}{2}\log\log\kappa\Bigr]  
\label{e:barJ_infty_theta}
\end{align}

A minimum of \eqref{e:barJ_infty_theta} solves  a moment matching problem:
\begin{proposition}[Moment matching under twisted dynamics]
\label[proposition]{t:calcOpt}
Suppose that  Assumption~\Ass{s:thetaCUSUM} holds, and
			that $\theta^\bullet\in\Re^d$    is a stationary point:    $ \nabla_\theta  \barJStarCoarse \, (   \kappa; \theta^\bullet ) = 0$ for some (and hence all)  $\kappa>0$.  Suppose moreover that    $\Fmax^0>\expa$ with $F= F_{\theta^\bullet}$ and 
  $  m_0^\theta \le  0$.
Then, $\cm_0^{+,\theta^\bullet}  =m_1^{\theta^\bullet}$ and more generally,
\begin{equation}
\twmarg_0^{+,\theta^\bullet} (F_\theta) = \marg_1(F_\theta)\, ,  \quad \theta\in\Re^d \, .
	\label{e:calcOpt}
\end{equation}
			\clThm
\end{proposition}

  \Cref{t:calcOpt} follows from \Cref{t:calcOptNon} in    the Appendix..    
follows quickly from \eqref{e:twistedMean} from which we obtain
\begin{equation}
\nabla_\theta
\Lambda_0  \,  (F_\theta)  =    \twmarg_0^{\theta} (\psi)
	\label{e:twistedMean_psi}
\end{equation}

 \smallskip

  There are several ways of characterizing an optimizer.  On  denoting $\barpsi^i_j  = \marg_i(\psi_j)$ for each $i, j$,

 \wham{1.}    Applying  \Cref{t:calcOpt},   
  an optimizer   $\theta^\circ$ must satisfy the moment matching constraints,
 \begin{equation}
	\twmarg_0^{+,\theta^\circ} \big( \psi_j \big) =  \barpsi^1_j       \,,   \quad1\le j\le d
	\label{e:calcOptLin}
\end{equation} 
However, 	\eqref{e:calcOptLin}
 is not sufficient for optimality.
 
 \wham{2.}
An optimizer is obtained as a solution to the convex program  \eqref{e:jStarConvexRelaxation}, which reduces to
\begin{equation}
 \begin{aligned}
\max \ \ &    \theta^\transpose \barpsi^1  
\\
\textrm{s.t.} \ \ 
& \Lambda_0(F_\theta) \le \expa \, ,   \  \  \theta^\transpose \barpsi^0 \le 0  \, ,   \   \theta\in\Re^d
\end{aligned}
 \label{e:jStarConvexRelaxationFD}
 \end{equation}

 \wham{3.}   A Lagrangian is defined for $(\theta,\zeta) \in \Re^d\times \Re^2_+$ by
 \begin{equation}
 \clL(\theta,\zeta) = -   \theta^\transpose \barpsi^1   + \zeta_1 [\Lambda_0( F_\theta) - \expa ] + \zeta_2 \theta^\transpose \barpsi^0 
 \label{e:LagrangeJFD}
\end{equation} 
 A dual of  \eqref{e:jStarConvexRelaxationFD} is then      
\[ 
 \sup_{\zeta \in\Re_+^2}
 \varphi^*(\zeta)
\,,
\]
 where     $\varphi^*(\zeta) = \inf \{  \clL(\theta,\zeta)  : \theta \in \Re^d \}$ is the  \textit{dual function}.

This dual is constructed by treating the primal \eqref{e:jStarConvexRelaxationFD} as a minimization problem  with objective  $ -   \theta^\transpose \barpsi^1  =-  \marg_1(F_\theta) $ and maintaining the constraints  in  \eqref{e:jStarConvexRelaxationFD}.

 Assumption \Ass{s:thetaCUSUM}~(iv) implies Slater's condition, so there is no duality gap. Given our modified objective this is expressed, 
 \[
  \sup_{\zeta \in\Re_+^2} \varphi^*(\zeta) = - \sup_\theta \theta^\transpose \barpsi^1  
 \]
 where the supremum on the right hand side is subject to the two inequality constraints in  \eqref{e:jStarConvexRelaxationFD}.

\begin{proposition}[Strong duality]
\label[proposition]{t:DualProps}
Suppose that \Ass{s:thetaCUSUM} holds.  Then,  with $\zeta^* = (1;0)$, 

\whamrm{(i)} 
  $\varphi^*(\zeta)= -\infty$ whenever $\zeta_1 + \zeta_2 \neq 1$,  and maximal for $\zeta^*  $ with value  
\begin{equation}
\begin{aligned}
\varphi^*(\zeta^*)  & =    -   \expa  +\inf_\theta \Gamma_0(\theta)     
\\
& \quad   \quad
\Gamma_0(\theta)   \eqdef  \Lambda_0(F_\theta)-\marg_1(F_\theta) 
\end{aligned}
\label{e:DualAndGamma}
\end{equation}
The lower bound holds, 
$ \varphi^*(\zeta^*) \ge  -   \expa  - 	\clK(\trajP_1\|\trajP_0) > -\infty$.

\whamrm{(ii)} 
$\theta^\circ $ solves the minimization problem \eqref{e:DualAndGamma} if and only if
$\nabla \Gamma_0(\theta^\circ) =0$, which is equivalently expressed 
  \begin{equation}
 \twmarg_0^{\theta^\circ} (\psi_j) =   \barpsi^1_j      \,,   \quad1\le j\le d
	\label{e:calcOptLinDual}
\end{equation}

\whamrm{(iii)}
Suppose that $\theta^\circ $ solves \eqref{e:DualAndGamma}.
  Let  $\theta^* = \theta^\circ + r^*\thone $ with   $r^* = \expa - \Lambda_0(F_{\theta^\circ}) $.
  Denote
  $\expa^0 =\clK(\trajP_0\|\twtrajP_0^F)  $ with $F = F_{\theta^\circ}$.  
  Then,  provided $ \expa  \le \expa^0$,
    \[
   \begin{aligned}
       \Lambda_0(F_{\theta^*}) =\expa\,, \quad
   \marg_0(F_{\theta^*}) & = -\expa^0 + \expa \le 0     &&  \text{ (feasibility)} 
\\
 \marg_1(F_{\theta^*}) &= - \varphi^*(\zeta^*) > 0  &&  \text{ (strong duality)} 
   \end{aligned}
 \]
 Consequently,  $\theta^* $ is a solution to the primal.
\clThm
  \end{proposition}

Observe that 	\eqref{e:calcOptLinDual} coincides with  \eqref{e:calcOptLin} only if $\thexp_+^{\theta}=1$ for $\theta = \theta^\circ$.       The construction in (iii) results in $\thexp_+^{\theta^*}=1$.

The relationship between the minimization of $\Gamma_0$ in \eqref{e:DualAndGamma} and the minimization of $\clK(\trajP_1\|\twtrajP_0^F) $ over $F$ explored in \Cref{t:jStarConvexRelaxationENT} is made clear through an application of 
\eqref{e:GammaK}:   
 \begin{equation}
\Gamma_0(\theta)    = \clK(\trajP_1\|\twtrajP_0^F) - \clK(\trajP_1\|\trajP_0) \,,\quad \textit{when $F=F_\theta$.}
\label{e:GammaKtheta}
\end{equation}		
 Hence the objectives differ by the fixed constant $ \clK(\trajP_1\|\trajP_0) $.

The following result will simplify computation in numerical examples that follow.

\begin{subequations}

\begin{proposition}[Optimal scalar offset]  
\label[proposition]{t:scalarOffsetOpt} 
Suppose that $G + r$  satisfies   Assumption~\Ass{s:extend} for some scalar $r$.      
For the one dimensional affine function class $\clG = \{ G+\theta : \theta \in\Re\}$,
a   value   $\theta^*$  minimizing \eqref{e:barJ_infty_theta}
 is     
 \begin{equation}
 \begin{aligned}
\theta^* & = \frac{\expa-\Lambda_0(\thexp_+G)}{\thexp_+}, 
\\
\textit{with} \ \ 
\thexp_+ & \in\argmax_{\thexp>0}\{\thexp\marg_1(G)-\Lambda_0(\thexp G)\}.
\end{aligned}
\label{e:OptOffset}
\end{equation}		
 Consequently, with  $F = G+\theta^*$ we have 
 \begin{equation}
 \begin{aligned}
  \Lambda_0(\thexp_+ F)  &  =  \expa \,,
  \\
     \cm_+ &  \eqdef
   \tfrac{d}{d \thexp}   \Lambda_0(\thexp F)  \Big| _{ \thexp = \thexp_+} =    \marg_1(F) 
 \end{aligned}
\label{e:ScalarOffsetConsequences}
\end{equation}
\clThm
\end{proposition}
\end{subequations}

\section{Numerical Results}
\label{s:numres}

The main results of the paper are illustrated here using two different models and various choices of $\{ F_n \}$ in the CUSUM statistic \eqref{e:CUSUM}.

 We do not have an analytic expression for the limiting overshoot $\oshoot_\infty$ defined in \eqref{e:oshootA4} and appearing   in the definition of $\aFine$ in \Cref{t:barcApprox}.  Consequently,  whenever needed this quantity was estimated via Monte-Carlo--see  \Cref{s:appndx_exp}  for details.

In most cases the identity $ \cm_+  = m_1$ is obtained on applying 
	\Cref{t:DualProps} or \Cref{t:scalarOffsetOpt},
 so that the definition in 
\eqref{e:RefinedPars} becomes $\gamma^2=\Upupsilon_0''(\thexp_+)/m_1^3$.

\subsection{Model 1.  Conditionally Independent Markov model}
\label{s:num_model1}

Model~1 is an instance of the C.I.\ Markov model described in \Cref{s:thetaCUSUM}, with observation space $\ystate=\Re$.  
 The change time is independent of $(\bfpreObs,\bfpostObs)$ and is chosen geometric so that $\Prob\{\tchange>n\}=\exp[-\expa(n+1)]$, $n\ge0$, with $\expa=0.02$. 
 
It is assumed that the independent Markov chains evolve as Gaussian linear processes, 
\begin{equation}
	X^i_{k+1}=A^iX^i_k+\sigma W^i_{k+1},
	\qquad i\in\{0,1\},\quad k\ge0
	\label{e:num_m1_ar}
\end{equation}
where $\{W^i_k\}$ is i.i.d.\ $N(0,1)$. 
Consequently the transition densities defining $L_\infty$ in \eqref{e:LLR_Markov} exist with 
$g_i(x,z)=\eta_\sigma(z-A^i x)$ and $\eta_\sigma$ is the $N(0,\sigma^2)$ density.  The numerical results use   $A^0=0.30$, $A^1=0.60$, and $\sigma=1$.

We take $\del=2$ and compare the outcomes of CUSUM using   different functions $F\colon\ystate^2\to\Re$ to define $F_{n+1} = F(Y_n, Y_{n+1})$ for $n\ge 1$.      \Cref{t:bestF} tells us that the optimal choice is  $F^* =  L_\infty + \expa$.    
We consider two others of the   form:
$F = \surL + r^*$  
in which $r^*$ was determined by \eqref{e:OptOffset}, and    for a density $\eta$ on $\Re$, 
\[
\surL(x,z) \eqdef \log\left({\surg_1(x,z)}/{\surg_0(x,z)}\right) 
\ \  \textit{with}  \ \ 
\surg_i(x,z)=\eta(z-A^i x) 
\qquad
\]

In each of the three cases the density $\eta$ was selected 
 with mean   zero and variance $\sigma^2 $:

\wham{Test 1a: Optimal.}
  $\eta=\eta_\sigma$, the density of an $N(0,\sigma^2)$ random variable.

  This is consistent with \eqref{e:num_m1_ar}:
    $\surg_i=g_i$, $\surL=L_\infty$, and $F=\Fstar=L_\infty+r^*$ with $r^*=\expa$.

\wham{Test 1b: Laplace.}
  $\eta$ is the Laplace$(0,b)$ density with $b=\sigma/\sqrt{2}$.

\wham{Test 1c: Student's $\bfmath{t}$.}
  $\eta$ is the Student's $t$ density with five degrees of freedom and scale $\sigma\sqrt{3/5}$ (see the Appendix for a formula).

\smallbreak

  Numerical values of $m_0,m_1,\thexp_+$, and $\gamma^2 $ were obtained in each of these three cases and   reported in   \Cref{s:appndx_exp}.
   In view of \eqref{e:barJFineBody} the cost approximation is linear in the value $(m_1\thexp_+)^{-1}$.   Computations confirm that Test~1a   maximizes the product    $ m_1\thexp_+ $, with  Test~1c coming in second place:   
 \[
	\text{Gaussian:}  \ \  0.090\,, \quad
	\text{Laplace:}   \ \ 0.070 \,, \quad
	\text{Student's $t$:} \ \    0.083     
\]

Simulation experiments were conducted to estimate the performance  $\barJ^*(  \kappa) $ and associated  thresholds $\barthresh^*(\kappa)  $ for a range of $\kappa$.   The two sets of plots shown in \Cref{f:model1_approx} compare these quantities with the approximations $  \barJFineStar ( \kappa) $, $   \barthreshFineStar(\kappa)   $    defined in \eqref{e:barFineBody}.     The plots illustrate the approximations  \eqref{eq:refined_bdds} established in \Cref{t:barcApprox}.

\subsection{Model 2. Relaxing  Statistical Independence}
\label{s:Model2Relax}

We consider an instance of the POMDP model of \Cref{s:POMDP}.
Let $\bfmZ$ be a Markov chain on $\zstate=\{0,1,2\}$ whose transition matrix is taken of the form
\begin{equation}
	Q
	=
	\begin{bmatrix}
		1-\epsy_0-q_{01} & q_{01} & \epsy_0\\
		q_{10} & 1-\epsy_1-q_{10} & \epsy_1\\
		0 & 0 & 1
	\end{bmatrix}.
	\label{e:num_model2_Ztransition}
\end{equation}
This is designed so that   state $2$ is absorbing, and it is assumed that   $0<\epsy_i<1$ for  each $i$ so that $\bfmZ$ is   uni-chain.  

 The observations evolve on $\ystate=\{0,1\}$ according to
$Y_k=h_0(Z_k)\oplus B_k$, where $\oplus$ denotes addition modulo two,
   $\bfmB$ is i.i.d.\ Bernoulli$(\epsy)$, independent of
$\bfmZ$,  and
$h_0(z)=\ind\{z\in\{1,2\}\}$.   Consequently,
\begin{equation}
\begin{aligned}
	g(y\mid z)
	&\eqdef 
	\Prob\{Y_k=y\mid Z_k=z\}
	\\
&	=
	(1-\varepsilon)\ind\{y=h_0(z)\}
	+\varepsilon\ind\{y\ne h_0(z)\}.
\end{aligned}
	\label{e:num_model2_g}
\end{equation}

The choice of Markovian state is not unique; here we take  $\Phi_k=(Z_k,Y_k)$, so that $h(\bst)=y$ for any pair $\bst = (z,y)$.  
Its transition matrix   is $P(\bst,\bst')=Q(z,z')g(y'\mid z')$ for any  $\bst = (z,y)$, $\bst' = (z',y') \in\bigstate = \zstate\times\ystate$ 

All of the assumptions of \Cref{t:POMDPhazard} hold:  $\bfPhi$ is uni-chain and aperiodic on 
$\bigstate=\bigstate_0\cup\bigstate_1$ with $\bigstate_i=\zstate_i\times\ystate$, so that 
 $\bigstate_1$   is absorbing as required.   
The change time is by definition
$\tchange=\min\{k\ge0:\Phi_k\in\bigstate_1\}
=\min\{k\ge0:Z_k=2\}$, so that  \Cref{t:POMDPhazard} verifies that the limit \eqref{e:hazardAss} exists, and    $ \expa=-\log\PFeig=0.0166$.   
Details regarding computation of the PF eigenvalue $\PFeig$ and other quantities used to obtain the transition matrix   $\oP$ defined in \eqref{e:twisted} are
reported in  \Cref{s:appndx_exp} along with   details regarding   computations required in the following.

We next describe several QCD tests used to illustrate the theory.  The first two are based on the information state introduced in
\Cref{s:POMDP}.  In Model 2 this reduces to 
$\condDist_k(z)=\Prob\{Z_k=z\mid Y_0^k\}$, $z\in\zstate$, 
which evolves on the simplex $ \clS^3 =\{\beta\in \Re_+^3 :  \sum \beta_i  =1\}$.

 \begin{figure}[h]
  	\includegraphics[width=0.65\hsize]{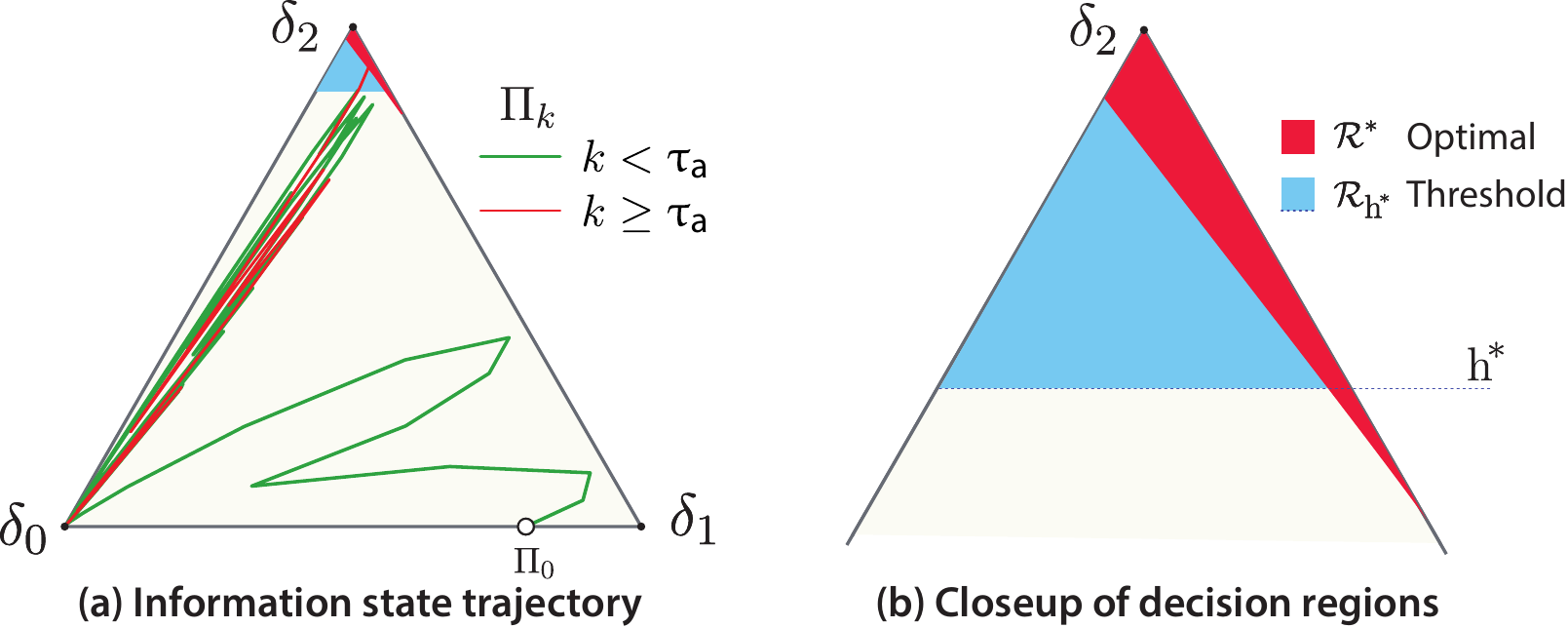} 
 	\centering
 	\caption{(a)  Typical sample path of
 		$\{\condDist_k:k\ge0\}$ evolving on the simplex.  
		(b) 
	Comparison of the decision regions for Tests~2a and 2b.   }
 	\label{f:model2_simplex_traj}
 \end{figure}

\wham{Test 2a:   Baseline.}  The optimal solution
$U_k^*=\fee^*(\condDist_k)=\ind\{\condDist_k\in\clR^*\}$ was
approximated, where
$\clR^*=\{\beta\in\clS^3:\fee^*(\beta)=1\}$.  
The approximation was obtained by  binning the two dimensional state space  $\clS^3$ and applying
policy iteration to obtain the optimal policy on the resulting finite state-space model.

\wham{Test 2b: Shiryaev's policy.}   This takes the form of \eqref{e:SyiryaevPolicy}, even though the conditions for optimality obtained in \cite{shi77,yak94} are not satisfied.    For $\threshSmall\in(0,1)$  denote $U_k^{\threshSmall}=\ind\{ p_k \ge \threshSmall\}   = \ind\{\condDist_k\in\clR_{\threshSmall}\}$, where  $\clR_{\threshSmall}=\{\beta\in\clS^3:\beta(2)\ge \threshSmall\}$.    For each $\kappa >0$ we choose $\threshSmallStar=\threshSmallStar(\kappa)$ to minimize the cost \eqref{e:MDD+kappaMDE}---these thresholds were estimated via   Monte-Carlo.

\Cref{f:model2_simplex_traj} compares the decision regions
$\clR^*$ and $\clR_{\threshSmallStar}$  for $\kappa=2$, along with a sample path of $\condDist_k$.

\wham{Test 2c: C.I.\ HMM approximation.}
Recall the C.I.\ HMM approximation in \Cref{t:POMDPhazard} in which $\bfmX^1 = \bfmX^\infty$, defined  as the stationary realization of $X_k^\infty = h(\Phi_k^\infty)$ with $\bfPhi$ evolving on $\bigstate_1$,
and $X_k^0 = h( \oPhi_k)$ with $\bfoPhi$ evolving on $\bigstate_0$.  

We consider   the choice   \eqref{e:HMMopt} designed based on this approximation:
\begin{equation}
	F^{(\text{2c})}_{n+1}
	=
 \log  \frac{p^1_{n+1\mid n}(Y_{n+1}   ) }{  \op^0_{n+1\mid n}(Y_{n+1}     )  }
	+\expa \,, 
	\label{e:HMMopt2b}
\end{equation}
in which  $p^1_{n+1\mid n}( y)= \marg_1(y) =  g(y\mid2)$, $y\in\ystate$, since $\{ X^\infty_k \}$ is i.i.d.\ in this example.   The conditional pmf  $\op^0_{n+1\mid n}$ is 
defined precisely as   in \eqref{e:CPY} with $i=0$,    based on the  observation model $ \oY_k = h(\oPhi_k)$ in which
the transition matrix of the Markov chain $\bfoPhi$ on $\bigstate_0$ is defined in \eqref{e:twisted}.

\wham{Test 2d: Finite memory approximation.}
A finite-memory approximation to the previous test was considered,
\begin{equation}
	F^{(\text{2d})}_{n+1}	=
	 \log  \frac{p^{1,\del}_{n+1\mid n}(Y_{n+1}   ) }{  \op^{0,\del}_{n+1\mid n}(Y_{n+1}     )  }
 	+ \expa \,,
	\label{e:HMMopt2c}
\end{equation}
defined such that $F^{(\text{2d})}_k$ is a function of $\strY{\del}_k = (Y_{k-\del+1},\dots, Y_k)$ for each $k\ge \del$.

 \begin{figure}[t]
 	\centering
 	\includegraphics[width=0.7\hsize]{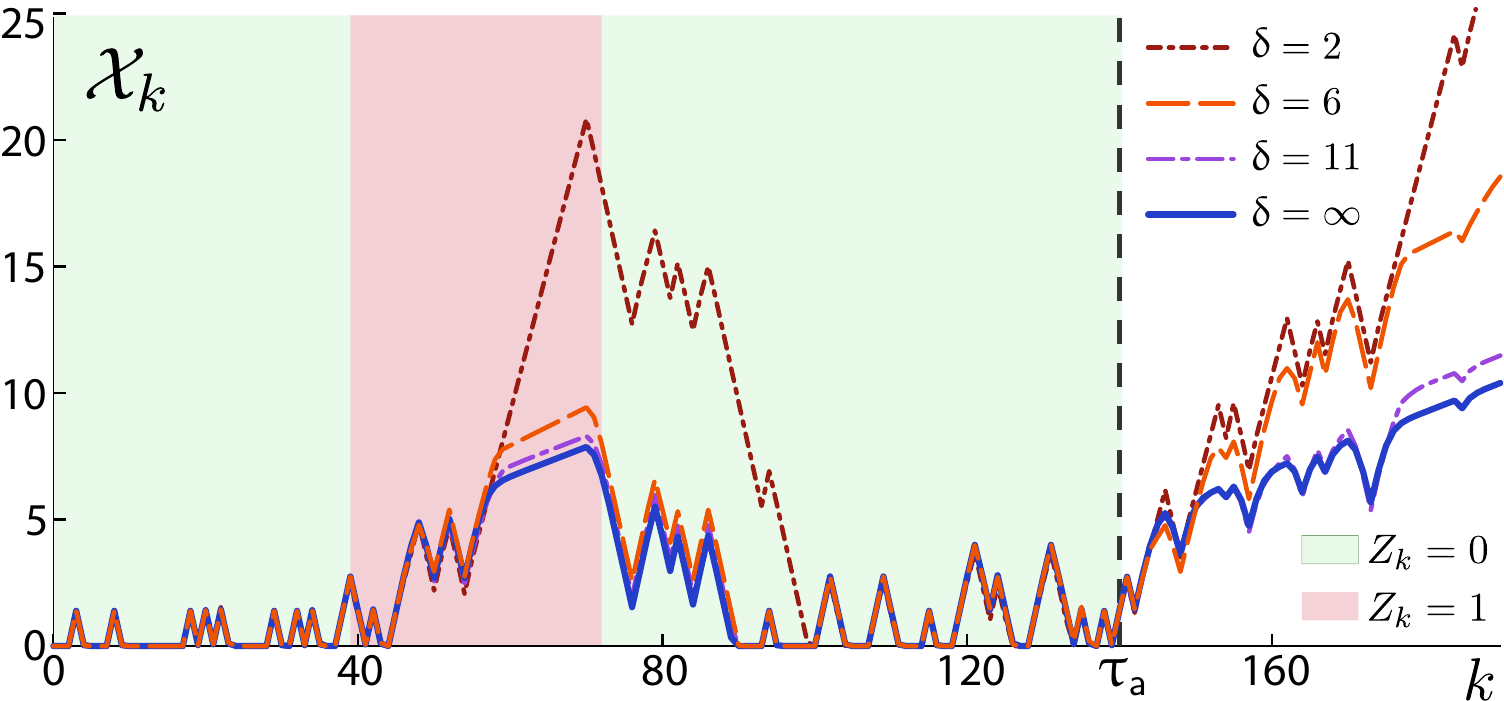}
\caption{Sample paths of the CUSUM statistic for   Test~2d for three choices of   $\del$ in \eqref{e:HMMopt2c},
	and for Test~2c (indicated by $\del=\infty$).
	} 	\label{f:model2_cusum_traj}
 \end{figure}

We take   $ p^{1,\del}_{n+1\mid n} (y) =  g(y\mid2)$ for any $\del\ge 1$ and, on observing $y_k = Y_k$ for $0\le k\le n$,  we define 
\[
\op^{0,\del}_{n+1\mid n}(y)
\eqdef
\Prob\{\oY^\infty_{n+1}=y\mid \oY^\infty_k = y_k\,, \ [n-\del + 2]_+ \le k\le n  \}
\]
where the superscript in 
 $ \oY^\infty_k = h(\oPhi_k^\infty)$ indicates that $\bfoPhi^\infty$ is stationary with marginal    $
\ometamarg$  (recall   \eqref{e:ExplicitApprox}).    In the special case $\del=1$ this is interpreted as  $  \op^{0,\del}_{n+1\mid n}(y) = \Prob\{  h (\oPhi_k^\infty) =y \}$,

Sample paths of  the CUSUM statistic are shown in \Cref{f:model2_cusum_traj} for several values of $\del$, with $\del = \infty$ referring to Test~2c.   
We observe that the pre- and post-change means satisfy the sign constraints $m_0<0$, $m_1>0$, 
though we can verify this analytically only for  $\del = \infty$ and sufficiently small $\expa>0$.

\wham{Test 2e: Optimization over a linear function class.}
We take $\del =2$ and let $\clG$ denote the set of all functions on $\ystate^2$.    Denote $F^{(\text{2e})}_{n+1}=\Fstar(Y_n,Y_{n+1})$, where $\Fstar$ is a solution to the primal \eqref{e:jStarConvexRelaxation}, which according to \Cref{t:CoarseConvex}  satisfies $\Lambda_0(\Fstar) = \expa $.

Given any function $F\colon\ystate^2\to\Re$     the twisted law  $\twtrajP_0$ is defined in \Cref{s:extend},  
from which we obtain the   bivariate marginals,   
\begin{equation*}
\cmarg_0^{[2]} (y,y') =  \Prob\{ \cY_k  = y\,, \ \cY_{k+1} = y' \}  
\end{equation*}
where the probabilities are with respect to  $\twtrajP_0$ and hence independent of $k$. 
Application of \Cref{t:completeOpt}  tells us that for the optimizing function $F^*$ we must have, for
 $y,y'\in\ystate$,
 \begin{equation}
\cmarg_0^{[2]} (y,y')  = \marg_1^{[2]} (y,y')
	=
	g(y\mid2)g(y'\mid2) 
\label{e:model2bivariateMatching}
\end{equation}
Computations in
 \Cref{s:appndx_exp}
   reveal that the total variation error is less than $10^{-8}$.

\begin{figure*}[h]
	\centering
	
	\includegraphics[width=0.95\hsize]{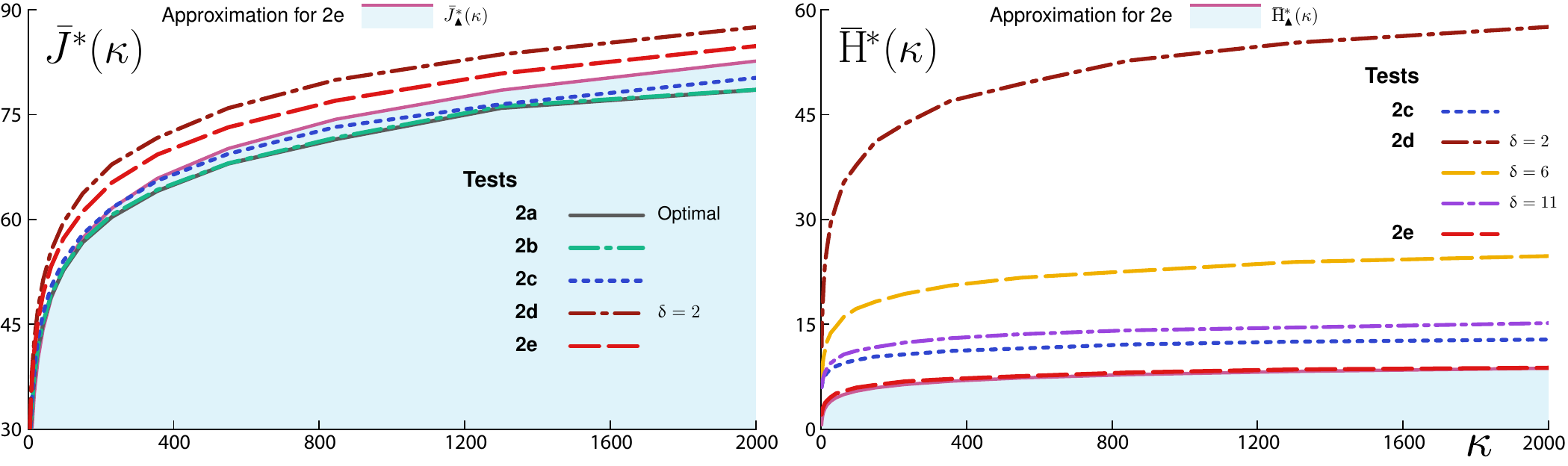}

	\caption{Performance and thresholds in Model~2. On the left is shown empirical performance for selected tests, and the   approximation \eqref{e:barJFineBody} obtained for Test~2e.  	Thresholds are shown on the right for selected tests,  as well as the approximation \eqref{e:barHFineBody}
for Test~2e. 
 }
	\label{f:model2_approx}
\end{figure*}

For purposes of computation we order the elements of $\ystate^2$ as
$\{(y_i,y'_i):1\le i\le4\} = \{ (0,0),(0,1),(1,0),(1,1) \}$,
and define $\psi_i(y,y')=\ind\{y=y_i,\ y'=y'_i\}$, so that $\clG=\{\theta^\transpose\psi:\theta\in\Re^4\}$.  
The optimization over $\clG$ to obtain $\Fstar = F_{\theta^*}$ for some $\theta^*\in\Re^4$
can be reduced to two dimensions.  
The proof of \Cref{t:Model2clG} is postponed to the Appendix.

\begin{lemma}
	\label[lemma]{t:Model2clG}
Denote  
$\thzz=(0;1;-1;0)$ and $\thone  = (1;1;1;1)$.  Then, 
 $	\Gamma_0(\theta+r\thone+s\thzz)=\Gamma_0(\theta)   
$ for  
	$\theta\in\Re^4$ and $r,s\in\Re$.  \clThm
\end{lemma}

%
 
Fix any $e^i\in\Re^4$, $i=1,2$, such that $\{\thone,\thzz, e^1,e^2\}$ forms an orthogonal basis for $\Re^4$.   It follows from \Cref{t:Model2clG} that to obtain a minimizer $\theta^\circ$ of $\Gamma_0$ it suffices to obtain a pair $(s^\circ,t^\circ)$ that minimizes $ \Gamma_0(se^1+te^2)$   over $ (s,t)\in\Re^2$,   and then take $\theta^\circ = s^\circ e^1+t^\circ e^2$.

\begin{wrapfigure}[18]{r}{0.45\textwidth}
	\includegraphics[width=\hsize]{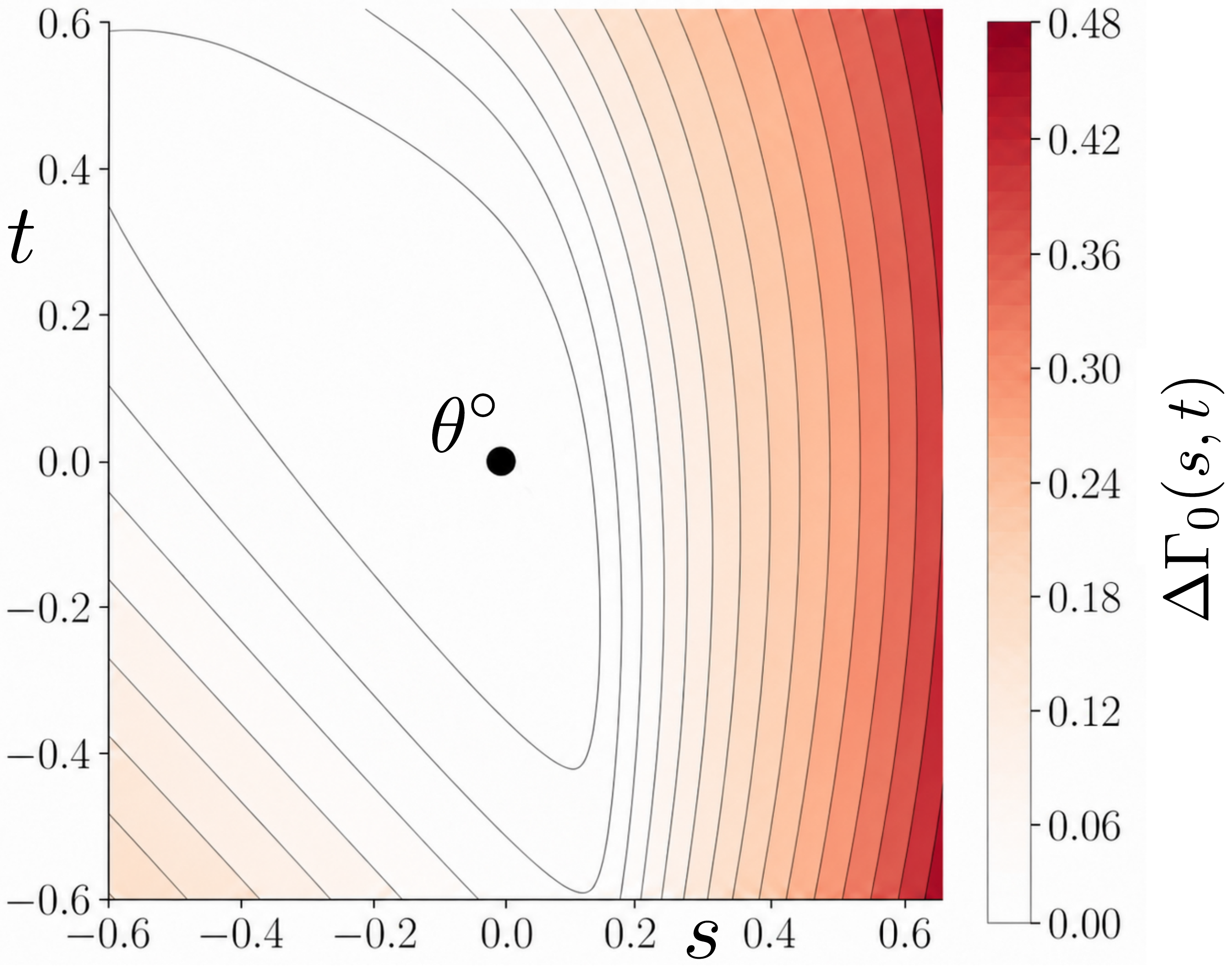}
\caption{Contour of $\Delta\Gamma_0$ in
	\eqref{e:num_model2b_contour_objective} for Test~2e.}
	\label{f:contour}
\end{wrapfigure}

  In computations we chose   $e^1 = (1;0;0;-1)/\sqrt{2}$ and  $e^2 =  (-1;1;1;-1)/2$.       \Cref{f:contour}  displays a contour plot of the error,
\begin{equation}
	\Delta\Gamma_0(s,t)
	=
	\Gamma_0(\theta^\circ+se^1+te^2)
	-
	\Gamma_0(\theta^\circ).
	\label{e:num_model2b_contour_objective}
\end{equation}

\Cref{t:DualProps}~(iii) implies that $\theta^* =	\theta^\circ+r^\circ \thone$ using  $r^\circ=\expa-\Lambda_0(F_{\theta^\circ})$ solves the primal \eqref{e:jStarConvexRelaxationFD}.   Moreover, another application of  \Cref{t:Model2clG} implies that the optimizer is not unique:
 $\theta^* + t\thzz $ also satisfies the conditions of 	\Cref{t:DualProps}~(iii), and hence also solves the primal, for any $t\in\Re$.

\wham{Comparison of the tests.}

As in Model~1,   experiments were conducted to estimate the performance  $\barJ^*(  \kappa) $ and associated  thresholds $\barthresh^*(\kappa)  $ for a range of $\kappa$.   The cost plot also includes the optimal policy $\fee^*$ described as test 2a.   
 
 The poor performance for Test~2d with small $\del$ is explained in part by the fact that   $\expa$ is used in \eqref{e:HMMopt2c} rather than the optimal offset characterized in \Cref{t:scalarOffsetOpt}.

 The approximations $  \barJFineStar ( \kappa) $ and $   \barthreshFineStar(\kappa)   $ are displayed only for Test~2e.   
 As in Model~1, the threshold approximations are nearly perfect.    There is a small   under-estimate of performance.

 In Tests 2c and 2e  we have by construction $ \thexp_+ = 1$  so that the approximations  \eqref{eq:refined_bdds} reduce to 
 \[
 \begin{aligned}
 \barthreshFineStar(\kappa) &=  \log\kappa+\frac{1}{2}\log\log\kappa+\bFine
 \\  
 \barJFineStar(\kappa)
			&=
\Bigl[\log\kappa+\frac{1}{2}\log\log\kappa+ \aFine\Bigr] \frac{1}{m_1}
\end{aligned}
 \]
 We do not yet have theory to justify these formulae for Test~2c since   \Cref{t:barcApprox} requires $F_n = F(\strY{\del}_n)$ with $ \del $ finite.  However, we can test the approximation empirically.    

The mean $m_1$ for   Test~2c may be expressed,
\[
m_1^{2c}  =  \lim_{N\to\infty}  \frac{1}{N}   \sum_{k=1}^N	F^{(\text{2c})}_k  =  \lim_{n\to\infty}  \frac{1}{n}  \InfoState_n^{(\text{2c})}
\]
where in each case the sample paths are generated with  $\tchange=0$ (equivalently,  $Z_0 =2$).     Monte-Carlo estimates based on the second characterization resulted in  $m_1^{2c} =  0.113$, and the value $m_1^{2e} = 0.111$ was obtained analytically.    The ratio  $0.113/0.111 = 1.02$   tells us that 
Test~2e has a cost of about 2\%~larger than for Test~2c when   $\kappa>0$ is large.

\section{Conclusions}
\label{s:conc}

The cost and threshold approximations \eqref{eq:refined_bdds} lead to a convex program for selection of the function used in a \textit{mismatched}  
CUSUM decision rule.  While the main results are initially cast   in a simple hidden Markov setting, the LDP setting allows extension to more complex models.    It is most surprising that the logarithmic approximation \eqref{e:MDEstep0.1} leads to the approximations \eqref{eq:refined_bdds} with vanishing error as $\kappa\to \infty$.   

The numerical results suggest that these asymptotic approximations remain accurate even for moderate values of \(\kappa\).

The optimization theory developed in \Cref{s:thetaCUSUM}
 yields extensions of classical 
likelihood-ratio constructions, such as the conclusion $\marg_1 =  \twmarg_0^{F^*}$ in \Cref{t:completeOpt}.  

As observed in prior literature on CUSUM for QCD,  an attractive feature of the framework is that it requires relatively little
knowledge of post-change behavior.  This is important in applications for
which normal operation is well understood, while faults, attacks, and other
abnormal events may take many forms.  Examples include transmission-line
faults and cyber attacks, where a reliable parametric model of post-change
behavior may be unavailable.    The introduction of mismatched CUSUM* also reduces the required knowledge regarding pre-change behavior.  

Current research is focused on efficient methods for the minimization of \eqref{e:barJ_infty_theta} using Monte-Carlo methods as motivated by the structure established in \Cref{s:optFiniteDim}.  This may be based on standard approaches to reinforcement learning,  the more recent feedback particle filter approach \cite{jostagmehmey22}, and may be  combined with prior techniques to approximate $\Lambda_0$ and its roots \cite{dufmet05a,dufmet05b,dufwil17,dufmey11}.   
Preliminary results in this direction may be found in
\cite{coomey24a,coomey24d}.

\bibliographystyle{abbrv}

\def\cprime{$'$}\def\cprime{$'$}

\appendix
		

\section{Appendix} 

This Appendix contains proofs of technical results and some details on numerical experiments.
Throughout this section we  let $o(1) $ denote an error term that vanishes as   $\thresh\to\infty$.

\subsection{LDP theory and    \Cref{t:barcApprox}}
\label{s:BigLDPproof}

Recall that in   \Cref{t:barcApprox}  it is assumed that $F_n = F(\Obs_n)$, for a  function $F\colon\ystate\to\Re$ 
(so that $\del=1$).

\begin{lemma}
\label[lemma]{t:MDD_renewal}    
Consider the POMDP model satisfying 
\Ass{s:POMDP},   or the C.I.\  model satisfying   \Ass{s:extend}~(i) and (ii).   Then, the mean detection delay satisfies  
$\MDD =    \thresh/m_1 + \oshoot_\infty  + o(1)  
$. 
\end{lemma}

\PF 
For the C.I.\ model this   is immediate from \eqref{e:oshootA4}.
For the POMDP model,  arguments   in \cite{fuhlai01} can be adapted to establish the result, with $\oshoot_\infty$ depending on the initial condition $\Phi_0$.  
\qed

We turn next to analysis of mean detection eagerness $\MDE$. 
We begin with justification of \eqref{e:MDEstep0}  
that invites application of techniques from LDP theory.

\begin{lemma}
\label[lemma]{t:MDE1}
The expression \eqref{e:MDEstep0} holds, and hence  for the C.I.\ model
\begin{equation}
\MDE  =    \sum_{n=0}^{ \infty }    \Prob_0 \{ \tstop\le  n  \}    \Prob\{ \tchange  >  n  \} 
\label{e:MDEstep2}
\end{equation} 
where $\Prob_0$ denotes the law    when $\tchange \equiv  \infty$.   
 \end{lemma}

 \wham{Proof}   
For any nonnegative integers \(a,b\),
\[
        (a-b)_+
        =
        \sum_{n=0}^{\infty}
        \ind_{\{b\le n<a\}}
        =
        \sum_{n=0}^{\infty}
        \ind_{\{b\le n\}}
        \ind_{\{a>n\}} .
\]
 Applying this identity with \(a=\tchange\) and \(b=\tstop\) gives the desired results.    
\qed    
   
Our next task is to approximate $ \Prob_0 \{ \tstop\le  n  \}   $, which involves the function $e_0$ defined in \eqref{e:rate1}.
The following is immediate from the upper bound   $ \Upupsilon_0(\thexp) \le \Fmax^0 \thexp$ when $\thexp>0$.
\begin{lemma}
\label[lemma]{t:e0infinite}
Subject to \Ass{s:extend} we have $e_0(t) =  t I_0(1/t)  = \infty$ for   $0<t < \Tmin \eqdef 1/( \Fmax^0)$.
\clThm
 \end{lemma}

 When the supremum defining $I_0(m)$ in  in \eqref{e:I0} is attained it is denoted   $\thexp(m)$:
\begin{equation}
			I_0(m) =   \sup_\thexp [  \thexp m - \Upupsilon_0(\thexp   )  ]  =   \thexp(m) m - \Upupsilon_0(\thexp(m)) 
\label{e:I0m}
\end{equation}
 We have   $\thexp_0  = \thexp(\cm_0)$ and  $\thexp_+ = \thexp(\cm_+ )$, giving
 \begin{equation}
			I_0(\cm_0)  =  \cm_0\thexp_0 \,, \ \  I_0(\cm_+)  =  \cm_+\thexp_+ - \expa
\label{e:I0cm0+}
\end{equation}

\begin{lemma}
\label[lemma]{t:MostLikelyPaths}   
Suppose that  Assumption~\Ass{s:POMDP} holds for the POMDP model or
Assumption~\Ass{s:extend} holds for the C.I.\ model.    We then have for $T> 0$,
\begin{equation}
  \log 
\Prob_0\Bigl\{ \sup_{0\le t\le T}   \xH_t     \ge 1  \Bigr  \}  =  - [   e_0(T)  + \epsy(T,\thresh) ]  \thresh   
\label{e:LDP}
\end{equation} 
in which the error term $\epsy$ satisfies   
\begin{equation}
\lim_{\thresh\to\infty}  \sup_{T_0\le T\le T_1}    \epsy(T,\thresh)   =  0
\label{e:LDerrorUni}
\end{equation}   
whenever 
 $\Tmin = 1/( \Fmax^0) <T_0 <T_1<\infty$. 
\end{lemma}

\wham{Proof}
Under either set of conditions we conclude that 
\Ass{s:extend}~(iv) holds.    
The next step is to apply the  contraction principle, see for example  \cite[Section~6.4]{ganocowis04},  
which implies the following limit for each $T$:
\[
\lim_{\thresh\to\infty} 
\thresh^{-1} 
  \log 
\Prob_0\Bigl\{ \sup_{0\le t\le T}   \xH_t     \ge 1  \Bigr  \}  =  -    e_0(T) 
\]
 
Observe that $\Prob_0\Bigl\{ \sup_{0\le t\le T}   \xH_t     \ge 1  \Bigr  \} $ is non-decreasing in $T$, and the function $e_0$ is continuous on $(\Tmin,\infty)$,   which implies uniform error bound \eqref{e:LDerrorUni}.
\qed

We have under the assumptions of \Cref{t:MostLikelyPaths}, 
\begin{equation}
\!\!\!
\lim_{\thresh\to\infty}  \frac{1}{
\thresh} 
  \log   \Big[     \Prob_0\{  \max_{0\le t\le s_n}  \xH_t \ge 1 \}    \Prob\{ \tchange  >  n  \}   \Big]
 =  -   G(s_n)
\label{e:preMDEsumApprox}
\end{equation}
with $G$ defined in \eqref{e:rate1G}.
Approximations of $\MDE$ are based on    \eqref{e:MDEstep1}, and  properties of $G$ summarized in the following:

\begin{subequations}

\begin{lemma}
\label[lemma]{t:Gprops}
Under the assumptions of \Cref{t:MostLikelyPaths} the function $G$  defined in \eqref{e:rate1G} has the following properties:

\whamrm{(i)}  It is convex on $\Re_+$,  continuously differentiable on the open interval $(\Tmin,  \infty )$ and     satisfies 
for $			 \Tmin < s < s^0 \eqdef 1/\cm_0 $,
\begin{equation}
\begin{aligned}
			G'(s)  &= I_0(1/s)  -  I'_0(1/s)/s  + \expa \,, 
			\\
			 G''(s) & =   I''_0(1/s)/s^3  >0
\end{aligned}
	\label{e:G12der}
\end{equation} 

\whamrm{(ii)} 
It is minimized at $s^*=1/\cm_+$,  and  
\begin{equation}
	G(s^*)   = \thexp_+ \,,  \qquad  G(s^0)  =  \thexp_0  +  \expa s^0
\label{e:Gts}
\end{equation}

   \whamrm{(iii)} 
   It admits the Taylor series approximation,
\[
\tilG(s) \eqdef G(s) - G(s^*) = \half (s-s^*)^2/\gamma^2    +O(s-s^*)^3
\]
	where $\gamma^2 =  1/G''\, (s^*) =  \Upupsilon_0''(\thexp_+)/\cm_+^3$.
\end{lemma} 
 \end{subequations}

\wham{Proof}
The derivative formulae 	\eqref{e:G12der}
follow   from  the definition of $G$.

To establish that $G$ is continuously differentiable first note that it is continuous on $(\Tmin,  \infty )$,  and affine on the interval $[s^0,\infty)$ with slope $\expa$.  It is therefore enough to establish the limit
\[
\expa = G'(s-) \eqdef 
\lim_{s\uparrow s^0} 
G'(s) = I_0(\cm_0)  -  \cm_0 I'_0(\cm_0)  + \expa
\]
	The identity  $I_0'(m) = \thexp(m) $ follows from convex duality, and   $I_0(m) = \thexp(m) m - \Upupsilon_0(\thexp(m))$ follows from the definition of $\thexp(m) $.   Consequently,   $ G'(s-) =  I_0(\cm_0)  -  \cm_0  \thexp_0 + \expa  = - \Upupsilon_0( \thexp_0 ) + \expa =  \expa $ as desired.

To establish (ii) first note that the function $G$ is strictly increasing for $s\ge s^0$,   which justifies consideration of only $s <s^0$. 
	
We make the change of variable $m =   1/s$,  so that  $m \ge \cm_0$ and   $e_0(s)     +  s \expa =  [ I_0(m) + \expa] /  m$
 whenever $s<s^0$.
     On differentiating with respect to $m$, the minimizer  $m^*$ must solve
\[
		0 =   \Big(  m I_0'(m) - [ I_0(m) + \expa]  \Big)/m^2  \Big|_{m = m^*}
\]	
Substituting  $I_0'(m) = \thexp(m) $ and   $I_0(m) = \thexp(m) m - \Upupsilon_0(\thexp(m))$ then gives
\[
		0 =   \Big(  m  \thexp(m)  - [  \thexp(m) m - \Upupsilon_0(\thexp(m))  + \expa]  \Big)  \Big|_{m = m^*}
\]
This combined with \eqref{e:I0cm0+}
establishes that $m^* = \cm_+$ is the unique minimizer.   
		
		Finally,
		from the definitions we have  $G(s^*) =  s^* I_0(1/s^*)  + s^*\expa $,   with $s^* = 1/\cm_+$.    
		On applying \eqref{e:I0cm0+}, which gives $I_0(\cm_+)  =   \cm_+ \thexp_+ -  \expa$, we obtain the first identity in \eqref{e:Gts}:
\[
		G(s^*)  =   [ \cm_+ \thexp_+ - \Upupsilon_0(\thexp_+) +\expa ]/\cm_+ =\thexp_+ 
\]
		The second identity   follows directly from   $   I_0(\cm_0) =     \thexp_0 \cm_0 - \Upupsilon_0(\thexp_0 )   =  \thexp_0 \cm_0 $.  
		
		This establishes (ii) and (iii) immediately follows.  
\qed

Our next results are successive refinements of 
\eqref{e:preMDEsumApprox}. 		
 
 Throughout the following we fix $\epsy \in (1/2,2/3)$ and denote $\delta_{\thresh} = 
\thresh^{-1+\epsy}$.  In view of \Cref{t:Gprops} the convex function $\tilG(s) \eqdef G(s) - G(s^*) $ admits the following useful bounds:
There is a fixed constant $b_0$ and function of two variables $ R_{\thresh}(s) $ such that  
\begin{equation}
\begin{aligned}
\tilG(s)  & \ge     \frac{\delta_{\thresh}^2}{2 \gamma^2 } +  b_0 \delta_{\thresh}   [    | s - s^*| -\delta_{\thresh}   ]  +  O(\delta_{\thresh}^3)      &&   |s-s^*| \ge \delta_{\thresh}  
\\
\tilG(s)  & =     \frac{  (s-s^*)^2}{2\gamma^2 } +   R_{\thresh}(s)  &&   |s-s^*|\le \delta_{\thresh}  
\\
&  \qquad  \textit{where}  \quad \sup_{|s-s^*|\le\delta_{\thresh}} |R_{\thresh}(s)| = O(\delta_{\thresh}^3)
\end{aligned}
\label{e:tilGbdds}
\end{equation}
where we have used $G''(s^*)=1/\gamma^2$.
The lower bound is obtained from two Taylor series approximations about $s^* \pm \delta_{\thresh} $, from which we obtain 
$b_0 = 1/ \gamma^2  + O(  \delta_{\thresh} )$.   
 
These bounds combined with \Cref{t:Gprops} and \eqref{e:preMDEsumApprox} imply that 
 we can disregard most of the terms in \eqref{e:MDEstep1}.   
 	 
\begin{lemma}
\label[lemma]{t:MDEsumApprox}
Under the assumptions of \Cref{t:MostLikelyPaths}  we have 
\begin{equation*}
  \sum_{n=0}^\infty   \exp\bigl(   -    \thresh  G(s_n)  \bigr)  =   [1 +  o(1) ]  \sum_{| s_n -s^*| \le  \delta_{\thresh} }     \exp\bigl(   -    \thresh  G(s_n)  \bigr)
 \end{equation*}
 \end{lemma}

\wham{Proof}
We equivalently establish \Cref{t:MDEsumApprox}   for the centered function defined by
 $\tilG(s) = G(s) - G(s^*)= G(s) - \thexp_+ $ for $s\ge 0$.

 We first obtain a lower bound on the sum over  $n$ satisfying $| s_n -s^*| \le  \delta_{\thresh} $:  
Observing that the total number of $n$ satisfying this constraint is no less than $ 2 [\thresh^\epsy  -1] $,  we may apply
 \eqref{e:tilGbdds} to obtain
\[
 \begin{aligned}
\sum_{| s_n -s^*| \le  \delta_{\thresh} }   \exp\bigl(   -    \thresh  \tilG(s_n)  \bigr) & 
 \ge  2 [\thresh^\epsy  -1]   \exp(  - b_1 \thresh  \delta_{\thresh} ^2 ) 
\end{aligned} 
\]
where $b_1 = \delta_{\thresh} ^{-2}  \max \{  \half G''\, (s^*) (s-s^*)^2  +   R_{\thresh}(s)  :   |s-s^*|\le \delta_{\thresh}  \}$ is uniformly bounded in $\thresh\ge 1$ for any $\epsy \in (1/2,2/3)$.

On summing over the complementary region we    apply the lower bound in
 \eqref{e:tilGbdds}:
\[
 \begin{aligned} 
 &
\sum_{| s_n -s^*|\ge  \delta_{\thresh} }    \exp\bigl(   -    \thresh  \tilG(s_n)  \bigr) 
\\
& \le   
2 
\exp(   -  \half  \thresh     \delta_{\thresh}^2/\gamma^2  )
 \sum_{n=0}^\infty  \exp\big(   -    \thresh  \{   b_0 \delta_{\thresh} + O (\delta_{\thresh} ^3) \}(n/\thresh) \big)      
 \\
 & =  2 \exp(   -  \half  \thresh     \delta_{\thresh}^2/\gamma^2  )
   \frac{ 1}{1- \exp\big(   - b_0 \delta_{\thresh} + O (\delta_{\thresh} ^3 ) \big)   }
\\
&   =   \frac{2}{b_0     }  \exp(   -  \half  \thresh^{-1+2\epsy}  /\gamma^2  )      \thresh^{1-\epsy}  [1+ o(1)]
\end{aligned} 
\]

The ratio of the preceding upper and lower bounds admits the order bound
\[
O\Big(
	H^{1-2\epsy}
	\exp\Big\{
		\Big(b_1-\frac{1}{2\gamma^2}\Big)
		H\delta_{\thresh}^2
	\Big\}
\Big)
=
O\Big(
	H^{1-2\epsy}
	e^{O(H\delta_{\thresh}^3)}
\Big)
\]
The right hand side is $o(1)$ since $	H^{1-2\epsy}
\to 0$ for
$\epsy>1/2$ and $H\delta_{\thresh}^3\to0$ for
$\epsy<2/3$.  This proves the lemma.

 This establishes  \Cref{t:MDEsumApprox}
 with $o(1) =  O(     \thresh^{1-2\epsy}   \exp(   -  \half  \thresh^{-1+2\epsy}  /\gamma^2  )    )$ 
\qed

The next step is to approximate the sum in \Cref{t:MDEsumApprox}   by an integral.      For this we again consider the centered function $\tilG$
to simplify calculations.   		
\begin{lemma}
\label[lemma]{t:MDEsumApprox2int}   
			The following holds for any $\epsy \in (1/2,2/3)$:
\[
\frac{1}{\thresh} \sum   \exp\bigl(   -    \thresh  \tilG(s_n)  \bigr) 
			=  \int   \exp\bigl(   -    \thresh  \tilG(s)  \bigr)  \,  ds   +   O(\delta_{\thresh}^2)  
\]
			where
			the sum is over $n$ satisfying $ | s_n -s^*| \le \delta_{\thresh}$,  and the 
			integral is over $ | s -s^*| \le \delta_{\thresh} $.
\qed
\end{lemma}

\PF  
		The proof begins with the  standard Riemann sum approximation of an integral:
\[
\Big| \Delta \sum_{i=1}^n  g(t_i)  -   \int_a^b  g(s)\, ds   \Big|   \le  \half    \| g'\|_\infty (b-a)  \Delta   
\]
		with  $\Delta= (b-a)/n$,  $t_i = a + i\Delta$ for each $i$,  and $  \| g'\|_\infty$ the maximum of $|g'(s)|$ over the range of $s$ considered. 
		
		Here we take   $g(s) =  \exp (   -    \thresh  \tilG(s)   ) $.  In the integral considered in the lemma  we have $\Delta = 1/\thresh$ and $b-a =  2 \delta_{\thresh}$, giving the   error bound $  \| g'\|_\infty   \thresh^{-2+\epsy}$.     To complete the proof it remains to show that $ \| g'\|_\infty  \le O(  \thresh^{\epsy})$. 
		
		The derivative is   $g'(s) =   - \thresh G'(s)   \exp (   - \thresh  \tilG(s)   )   $.  
		In view of optimality we have $G'(s^*) =0$,  and hence 
\[
		|g'(s) |  = \thresh |  G'(s) - G'(s^*) |   \exp (   - \thresh  \tilG(s)   )   \le  \thresh \ell_{G'} |s - s^* |  
\]
		where $\ell_{G'}$ is a Lipschitz bound for $G'$ in a neighborhood of $s^*$.    Given the range of $s$ we arrive at the desired bound, 
		$|g'(s) |     \le   \ell_{G'}\thresh^\epsy   $.   
\qed
		
		We next obtain a computable approximation for the integral.    
\begin{lemma}
\label[lemma]{t:GintQuad}
			The following holds for any $\epsy \in (1/2,2/3)$:
\[
\thresh\int_{|s-s^*|\le\delta_{\thresh}}  
   \exp\bigl(   -    \thresh  \tilG(s)  \bigr)  \,  ds =   	[1 +   o(1) ]   \sqrt{\thresh}\sqrt{2\pi\gamma^2} \,,
\]
where  
			$\gamma^2 =  \Upupsilon_0''(\thexp_+)/\cm_+^3$ is defined in \Cref{t:barcApprox}.
\end{lemma}
		
\PF    We apply the approximation of $\tilG  $ appearing in  \eqref{e:tilGbdds} from which we obtain for $ | s -s^*| \le \delta_{\thresh}$,
\[
\begin{aligned}
   \exp\bigl(   -    \thresh  \tilG(s)  \bigr)   
&=   
   \exp\bigl(   -   \half \thresh    (s-s^*)^2/\gamma^2  \bigr)      \exp\bigl(       O(\thresh \delta_{\thresh}^3)   \bigr)   
 \end{aligned}
\]
We have  $ \exp\bigl(       O(\thresh \delta_{\thresh}^3)   \bigr)    = 1 +  O(\thresh \delta_{\thresh}^3) = 1+o(1) $ on noting that  $\thresh\delta_{\thresh}^3= \thresh^{-2+3\epsy}\to0$ as $\thresh\to\infty$ under the assumption that  $\epsy<2/3$.

		Consequently, writing $\sigma^2=\gamma^2/\thresh$,
\[
 \begin{aligned}
\thresh\int_{|s-s^*|\le\delta_{\thresh}}  & \exp\bigl(-\thresh\tilG(s)\bigr)\,ds
\\
&		=
		[1 +  o(1) ]  \thresh\int_{|s-s^*|\le\delta_{\thresh}}
\exp\bigl(-\tfrac{(s-s^*)^2}{2\sigma^2}\bigr)\,ds .
\end{aligned}
\]
		If $W\sim N(0,1)$, the last integral equals
\[
\sqrt{2\pi\sigma^2}\,\Prob\{|W|\le \delta_{\thresh}/\sigma\}.
\]
The probability above tends to one rapidly as $\thresh\to\infty$ since $\delta_{\thresh}/\sigma=\gamma^{-1}\thresh^{-1/2+\epsy}\to\infty$   under the assumption that  $\epsy>1/2$.  Therefore, as claimed, 
\[
\thresh\int_{|s-s^*|\le\delta_{\thresh}}  \exp\bigl(-\thresh\tilG(s)\bigr)\,ds
		=
			[1 +   o(1) ]   \sqrt{\thresh}\sqrt{2\pi\gamma^2}.
\]
 \qed

\subsection{Combining and optimizing the approximations}
\label{s:BigProof}		
 
 We arrive at the following approximations,
\begin{equation}
 \begin{aligned}
\MDD(\thresh)  &  = \MDDfine(\thresh) + o(1) 
\\
 \MDE(\thresh) & =
 [1+o(1)]   \MDEfine(\thresh) 
  \end{aligned}
\label{e:MDD-MDEapprox}
\end{equation}  
\begin{equation}
 \begin{aligned}
\textit{where} \ \ 
\MDDfine(\thresh) & \eqdef  \thresh/m_1 + \oshoot_\infty   
\\
\MDEfine(\thresh) 
&\eqdef
\sqrt{\thresh}\sqrt{2\pi\gamma^2}\exp(-\thresh\thexp_+)
 \end{aligned}
\label{e:barJinftyMDE}
\end{equation}
where the large-$\thresh$ approximation for $\MDD $ follows from \Cref{t:MDD_renewal} and
the approximation for $\MDE$ is obtained on combining its representation \eqref{e:MDEstep2} as a sum,
the approximation \eqref{e:preMDEsumApprox} for each term in this sum,
and the approximation for the sum obtained on combining  \Cref{t:MDEsumApprox,t:MDEsumApprox2int,t:GintQuad}.

Recall that $\barJ(\thresh, \kappa)  $ is the value  of  \eqref{e:MDD+kappaMDE} using CUSUM with threshold $\thresh>0$.  From the foregoing, 
\begin{equation}
 \barJ(\thresh,\kappa)    = \MDDfine(\thresh) + o(1) + [1+o(1)] \kappa \MDEfine(\thresh) 
\label{e:barJa}
\end{equation}
which motivates the notation  
\begin{equation}
\barJFineB(\thresh,\kappa)  = \MDDfine(\thresh) +\kappa \MDEfine(\thresh) 
\label{e:barJFineTwoTerms}
\end{equation} 
 The approximation $\barJFineB$ may be far from $\barJ$ when $\kappa \MDEfine(\thresh) $ is large.   Fortunately, on minimizing $\barJFineB$ over $\thresh$ we find that this quantity is bounded as a function of $\kappa$ (recall discussion preceding  \eqref{e:whyTight}).
 
The minimizer and minimizing value are denoted 
\[
 \barthreshFineB(\kappa)
    =\argmin_{\thresh\ge \thresh^0}
       \barJFineB(\thresh,\kappa) \,, 
\quad
\barJFineB(\kappa)
   =\min_{\thresh\ge\thresh^0}
       \barJFineB(\thresh,\kappa).
\]
where $ \thresh^0 >0$ is introduced to ensure a unique minimizer.  
\Cref{t:opt_f_kappa} below  implies that we may take $\thresh^0 = x^0/ \thexp_+$ with $x^0 \eqdef (1+\sqrt{2})/2$.

Approximation of these quantities is obtained through  the change of variables $x=\thexp_+\thresh$  and objective  
\[
		f_\kappa(x)=x+\kappa\exp(\bFine)\sqrt{x}e^{-x} \,,
\]
\begin{equation}
\begin{aligned}
\text{so that 
} \quad \barJFineB(\thresh,\kappa)  & =  \frac{1}{m_1\thexp_+} f_\kappa(x)+ \oshoot_\infty  
\\
\kappa \MDEfine(\thresh)   & =  \frac{1}{m_1\thexp_+}  [f_\kappa(x) - x]
\end{aligned}
\label{e:fJthresh}
\end{equation}
Computation of   $\barthreshFineB(\kappa)$  is thus reduced to minimizing $		f_\kappa$.

Any minimizer $\xkFine> x^0$ satisfies $f'_\kappa(\xkFine)=0$, which is equivalently expressed in two equivalent forms, where the second follows from the first on taking logarithms:   
\begin{subequations}
\begin{align}
	1 & =
	\kappa \exp(\bFine) \sqrt{\xkFine}e^{-\xkFine} 
	\Bigl[1-\frac{1}{2\xkFine}\Bigr] \,, 
	\label{e:fineFOC}
\\ 
0 & =
	\log(\kappa ) +  \bFine - \xkFine + \half \log(\xkFine ) + \log  
	\Bigl[1-\frac{1}{2\xkFine}\Bigr] \, .
	\label{e:fineFOClog}
\end{align}
		 		\end{subequations}

 \begin{lemma}
\label[lemma]{t:opt_f_kappa} 
$ f_\kappa$ is  strictly convex on $(x^0,\infty)$ for any $\kappa$, with $x^0 =(1+\sqrt2)/2$.  
  For sufficiently large $\kappa$ its minimum  $\xkFine$  is unique, and  admits the approximation $\xkFine =  \log\kappa+\frac12\log\log\kappa+ \bFine +  o(1) $.  
Moreover,  denoting $f_\infty^\Delta(a) = a+e^{-a}-1 $ for  fixed $a\in\Re$, we have
\begin{equation} 
\lim_{\kappa\to\infty} [
 f_\kappa(\xkFine+a)-f_\kappa(\xkFine) ]
 =
f_\infty^\Delta(a) 
\label{e:fkappaErr}
\end{equation} 
The function $f_\infty^\Delta$ is non-negative and strictly convex, with $a=0$ the unique minimizer.
\end{lemma}

\PF    
We have  $\xkFine= [1+o(1)] \log\kappa $ from
  \eqref{e:fineFOClog}, 
 and  the desired approximation is obtained 
on
 substitution  into \eqref{e:fineFOClog}:
\[
\xkFine
		=
\log\kappa+\frac12\log\xkFine+ \bFine +  o(1) 	=
\log\kappa+\frac12\log\log\kappa+ \bFine +  o(1) 
\]

A second derivative formula is easily obtained, 
\[
f_\kappa''(x)
=
\kappa e^\beta\sqrt{x}e^{-x}
\left(1-\frac1x-\frac{1}{4x^2}\right) 
\]
which establishes strict convexity of  $f_\kappa$   for
$x> x^0$.    Moreover,  for  sufficiently large $\kappa$ we have
$f_\kappa'(x^0)<0$, while
$f_\kappa'(x)\to1$ as $x\to\infty$.
Strict convexity therefore gives a unique minimizer
$\xkFine>x^0$ for such $\kappa$, characterized by \eqref{e:fineFOC}.

Finally, we obtain from \eqref{e:fineFOC}  $
\kappa e^\beta\sqrt{\xkFine}e^{-\xkFine}=1+o(1)$, from which we obtain \eqref{e:fkappaErr}: for  fixed  $a\in\Re$,
\[
\begin{aligned}
f_\kappa &(\xkFine+a)-f_\kappa(\xkFine)
\\
&=
a+\kappa e^\beta\sqrt{\xkFine}e^{-\xkFine}
\left[
\sqrt{1+\frac{a}{\xkFine}}e^{-a}-1
\right]
 \longrightarrow f_\infty^\Delta(a) 
\end{aligned}
\]
\qed

		\begin{subequations}
\begin{lemma}
\label[lemma]{t:fineApproxAlmost}
Under the assumptions of \Cref{t:MostLikelyPaths} we have the following approximations:
\begin{align} 
\barthreshFineStar(\kappa) & = \barthreshFineB(\kappa)  +o(1)\, ,
\label{e:barHFineApprox}
		\\
\barJFineStar(\kappa) & = \barJFineB(\kappa)    + o(1) \, .
		\label{e:barJFineApprox}
\end{align}
Moreover, both
$
\kappa \MDEfine( \barthreshFineB(\kappa)  ) $ 
and 
$
\kappa \MDEfine( \barthreshFineStar(\kappa)  ) $ 
are bounded in $\kappa$.
\end{lemma}

\label{e:fineApprox}
		\end{subequations}

\PF
We have by definition $
   \barthreshFineB(\kappa)=\xkFine/\thexp_+$
so \Cref{t:opt_f_kappa} gives \eqref{e:barHFineApprox}.
Moreover, \eqref{e:fineFOC} gives
\[
 \kappa e^{\bFine}\sqrt{\xkFine}e^{-\xkFine}
   =1+o(1).
\]
Hence, by \eqref{e:fJthresh},
\[
\begin{aligned}
 \barJFineB(\kappa)
   &=\frac{\xkFine+1+o(1)}{m_1\thexp_+}
      +\oshoot_\infty
    =\barJFineStar(\kappa)+o(1),
\\
 \kappa\MDEfine(\barthreshFineB(\kappa))
   &=\frac{1+o(1)}{m_1\thexp_+}.
\end{aligned}
\]
This proves \eqref{e:barJFineApprox} and  boundedness of $ \kappa\MDEfine(\barthreshFineB(\kappa))
$.

Boundedness of 
$
\kappa \MDEfine( \barthreshFineStar(\kappa)  ) $  then follows from \eqref{e:barHFineApprox}:
\[
\frac{\MDEfine(\barthreshFineStar(\kappa))}
     {\MDEfine(\barthreshFineB(\kappa))}
=1+o(1)
\]
\qed		


\def\threshB{h}

\wham{Proof of \Cref{t:barcApprox}.}
In view of \Cref{t:fineApproxAlmost}, to prove 
  \eqref{eq:refined_bdds} it is sufficient to establish the approximations, 
\begin{equation}
 \barJ^*(\kappa)
 =
 \barJFineB(\kappa)+o(1) \,,   \qquad
\barthresh^*(  \kappa) 
=
\barthreshFineB(\kappa)+   o( 1)	
\label{eq:refined_bdds_goal}
\end{equation}  
For this we streamline notation:
The respective optimizers of    $\barJ$ and $\barJFineB$ are denoted
\[
 \threshB^*=\barthresh^*(\kappa),
 \qquad
 \threshB^\circ=\barthreshFineB(\kappa).
\]

\wham{Step 1:}  \textit{$\threshB^*\to\infty$ as $\kappa\to\infty$}.
Since $\threshB^\circ\to\infty$ as $\kappa\to\infty$ and
$\kappa\MDEfine(\threshB^\circ)=O(1)$, \eqref{e:barJa} and
\Cref{t:fineApproxAlmost} give
\[
  \barJ(\threshB^\circ,\kappa)
   =\barJFineB(\kappa)+o(1)
   =O(\log\kappa).
\]
Optimality of $\threshB^*$ gives $  \barJ(\threshB^*,\kappa) \le
 \barJ(\threshB^\circ,\kappa)$, and therefore  
\[
 \MDE(\threshB^*)
 \le \frac{ \barJ(\threshB^*,\kappa)}{\kappa}
 =O\left(\frac{\log\kappa}{\kappa}\right)\longrightarrow0.
\]
Since $\MDE(\threshB)$ is non-increasing in $\threshB$, while
\eqref{e:MDD-MDEapprox} implies $\MDE(\threshB_0)>0$ for every
sufficiently large fixed $\threshB_0$. 

If a subsequence   $ \{\threshB^*_i =  \threshB^*(\kappa_i) : i\ge 1\}$ were bounded,
with $\kappa_i\to\infty$, fix $\threshB_0>0$
sufficiently large  exceeding this bound.  Monotonicity
would give $\MDE(\threshB^*_i)\ge\MDE(\threshB_0)>0$ for each $i$,
contradicting $\MDE(\threshB_i^*)\to0$.  Hence $\threshB^*\to\infty$.

\wham{Step 2:}  \textit{$\kappa\MDEfine(\threshB^*)$ is bounded}.
For   fixed $a>0$ we have
\[
\frac{\MDEfine(h^*+a)}{\MDEfine(h^*)}
=
\sqrt{1+\frac{a}{h^*}}\,
e^{-\thexp_+ a}
=
[  1+o(1)] e^{-\thexp_+ a}  
\]
 Optimality of $\threshB^*$ and
\eqref{e:MDD-MDEapprox} give
\[
\begin{aligned}
0
&\le
\barJ(\threshB^*+a,\kappa)-\barJ(\threshB^*,\kappa)
\\
&=
\frac{a}{m_1}
-
\bigl[1-e^{-\thexp_+a}+o(1)\bigr]
\kappa\MDEfine(\threshB^*)+o(1).
\end{aligned}
\]
Consequently $\kappa\MDEfine(\threshB^*)=O(1)$, as claimed.
This combined with \eqref{e:barJa} gives
\begin{equation}
\barJ^* (\kappa)   =
 \barJ(\threshB^*,\kappa)
   =\barJFineB( \threshB^*, \kappa)+o(1).
\label{e:Jalmost}
\end{equation}

\wham{Step 3:}  \textit{Cost approximation}.
Optimality of 
  $\threshB^*, \threshB^\circ$ combined with   \eqref{e:Jalmost}   implies the cost approximation in \eqref{eq:refined_bdds_goal}: 
\[
\begin{aligned}
\barJFineB(\kappa)
\le
\barJFineB(\threshB^*,\kappa)
& =
\barJ(\threshB^*,\kappa)+o(1)
\\
& \le
\barJ(\threshB^\circ,\kappa)+o(1)
\\
&
=
\barJFineB(\kappa)+o(1).
 \end{aligned}
\]

 \wham{Step 4:}  \textit{Threshold approximation}.
Let $x^*=\thexp_+ \threshB^*$.  The preceding   and
\eqref{e:fJthresh} imply
\[
  \tfrac{1}{m_1\thexp_+}
 \big[ f_\kappa(x^*)-f_\kappa(\xkFine) \big]
 =
  \barJFineB(\threshB^*,\kappa)   
-  \barJFineB(\threshB^\circ,\kappa)   =o(1). 
\]
Given $\epsy>0$, convexity gives, for all sufficiently
large $\kappa$,
\[
\inf_{\substack{x\ge x^0\\ |x-\xkFine|\ge\epsy}}
\bigl[f_\kappa(x)-f_\kappa(\xkFine)\bigr]
=
\min_{\sigma\in\{-1,1\}}
\bigl[
f_\kappa(\xkFine+\sigma\epsy)-f_\kappa(\xkFine)
\bigr].
\]
By \eqref{e:fkappaErr}, the right hand side converges to
\[
\min\{f_\infty^\Delta(-\epsy),
      f_\infty^\Delta(\epsy)\}>0.
\]
Since
$f_\kappa(x^*)-f_\kappa(\xkFine)=o(1)$, it follows that
$|x^*-\xkFine|<\epsy$ for all sufficiently large $\kappa$.
Hence $x^*-\xkFine=o(1)$, and therefore
$\threshB^*=\threshB^\circ+o(1)$.
 This completes the proof of \eqref{eq:refined_bdds_goal} and 
   \eqref{eq:refined_bdds}.
\qed

\subsection{Optimization and duality.}
\label{s:FstarInf}

We begin with justification of the   convex program  \eqref{e:jStarConvexRelaxation}.
\wham{Proof of \Cref{t:CoarseConvex}:}     
\Cref{t:Jscalings}   justifies the equivalence with the constrained optimization problem in which $\marg_1(F) $ 
is maximized subject to the equality constraint $ \Lambda_0(F) = \expa$ (so that $\thexp_+(F) =1$).  
 
Next we justify the relaxation  $ \Lambda_0(F) \le \expa$:   suppose that $F$ is feasible with $ \Lambda_0(F) < \expa$,
and denote  $F_r = F+ r\Fmix$.

For sufficiently small $r>0$, we can maintain $ \Lambda_0(F_r) < \expa$ and  
under the assumptions imposed on $\Fmix$ assumed in the proposition,
\whamrm{(i)}  $ \marg_0(F_r) =  \marg_0(F) + r \marg_0(\Fmix) <0 $ (so that     $F_r $  remains feasible)
\whamrm{(ii)}  The objective is increased:   $ \marg_1(F_r) =  \marg_1(F) + r \marg_1(\Fmix)   >  \marg_1(F) >0$.  
 
 Hence an optimizer  $F^*$ of  \eqref{e:jStarConvexRelaxation} will satisfy $ \Lambda_0(F^*) = \expa$.

Assuming an optimizer $F^*$  exists, observe that $ r\Fmix \in\clG$ is feasible for sufficiently small $r>0$  (to ensure  $ \Lambda_0(r\Fmix ) \le \expa$).    Optimality then implies positivity of the objective:   $\marg_1(F^* ) \ge  \marg_1(  r\Fmix )    > 0  $.  
\qed

We now pivot to analysis of the finite-dimensional optimization problem \eqref{e:barJ_infty_theta}.    
 Recall the   Lagrangian $ \clL(\theta,\zeta)  $ defined in  \eqref{e:LagrangeJFD}    and the dual function $\varphi^*(\zeta) = \inf \{  \clL (\theta,\zeta)  : \theta \in \Re^d \}  $  with domain $\Re_+^2$.  The following lemma addresses much of \Cref{t:DualProps}:
 
 \begin{lemma}
\label[lemma]{t:Lagrange}
Under \Ass{s:thetaCUSUM} we have for $\zeta \in \Re_+^2$,
\whamrm{(i)}
  $\varphi^*_r(\zeta)= -\infty$ whenever $\zeta_1 + \zeta_2 \neq 1$.
 
 \whamrm{(ii)}
 Otherwise  $\varphi^*(\zeta) = - \expa \zeta_1 +\varphi_0^*(\zeta) $, where $\varphi_0^*(\zeta) = $
  \[
\begin{aligned}
   \inf_\theta \Bigl\{    \zeta_1 [ \Lambda_0( F_\theta) - \marg_1 (F_\theta)   ] + (1-\zeta_1) [ \marg_0 (F_\theta) - \marg_1 (F_\theta) ]  \Big\}
  \end{aligned}
 \]

 \whamrm{(iii)} 
$\theta^\circ $ solves the minimization problem \eqref{e:DualAndGamma} if and only if
$\nabla \Gamma_0(\theta^\circ) =0$, which holds 
 if and only if the moment matching constraints 
\eqref{e:calcOptLinDual} are satisfied.
\end{lemma}

 \PF 
 Part (ii) is immediate from the definition of the Lagrangian.
 
 To see (i) recall from Assumption \Ass{s:thetaCUSUM}~(iv)   that there is $\thone \in\Re^d$ satisfying  ${\thone}^\transpose \psi\equiv 1$.
 For any $t\in\Re$ we have  $ \Lambda_0( F_\theta) = t$ when  $\theta = t \thone$.     Consequently,  if   $\zeta_1 + \zeta_2 \neq 1$,  
\[
\varphi^*(\zeta) \le  \inf_t  \clL(t \thone,\zeta)  = 
  - \expa \zeta_1  + 
 \inf_t \{  - t + \zeta_1 t + \zeta_2 t\} = -\infty
 \]
 
 Part (iii) follows from    
 \eqref{e:twistedMean_psi}:   
$
 \nabla \Gamma_0(\theta) =  \twmarg_0^{\theta} (\psi) -\marg_1 (\psi) $.   \qed

\begin{lemma}
\label[lemma]{t:pi0Fneg}
Suppose that Assumption~\Ass{s:thetaCUSUM}  holds, and that $\theta^\circ$ minimizes $\Gamma_0$ and satisfies $\Lambda_0(  F_{\theta^\circ}) =0$.    
Then 
$ \marg_1(F_{\theta^\circ}) >0$ and 
 $ \marg_0(F_{\theta^\circ}) = -  \clK(\trajP_0\|\twtrajP_0^{\theta^\circ})   < 0$.   
 \end{lemma}

 \PF    
We have  $- \marg_1(F_{\theta^\circ})  = \Gamma_0  (F_{\theta^\circ}) \le  \Gamma_0  (\theta) $ for any $\theta$ by optimality of $\theta^\circ$ and the assumption $ \Lambda_0(  F_{\theta^\circ}) =0$.    
Hence to establish $ \marg_1(F_{\theta^\circ}) >0$ it is sufficient to show that $ \Gamma_0  (\theta) <0$ for some $\theta$. 

Recall $ \thmix $ is defined in Assumption~\Ass{s:thetaCUSUM},
satisfying $ \marg_0(F_{\thmix} ) < 0$ and $ \marg_1(F_{\thmix} ) >0$.
For sufficiently small $r>0$,  
\[
 \Gamma_0  (r \thmix)  =  r \{  \marg_0(F_{\thmix} )  -  \marg_1(F_{\thmix} ) \}  + O(r^2) <0 \, , 
\]
which establishes $\marg_1(F_{\theta^\circ}) >0$.

Consideration of  $ \marg_0(F_{\theta^\circ}) $ is based on   the identity,
\[
 \Lambda_0(  F_{\theta^\circ})  - \marg_0(F_{\theta^\circ})  
 =
	\clK(\trajP_0\|\twtrajP_0^{\theta^\circ}) 
\]
which is obtained from \Cref{t:m1_ent_F} on replacing
$\trajP_1$ with   $\trajP_0$  in
	\eqref{e:m1_ent_F}.
  It remains to show that 
 $\clK(\trajP_0\|\twtrajP_0^{\theta^\circ})  >0 $.   
 For this we apply the moment matching conditions.

    \Cref{t:Lagrange}~(iii) combined with  \eqref{e:twistedMean_psi} imply that the optimizer  $\theta^\circ$ satisfies,
  \[
   \tfrac{d}{dr}  \Lambda_0( r F_{\theta^\circ}) \big|_{r=1} =  \twmarg_0^{\theta^\circ} ( F_{\theta^\circ})  =    \marg_1(F_{\theta^\circ})   >0 
 \]
  and applying 	\eqref{e:twistedMean},
  \[
   \tfrac{d}{dr}  \Lambda_0( r F_{\theta^\circ}) \big|_{r=0} =  \marg_0 ( F_{\theta^\circ})   
  \]   
Suppose that in fact   $ \marg_0(F_{\theta^\circ})  = 0$.     Then the convex function $f(r) = \Lambda_0( r F_{\theta^\circ})    $
is minimized at $r=0$, with $f(0) =0$, and is hence   everywhere non-negative.    The value $r=1$ is also a minimizer, since $f(1) =   \Lambda_0(  F_{\theta^\circ})  = 0 $ by assumption.    
   Convexity implies   that   $\Lambda_0( r F_{\theta^\circ})  = r  \marg_0(F_{\theta^\circ})  =0 $  for all $0\le r\le 1$, and thus
\[
0 = \tfrac{d}{dr}  \Lambda_0( r F_{\theta^\circ}) \big|_{r=1} =    \marg_1(F_{\theta^\circ})   
\]
  This contradicts     $ \marg_1(F_{\theta^\circ})  >0$, completing the proof.
           \qed

\wham{Proof of \Cref{t:DualProps}:}     
The form of the dual \eqref{e:DualAndGamma}  is standard. 
 Moreover, \eqref{e:GammaKtheta} implies that 
 $\Gamma_0( \theta) \ge      -
	\clK(\trajP_1\|\trajP_0) $, 	 so that $ \varphi^*(\zeta^*)\ge    -   \expa  - 	\clK(\trajP_1\|\trajP_0)$  with $\zeta^* = (1;0)$.  
  This combined with \Cref{t:Lagrange} establishes (i).
  
Part (ii) was established in 	\eqref{e:calcOptLinDual}.
    
Part (iii) concerns the parameter $\theta^* = \theta^\circ + r^*\thone $,  where $r^* = \expa - \Lambda_0( F_{\theta^\circ})$ is chosen so that $\Lambda_0( F_{\theta^*})=\expa$.   

We consider first $\expa=0$  and denote $\theta^\bullet =  \theta^\circ + r^\bullet \thone $
with $r^\bullet= - \Lambda_0( F_{\theta^\circ})$, so that $ \Lambda_0( F_{\theta^\bullet}) =0$.    
\Cref{t:pi0Fneg} then gives the desired conclusions:
	$ \marg_1(F_{\theta^\bullet}) >0$ and 
 $ \marg_0(F_{\theta^\bullet}) = - \expa^0   < 0$, with $\expa^0 = \clK(\trajP_0\|\twtrajP_0^{\theta^\bullet}) $. 

In general, with $r^* = \expa  + r^\bullet $ and $\theta^* = \theta^\circ + r^*\thone $,
\[
   \begin{aligned}
 \marg_1(F_{\theta^*})    &= \expa + \marg_1(F_{\theta^\bullet}) >0 
\\
  \marg_0(F_{\theta^*})  & = \expa - \expa^0 \le  0\,, \ \ 0\le \expa \le \expa^0
  \end{aligned}
 \]
We conclude that $\theta^*$      is feasible and satisfies complementary slackness with $\zeta^* = (1;0)$, which implies optimality.   
   \qed

\wham{Proof of \Cref{t:scalarOffsetOpt}}     

Consider the basis $\psi=(G;1)$, so that
$F_\theta=\theta_1G+\theta_2$.
By \Cref{t:DualProps}, an optimizer   is $ F^\bullet=\thexp_+G+r^\circ$ with 
$ r^\circ=\expa-\Lambda_0(\thexp_+G)$ and 
\[
 \thexp_+\in\argmin  \{  \Lambda_0( \thexp G ) -  \thexp \marg_1(G) \} 
\]
Applying \Cref{t:Jscalings} we conclude that $F^* =  F^\bullet /\thexp_+$ 
also minimizes $\barJStarCoarse$ over this two dimensional
function class.   This establishes \eqref{e:OptOffset} since by construction  $ F^* = G+\theta^*$  with 
$ \theta^*= [\expa-\Lambda_0(\thexp_+G)]/ \thexp_+$.

Turning to \eqref{e:ScalarOffsetConsequences}, the first order condition for optimality of $\thexp_+$ in the maximization \eqref{e:OptOffset} may be expressed 
$   \tfrac{d}{d \thexp} \Lambda_0 (\thexp F)    =    \marg_1(F) $ when $\thexp = \thexp_+(F)$.    
Writing  $\Lambda_0 (\thexp F) = \Lambda_0 (\thexp G) 
 + \thexp \theta^*$ we conclude that $   \tfrac{d}{d \thexp} \Lambda_0 (\thexp F)    =    \marg_1(F) $ at $\thexp=\thexp_+$  as claimed.

The expression $   \Lambda_0(\thexp_+ F)   =  \expa $ is by definition of $\theta^*$.
\qed

We now return to the possibly infinite dimensional   linear function class $\clG$.

\wham{Proof of \Cref{t:jStarConvexRelaxationENT}:}     

The identity \eqref{e:GammaK} implies that $\clK(\trajP_1\|\twtrajP_0^{F+r})  = \clK(\trajP_1\|\twtrajP_0^F) $ for any $r\in\Re$ and any $F\in\clG$.  This justifies the claim that   $\Lambda_0(F^\circ) = 0$ can be assumed without loss of generality. 

Under the assumptions of the proposition, the remaining conclusions are obtained by reducing $\clG$ to  the three dimensional function class spanned by $\{1, F^\circ , \Fmix \}$, and applying \Cref{t:DualProps}.
\qed

\wham{Proof of \Cref{t:completeOpt}:}     
For any $G\in \clG$ consider
the three dimensional function class $\clG_0$ spanned by $\{1, F^\circ, G \}$.
Then $\theta^* = (0; 1; 0)$ is an optimizer of $\Gamma_0$, and applying \Cref{t:DualProps} we conclude that 
\[
\twmarg_0^{F^*} (G) = 
 \twmarg_0^{\theta^*} (G) = \marg_1 (G)    
\]
 This establishes that
  $\twmarg_0^{F^*}  = \marg_1$ under the assumption that $\clG$ is dense.
  \qed

We conclude this subsection with a simple result required in the proof  of \Cref{t:Model2clG}:

For a function $H\colon\ystate\to\Re$, define
\[
H_+(y,y')=H(y'),
\qquad
H_-(y,y')=H(y).
\]

\begin{lemma}
	\label[lemma]{t:Null2}
	Suppose that $\del=2$.  Let $F\colon\ystate^2\to\Re$ and let
	$H\colon\ystate\to\Re$ be bounded.  Suppose that the linear function
	class spanned by $\{F,H_+,H_-,1\}$ satisfies \Ass{s:thetaCUSUM}.  Then
	\[
	\Lambda_0(F+H_+-H_-) = \Lambda_0(F).
	\]
\end{lemma}

\PF 
Letting $G = F+H_+-H_-$  we find that 
 the difference between the corresponding partial sums in
\eqref{e:logCGF} telescope:   with $\del=2$,
\[
\Big(  \sum_{k=0}^{n-1} G(\strY{\del}_k) \Big)  -  \Big(  \sum_{k=0}^{n-1} F(\strY{\del}_k) \Big)
=
	H(\preObs_{n-1})-H(\preObs_0).
\]
This difference is bounded in absolute value by
$2\|H\|_\infty$, so  $\Lambda_0(G) = \Lambda_0(F) $  from the definition \eqref{e:logCGF}.
\qed

\wham{Proof of \Cref{t:Model2clG}:}
We have  $
F_{\thzz}(y,y') = {\thzz}^\transpose\psi(y,y') = y'-y$,  for any $y,y'\in\ystate$,
by construction of the vector $\thzz$ and the basis $\psi$.
\Cref{t:Null2} gives $\Lambda_0(F_{\theta+s\thzz}) = \Lambda_0(F_\theta)$ for any $s$,  and   stationarity gives  $\marg_i^{[2]} (F_{\theta+s\thzz}) =\marg_i^{[2]} (F_{\theta})$ for each $i$.  From the definitions we conclude that  
\[
\Gamma_0(\theta+s\thzz)=\Gamma_0(\theta)\, , \quad s\in\Re
\]
The proof is completed on recalling that   $F_{\theta+r\thone} = F_{\theta} + r$ and 
hence
 $\Gamma_0(\theta+r\thone)= \Gamma_0(\theta) $ for $r\in\Re$.    
\qed

\subsection{Proof of \Cref{t:calcOpt}}

We proof a generalization of  \Cref{t:calcOpt} based on a finite-dimensional function class 
$\clG = \{ F_\theta  : \theta\in\Re^d \}$,  with   $F_\theta$ is continuously differentiable in $\theta$.  
We denote
  $\psi_\theta =\nabla_\theta F_\theta$ and  $\Upupsilon_0^\theta(\thexp) = \Lambda_0(\thexp F_\theta )$ for each $\theta\in\Theta$ and $\thexp\in\Re$.   Let $  \twmarg_0^{\thexp ,\theta}$ denote the probability measure
\eqref{e:twmarg0} obtained using $\thexp $ with function $F_\theta$.      When using $\thexp_+^\theta$ we simplify the notation to
		$  \twmarg_0^{ +,\theta}$,  and write $\cm_0^{+,\theta} =  \twmarg_0^{ +,\theta} (F_\theta)$.

The following assumptions are analogous to those of Assumption~\Ass{s:thetaCUSUM}:

  \wham{Assumptions   for the nonlinear function class.}
For each $\theta\in\Re^d$,  
\whamrm{(i)} 
  The mean $m_\theta^i = \marg_i(F)$  exists for each $i$,

  \whamrm{(ii)} 
  If  $\Fmax^0> 0$ 
  with $F=F_\theta$ 
  and 
  $m_\theta^0 \le  0$    then parts (iii)--(v) of  Assumption~\Ass{s:extend} hold for this $F$.

  \whamrm{(iii)} 
    If    $m_\theta^1 >  0$    then   \eqref{e:oshootA4} holds for some $ \oshoot_\infty \in\Re$.

  \whamrm{(iv)}   There is  $\thone \in\Re^d$ such that $  \psi_\theta \transpose \thone \equiv 1$,  and   $\theta^\diamond \in\Re^d$ such that   
$ \marg^0(F_{\theta^\diamond} ) < 0$ and $ \marg_1(F_{\theta^\diamond} ) >0$.

  \smallskip
  
  In addition,
    \whamrm{(v)}  
   The functions $\Upupsilon_0^\theta (\thexp) $, 
    $\Lambda_0 ^\theta (\thexp \psi_\theta^i) $ and $ \twmarg_0^{\thexp ,\theta}  (  \psi_\theta^i) $ are each continuous over  
    $(\theta,\thexp)\in   \Theta \times \Re$ for each $1\le i\le d$.    
       
  $\psi_\theta(y)=\nabla_\theta F_\theta(y)$ exists and is continuous in $\theta$ for each $y\in\ystate$.
\qed

\begin{proposition}
\label[proposition]{t:calcOptNon}
Suppose that     $\theta^\bullet\in\Theta$ is a stationary point of \eqref{e:barJ_infty_theta}:   $ \nabla_\theta  \barJStarCoarse \, (   \kappa; \theta^\bullet ) = 0$ for some (and hence all)  $\kappa>0$.  Suppose moreover that    $\Fmax^0>0$ 
and $ \marg_0(F)\le 0$  
 with $F= F_{\theta^\bullet}$.   Then,
\begin{equation}
	\cm_0^{+,\theta^\bullet}  =m_1^{\theta^\bullet}
	\ \ \textit{and} \ \ 
	\twmarg^0_{+,\theta^\bullet} \big( \psi_{\theta^\bullet} \big)     =   \marg_1\big( \psi_{\theta^\bullet} \big)
	\label{e:calcOptNon}
\end{equation}
			\qed
\end{proposition}

The proof requires interpretation of the gradient,
   $ \nabla_\theta  \log
			\barJ_\infty^* (     \kappa; \theta)
			=
			-     \nabla_\theta  \log    m_1^\theta   -   \nabla_\theta  \log \thexp_+^\theta$.     
\Cref{t:twistedMean_psi} provides a first step, whose proof is   immediate from the assumptions: 		
 
\begin{lemma}
\label[lemma]{t:twistedMean_psi}
Under the assumptions of this subsection, consider  any $\theta$ for which
     $\Fmax^0>0$ 
and $ \marg_0(F)\le 0$  
 with $F= F_{\theta}$.
Then,
\[
\frac{d}{d\thexp}
\Upupsilon^\theta_0  \,  (\thexp)  =   \twmarg_0^{\thexp ,\theta} (F_\theta)
\qquad 
\nabla_\theta
\Upupsilon^\theta_0  \,  (\thexp)  =    \twmarg_0^{\thexp ,\theta} (\psi_\theta)
\]			\clThm
\end{lemma}

 We next obtain the components of the gradient  $ \nabla_\theta  \log
			\barJ_\infty^* $.

\begin{lemma}
	\label[lemma]{t:calc1}
For $\theta$ satisfying the assumptions of \Cref{t:twistedMean_psi} we have
	\begin{equation}
		\begin{aligned}
		 \nabla_\theta  \log    m_1^\theta   = \frac{1}{m_1^\theta} \marg_1\big( \psi_{\theta} \big)  
\qquad
 \nabla_\theta  \log \thexp_+^\theta	 =  -  
			\frac{1}{\cm_0^{+,\theta}} \twmarg_0^{ +,\theta} \big( \psi_{\theta} \big)
		\end{aligned}
		\label{e:barJthetaoptGrad}
	\end{equation}
\end{lemma}

\PF 
The first identity  follows from the definitions. 
To obtain the desired representation for $ \nabla_\theta  \log \thexp_+^\theta	$  we use the defining property that  $\Upupsilon^\theta(\thexp_+^{\theta})  = \expa$ for all $\theta\in\Theta$.
		Applying  \Cref{t:twistedMean_psi}
		gives   for each $\theta\in\Theta$,
\[
		0 = \nabla_\theta  \Upupsilon^\theta(\thexp_+^{\theta})     = \nabla_\theta  \Lambda_0(   \thexp_+^{\theta}  F_\theta ) = 
\twmarg_0^{+,\theta} \bigl(   \nabla_\theta  \{  \thexp_+^{\theta}  F_\theta \} \bigr)
\]
		Applying the product rule,
\[
		0 =   \cm_0^{+,\theta}  \nabla_\theta   \thexp_+^{\theta}    +   \thexp_+^{\theta}    \twmarg_0^{+,\theta} \big( \psi_{\theta} \big)   
\]  
which implies the desired representation for $ \nabla_\theta  \log \thexp_+^\theta	$.
\qed

\wham{Proof of \Cref{t:calcOpt}:}
If $\theta^\bullet\in\Re^d$ satisfies the assumptions of the proposition then 
\Cref{t:calc1}   implies the following for any vector $u\in\Re^d$:
\[
\begin{aligned}
\frac{1}{m_1} \marg_1\big(u^\transpose  \psi_{\theta^\bullet} \big)  =  
\frac{1}{\cm_0^{+,\theta^\bullet}} \twmarg_0^{+,\theta^\bullet} \big( u^\transpose  \psi_{\theta^\bullet} \big)  
\end{aligned} 
\]
On choosing $u=\thone$ this gives  $u^\transpose  \psi_{\theta^\bullet}  \equiv 1$ and hence $m_1 =  \cm_0^{+,\theta^\bullet}$.   
Since $u$ is arbitrary we conclude that $	\twmarg^0_{+,\theta^\bullet} \big( \psi_{\theta^\bullet} \big)     =   \marg_1\big( \psi_{\theta^\bullet} \big)
$, and hence   the desired conclusion \eqref{e:calcOptNon}.
\qed

\subsection{POMDPs and quasi-stationarity}
\label{s:LDRRW} 
		
Here we provide a proof of \Cref{t:POMDPhazard}.

We begin with the required P-F theory for $\haP$ appearing in \eqref{e:twisted}.  It is simplest to adopt matrix notation,    $\bigstate = \{1,\dots, N\}$,  and $\bigstate_0 = \{ 1,\dots, N_0 \}$ with $1\le N_0<N$.   
Hence the matrix
$\haP$  is the $N_0\times N_0$ matrix with $\haP(i,j) = P(i,j)$ for $1\le i,j \le N_0$.      The following can be found in \cite{num84}:
		
\begin{lemma}
\label[lemma]{t:PFlemma}
			The following hold under the assumptions of \Cref{t:POMDPhazard}:
\whamrm{1)}
There is a solution to the eigenvector equation $\haP \xi = \PFeig \xi$, with  $\PFeig>0$   and  $\xi_i > 0$ for each $i \in\bigstate_0$.

\whamrm{2)}   
The \textit{twisted transition matrix}   in \eqref{e:twisted} becomes   $\cP_{i,j} = \PFeig^{-1} \xi_i^{-1} \haP_{i,j}   \xi_j  :  i,j\in\bigstate_0$.

\whamrm{3)}    
The transition matrix $\cP$    defines a uni-chain and aperiodic Markov chain with state space $\bigstate_0$.  Consequently it  has a unique invariant pmf $\ometamarg$, and 
for each $i,j\in\bigstate_0$,
\[
\ometamarg(j) = \lim_{n\to\infty}  \cP^n_{i,j} =  \frac{\xi_j }{  \xi_i }  \lim_{n\to\infty} \frac{1}{ \PFeig^n}  
P_{i,j} ^n
\]   
			where convergence holds at a geometric rate.  
\end{lemma}

\PF 
		Assumption~\Ass{s:POMDP}~(iii) implies a weak form of   irreducibility and aperiodicity:  for some $i^\circ$,
\[ 
		\haP^n_{i, i^\circ}  >0 \,,\qquad n\ge n^\circ 
\]
This is sufficient to ensure the existence of the P-F eigenvalue-eigenvector pair  \cite{konmey03a}.		
It follows that  $\cP^n_{i, i^\circ}
  >0 $
  for all  $ n\ge n^\circ$ and $i \in \bigstate_0$, so that this transition matrix  defines a uni-chain  and aperiodic Markov chain.     \qed

\wham{Proof of \Cref{t:POMDPhazard}}  
Part (i) is immediate from \Ass{s:POMDP}~(ii) which ensures ergodicity of $\bfPhi$.    
The multiplicative ergodic theorem \eqref{e:logCGF1} may be found in \cite{konmey03a}.

The transition matrix  $\oP$ defined in \eqref{e:twisted} is precisely $\cP$ in \Cref{t:PFlemma}.   From the lemma we know it is uni-chain and aperiodic, with unique invariant pmf on $\bigstate_0$, denoted $\ometamarg$.  

The $n$-fold transition matrix admits the interpretation,  
\begin{equation}
\oP^n(\bst,\bst')  =    \frac{1}{\PFeig^n} \frac{\xi(\bst')}{\xi(\bst)}  \Expect[  \ind \{ \tchange>n  \} \ind \{ \Phi_n = \bst'   \} 
			 \mid  \Phi_0 = \bst]
\label{e:twisted-n}
\end{equation}
\begin{equation}
 \begin{aligned}
 \! \! \!
 \textit{giving
}
\
\Prob \{ \tchange>n  \}  &=    \xi(\bst) \Big[\sum_{\bst'} \oP^n(\bst,\bst')      \frac{1}{\xi(\bst')}   \Big] \PFeig^n 
\\
&=  \xi(\bst)   [ \ometamarg(1/\xi)  + \epsy_n]    \PFeig^n 
 \end{aligned}
\label{e:twisted-n_tchange}
\end{equation}
where the final equality follows from ergodicity of $\bfoPhi$. 
Part (iii) easily follows with   $\expa = - \log(\PFeig)$.

 The proof of (ii) follows from an interpretation of 
the probability $\Prob\{  \oPhi_{n-k+1} = \bst_k \, , 1\le k\le m \}$ similar to \eqref{e:twisted-n},
which converges to $\Prob\{  \oPhi^\infty_{n-k+1} = \bst_k \, , 1\le k\le m \}$ at a geometric rate as $n\to\infty$.   
 We provide a proof only for $m=1$:    Applying \eqref{e:twisted-n} once again we have for any initial condition $\Phi_0=\bst \in\bigstate_0$, 
 \[
 \begin{aligned}
\Prob & \{  \Phi_{n-\ell}  = \bst_1  \,,   \tchange >   n  \}  
 \\
 &= 
 P^{n-\ell}(\bst, \bst_1 ) P^\ell (\bst_1 ,  \bigstate_0) 
 \\
& =
 \Big(
  \PFeig^{n-\ell}  \xi(\bst)   \frac{1} {\xi(\bst_1)} \oP^{n-\ell} (\bst,\bst_1)    
  \Big)
  \\
  & \ \ 
  \times
  \Big( 
   \xi(\bst_1) \Big[\sum_{\bst'} \oP^\ell(\bst_1,\bst')      \frac{1}{\xi(\bst')}   \Big] \PFeig^\ell   \Big)
   \\
  & = \PFeig^n   \xi(\bst)   \Big(   \ometamarg ( \bst_1) + \epsy_{n-\ell}(\bst) \Bigr)    \Big( \ometamarg(1/\xi)    + \epsy_{\ell}(\bst_1) \Bigr) 
  \end{aligned}
\]
where $\epsy_{\ell}(\bst_1) \to 0$ at a geometric rate as $\ell\to\infty$, and $\epsy_{n-\ell}(\bst) \to 0$ 
 a geometric rate as $n\to \infty$ for   fixed $\ell$.

Substituting \eqref{e:twisted-n_tchange} then gives $\Prob \{  \Phi_{n-\ell}  = \bst_1  \,,   \tchange >   n  \}  =$
 \[
 \begin{aligned} 
   \Big(   \ometamarg ( \bst_1) + \epsy_{n-\ell}(\bst) \Bigr)    \Big( \Prob \{ \tchange>n  \} [ 1+   \epsy_{\ell}(\bst_1)     \xi(\bst) /\ometamarg(1/\xi)  ]   \Bigr) 
  \end{aligned}  
\]
This implies \eqref{e:POMDPergodic_ii} when $\ometamarg( \bst_1)>0$.

The derivation of the expression for  $  \Lambda_0(G) $ in \eqref{e:logCGF0} proceeds as follows.    
Assume  that $\del=1$ for simplicity, and consider the positive matrix
$\haP_g(\bst,\bst') =  \exp(G(h(\bst))) \oP(\bst,\bst') $ with P-F eigenvalue-eigenvector pair $\lambda_g, \cg $, so that as in the preceding we may define a transition matrix on  $\bigstate_0$ via $\cP_{\!g} (\bst,\bst')  =[ \lambda_g \cg (\bst)]^{-1} \cg (\bst')  \haP_g(\bst,\bst') $.   

The remainder of the proof is  based on interpretation of the $n$-fold product $\cP_{\!g}^n$, similar to  \eqref{e:twisted-n}, and applying the ergodic theorem.

\begin{subequations}

Letting $\ometamarg_g$ denote the invariant pmf for $\cP_{\!g}$, we have for any   $f\colon\bigstate_0\to\Re$,
\[
\lim_{n\to\infty} 
\cP_{\!g}^n f\, (\bst)  = \ometamarg_g(f) \,,   \quad \bst \in \bigstate_0 \, .
\]
We next interpret the left hand side:
With $\oX_n =  h (\oPhi_n) $, and all expectations conditioned on $ \Phi_0 = \bst$,
\begin{equation}
\begin{aligned}
\cP_{\!g}^n f\, (\bst)   = \Expect[f(\oPhi_n)  \mid \oPhi_0 = \bst] 
	 &= \frac{1}{\lambda_g^n}  \frac{1}{\cg (\bst)}  \Expect\Big[ \exp\Big(  \sum_{k=0}^{n-1} G( \oX_k   ) \Big)   \cg ( \oPhi_n ) f( \oPhi_n ) \Big]
\\
&=    \frac{1}{\PFeig^n \lambda_g^n}  \frac{1 }{ \xi(\bst) \cg (\bst)}  \Expect\Big[ \exp\Big(  \sum_{k=0}^{n-1} G( Y_k   ) \Big)   \cg ( \Phi_n ) \xi( \Phi_n )  f( \Phi_n )   \ind \{ \tchange>n  \}\Big]
\end{aligned}
\label{e:cPgn}
\end{equation}
This combined with the ergodic theorem gives, 
\begin{align}  
\ometamarg_g  ( 1/ \cg   ) & = \frac{1}{\cg (\bst) }   \lim_{n\to\infty}   \frac{1}{\lambda_g^n}  \Expect\Big[ \exp\Big(  \sum_{k=0}^{n-1} G( \oX_k   ) \Big)    \Big]   &&    f =  1/ \cg 
\label{e:metaCDFo}
\\
\ometamarg_g  (  1/ [\cg \xi]  ) 
   &=  \frac{1 }{ \xi(\bst) \cg (\bst)} 
  \lim_{n\to\infty}  \frac{1}{\PFeig^n \lambda_g^n}   \Expect\Big[ \exp\Big(  \sum_{k=0}^{n-1} G(  Y_k   ) \Big)  
			  \ind \{ \tchange>n  \}\Big]     &&  f =  1/ [\cg \xi] 
\label{e:metaCDF}
 \end{align}
\label{e:meta2f}
\end{subequations}

On taking logarithms,    \eqref{e:metaCDFo} gives \eqref{e:logCGF0twist} with $\Lambda_0(G) = \log(\lambda_g)$:
\[
\lim_{n\to\infty} \frac{1}{n} \log  \Expect\Big[ \exp\Big(  \sum_{k=0}^{n-1} G( \oX_k   ) \Big)    \Big]  
				= \log(\lambda_g)
\] 
On combining \eqref{e:metaCDF}
with   
\eqref{e:twisted-n_tchange} we obtain
\begin{equation*}
\ometamarg_g( 1/ [\cg \xi] )  =  \frac{    \ometamarg(1/\xi) }{ \cg (\bst)}   \lim_{n\to\infty}  \frac{ 1  }{  \lambda_g^n}  \Expect\Big[ \exp\Big(  \sum_{k=0}^{n-1} G(  Y_k   ) \Big)  
			\mid \tchange>n \Big] 
 \end{equation*}
This implies  \eqref{e:logCGF_POMDP0} with the same limit, $\Lambda_0(G) = \log(\lambda_g)$.
 \qed

\subsection{Experiment details}
\label{s:appndx_exp}

This section records the simulation settings used for \Cref{s:numres}.  Throughout, $N=6\times10^6$ denotes the number of independent runs used for each threshold search or policy evaluation.  A run is terminated when the stopping rule fires; if no stop has occurred by time $\kmax $, the run is truncated at $\kmax $.   
Over the $N$ independent runs we obtain estimates of the CUSUM* threshold
$\barthreshFineStar(\kappa)$ and associated cost $\barJFineStar(\kappa)$ for a selected range of $\kappa$.   
For each $\kappa$ the value $\barthreshFineStar(\kappa)$ was obtained via grid search.   The grid and the value of $\kmax$  varied depending on the test: 
\[
\begin{array}{c|c|c|c}
\text{Test}
& \text{$\thresh$ grid range}
& \text{grid increment}
& \kmax\\
\hline
\rule{0pt}{3.2ex}
1\text{a--}1\text{c}
& [2.5,9.0]
& 2\times10^{-2}
& 2\times10^3\\
2\text{c--}2\text{d}
& [3,70]
& 2\times10^{-3}
& 2\times10^4\\
2\text{e}
& [2.09,9.0]
& 5\times10^{-4}
& 2\times10^3
\end{array}
\]

 Tests~a and b used $\kmax =  2\times10^4$ in Monte-Carlo estimates of cost   \eqref{e:MDD+kappaMDE},
  based on the respective  decision regions   $\clR^*$,  $\clR_{\threshSmallStar}$.

We also required estimates of  
$\oshoot_\infty$ defined in \eqref{e:oshootA4} and appearing in the definition of $\aFine$ in \Cref{t:barcApprox}.  
For a range of thresholds $\thresh \ge 1$   estimates of 
 $\oshoot(\thresh) $ defined in  \eqref{e:oshootDef} were obtained via   Monte-Carlo.  Independent sample paths were generated with $\tchange =0$   to obtain the   sample mean estimate $\widehat{\MDD}(\thresh)$, from which we define   
\begin{equation}
\widehat{\oshoot}(\thresh) \eqdef \widehat{\MDD}(\thresh)-\frac{\thresh}{m_1}.
\label{e:oshootEst}
\end{equation}
Repeating the computation over increasing thresholds provided evidence of convergence and hence a numerical proxy for $\oshoot_\infty$.
\notes{sm2ac:  note simplifications here}

\Wham{Model 1}
The observations are defined by the C.I.\ Markov model,  with linear Gaussian dynamics \eqref{e:num_m1_ar},   $A^0=0.30$, $A^1=0.60$, $\sigma=1$, and $\expa=0.02$.  The corresponding stationary variances in \eqref{e:num_m1_ar} are
\[
\frac{\sigma^2}{1-(A^0)^2}\approx1.0989,
\qquad
\frac{\sigma^2}{1-(A^1)^2}=1.5625.
\]
For each of  the three tests  and each threshold $\thresh$, sample averages of $\MDD(\thresh)$ and $\MDE(\thresh)$ were stored.  These averages determine the empirical estimate of \eqref{e:MDD+kappaMDE}, and CUSUM* was obtained by minimizing this estimate over a uniform threshold grid on $[2.5,9.0]$ with increment $0.02$.

The constants used in \eqref{e:barHFineBody}--\eqref{e:barJFineBody} are included in the following table:
\[
\begin{array}{c|c|c|c|c|c|c|c}
\!\!
\text{Test}\backslash \text{par}
 & r^*  & \thexp_+ & m_0 & m_1
	& \gamma^2 & \bFine & \oshoot_\infty \\
	\hline
	\rule{0pt}{3.2ex}
	\text{1a}
	& 0.020 & 1.00 & -0.03 & 0.09
	& 327 & 1.41 & -1.64\\
	\text{1b}
	& 0.024 & 0.69 & -0.03 & 0.10
	& 251& 1.21 & -9.12\\
	\text{1c}
	& 0.022 & 0.83 & -0.03 & 0.10
	& 257& 1.30 & -5.08
\end{array}
\]
\Cref{t:bestF}  tells us that $\thexp_+ =1$ and $ r^*  = \expa = 0.020$ for   Test~1a.     For Laplace and Student's $t$, write $G=\surL$ and compute $\thexp_+$ through the optimization \eqref{e:OptOffset}.    The value of  $r^*$ then coincides with $\theta^*$ given  \eqref{e:OptOffset}:  $r^* =  [\expa-\Lambda_0(\thexp_+G)]/\thexp_+$.

Estimates of $\oshoot(\thresh)$ was obtained via    \eqref{e:oshootEst},  based on $5\times10^5$ independent runs.    
For each of Tests~1a--1c,   the estimates are approximately decreasing in $\thresh$:   
\[
\begin{array}{c|ccccc}
 & \multicolumn{5}{c}{\widehat{\oshoot}(\thresh)}\\
\!\!
\text{Test}\backslash\thresh
 & 5 & 7 & 9 & 11 & 13\\
\hline
\rule{0pt}{3.2ex}
1\mathrm{a} & -2.08 & -1.82 & -1.69 & -1.64 & -1.64\\
1\mathrm{b} & -9.99 & -9.37 & -9.205 & -9.18 & -9.12\\
1\mathrm{c} & -5.71 & -5.32 & -5.20 & -5.13 & -5.08
\end{array}
\]


\medskip
\def\haclP{\widehat{\mathcal{P}}}

\Wham{Model 2}
This HMM model is defined by  \eqref{e:num_model2_Ztransition} and  \eqref{e:num_model2_g}.  
In all Monte-Carlo experiments to estimate CUSUM* cost and thresholds,  the Markov chain was initialized using  $\Prob\{Z_0=0\}=\Prob\{Z_0=1\}=1/2$.
 A different initialization was used in \Cref{f:model2_simplex_traj} for the sake of illustration.

\notes{I'm leaving this out since we can see the x-axis to see the range of kappa:
\\
 The cost values in \Cref{f:model2_approx} were computed on a common $\kappa$ grid with $17$ values spanning $[2,2\times 10^3]$.
}

We divide the remainder of the discussion on Model~2 into two parts, beginning with Tests~2a and 2b that are functions of the information state.

\wham{Model 2:  tests based on the information state}

\def\BellOpInfo{\textsf{B}}
\def\haBellOpInfo{\widehat{\textsf{B}}}

 \def\hafilt{\widehat{\textsf{T}}}  
\def\clD{\mathcal{D}}

\whamit{Test 2a.} 

The optimal decision region is $\clR^*  = \{ \beta \in  \clS^3  :   J^* (\beta) =  \clC_\bullet(\beta) \}$, where $ J^* \colon  \clS^3  \to\Re_+$ is the value function.  This region  is convex since $J^*$ is a concave function.   

Approximation of $J^*$ was obtained through approximation of the dynamic programming equation 
$J^*(\beta)=\min \{  \clC_\bullet(\beta) \,,   \BellOpInfo J^*\, (\beta) \} $ where the terms on the right hand side are defined as follows:  cost functions  on $ \clS^3  $ are defined in consideration of  \eqref{e:costPOMDP}: 
for   $\beta\in\clS^3$,
\[
\clC_\circ(\beta) = \beta(2)\,,   \qquad 
 \clC_\bullet(\beta) =  \kappa\sum_{  z \in \zstate_0 }   \beta(z) \Expect[\tchange\mid Z_0 = z ] 
\]
The  nonlinear filtering equation  is expressed $\condDist_{k+1}  =	\filt(\condDist_k, Y_{k+1})$ for a mapping $\filt \colon \clS^3\times \ystate \to \clS^3$,   and the Bellman operator is then
\begin{equation}
\BellOpInfo J^*\, (\beta)  = 
\clC_\circ(\beta)+\sum_{y\in\{0,1\}}p_y(\beta)\,J^*\bigl(\filt(\beta,y)\bigr) 
	\label{e:num_binned_bellman}
\end{equation}

Approximation of $\clR^*$ was based on the  selection of grid-points $\{ \beta^i : 1\le i \le N_S \}\subset \clS^3$ to construct an MDP with state space $\clX = \{1,\dots,N_S\}$ and  value function $\haJ^*\colon\clX\to\Re_+ $ defined   so that $\haJ^* (i) \approx J^*\, (\beta^i)$ for each $i$.  
For any $\beta\in\clS^3$ denote
\[
\big [ \beta  \big ] _{\clX} = \argmin_{1\le i\le N_S }    \|  \beta - \beta^i \|   
\] 
with tie-breaking performed using lexicographic ordering if the minimum is not unique.  The true value function is continuous, and hence $J^*(\beta) \approx J^*(\beta^i)$ when  $i = \big [ \beta  \big ] _{\clX} $ and the grid is sufficiently fine.   

Recalling $   \filt(\beta,y) \in \clS^3$ for any $y\in\ystate$,  $\beta\in  \clS^3$,  we define the approximate information-state dynamics via
\[
	\hafilt(\beta^i,y)	=  [ \filt(\beta^i,y) ] _{\clX}  \in \clX \,, \quad \textit{for each $  i\in \clX$.}
\]
and the approximating value function is the solution to   the dynamic programming equation,
\begin{equation}
 \begin{aligned}
\haJ^*(i) &=\min \{  \clC_\bullet(\beta^i) \,,   \haBellOpInfo \haJ^*\, (i) \} \,,
\\
\haBellOpInfo J^*\, (i) & = 
\clC_\circ(\beta^i)+\sum_{y\in\{0,1\}}p_y(\beta^i)\,\haJ^*\bigl(\hafilt(i,y)\bigr) 
 \end{aligned}
	\label{e:num_binned_bellmanApp}
\end{equation}
The region $\clR^* $ in  	\Cref{f:model2_simplex_traj}
  is  the convex hull  of  $ \{ \beta^i  :   \haBellOpInfo \haJ^*\, (i)    \le  \clC_\bullet(\beta^i) \}$.

To complete the description of this construction we describe selection of  $\{ \beta^i : 1\le i \le N_S \}$.

As seen in \Cref{f:model2_simplex_traj}, the optimal decision region is a small neighborhood of $\delta_2$ even for the small value $\kappa=2$ used to obtain this plot.  For larger $\kappa$ the regions $\clR^*$ and $\clR_{\threshSmallStar}$  are more closely concentrated near  $\delta_2$.  Consequently, non-uniform binning was designed to obtain a small value of $\|\beta^i  - \delta_2 \|$ for the majority of $i$.

One grid-point was $\beta^1 = \delta_2$.  For the remainder, 
we make the change of coordinates $r=\beta(1)/[\beta(0)+\beta(1)]$ 
and $s=\beta(2)$ for $\beta\in\clS^3$ satisfying $s\neq 1$,
and then for integers $N_s, N_r > 1$     and $0\le\ell\le  N_r$,
\[
s_j=1-(1-j/N_q)^4,\quad 1\le j\le N_s\,, 
\qquad
r_\ell=\ell/N_r   
\]
Hence the total number of grid-points is $N_S =   N_s (N_r +1) +1$.   
Numerical computations used   $N_s=500$ and $N_r=250$, resulting in  $N_S = 500\times 251+1=125,501$.

\whamit{Test 2b.} 
The optimal thresholds $\threshSmallStar(\kappa)$ were estimated by minimizing over a grid of $N=2001$ values.    
The estimates of \eqref{e:MDD+kappaMDE} used $\kmax =2\times 10^4$, and the threshold search used $\kmax =2.5\times 10^4$.

As in Test 2a we are most interested in thresholds close to unity when   $\kappa>0$ is large.    
Threshold values $\{ \threshSmall_i : 1\le i\le N \}$ were taken to be
uniformly spaced in $-\log(1-\threshSmall)$,  meaning that  $\gamma = 
- \log(1-\threshSmall_{i+1})  + \log(1-\threshSmall_{i})  = - 
\log([1-\threshSmall_{i+1}]/[1-\threshSmall_{i} ] $ is independent of $i$.    Equivalently,
\[
\threshSmall_{i+1} = 1 +   e^{-\gamma } ( \threshSmall_{i}  -1)\,  \quad 1\le i\le N
\]
This was initialized with $ \threshSmall_1 = 0.88$ to obtain values on the interval  
 $ [0.88,1-10^{-8}]$ with the majority close to unity.

\wham{Model 2:  tests based on CUSUM}

The special structure of the model allows us to obtain $\oP$ in \eqref{e:twisted} 
in the form $\oP(\bst,\bst')=\oQ(z,z')g(y'\mid z')$ where $g$ is defined in \eqref{e:num_model2_g} and $\oQ$ is defined as in 
\eqref{e:twisted}:
\[
\oQ(z,z') = \frac{1}{\delta} \frac{\xi(z')}{\xi(z)}  Q(z,z' )   \,,\quad z,z'\in\zstate_0 =\{0,1 \}\,, 
\]
where we compute $\delta=0.9836$ and $\xi=(1 ; 0.6413)$, and hence $\oQ$ and 
 its invariant pmf  are
\begin{equation}
\oQ =	\begin{bmatrix}
		0.9974 & 0.0026\\
		0.0951 & 0.9049
	\end{bmatrix},
	\qquad
	\owmetamarg\approx(0.9733,0.0267) 
	\label{e:num_model2_oP}
\end{equation}

 \Cref{t:POMDPhazard}   gives   $ \expa=-\log\delta=0.0166$,  and    for the approximate C.I.\ model,
\[
	\marg_0(1)	= \sum_{z\in\zstate_0}\owmetamarg(z)g(1\mid z)   \approx0.2160,
\] 
while by definition $\marg_1(1)=g(1\mid2)=0.8$.

\whamit{Test 2c.}

Estimates of the optimal threshold $\barthresh^*(\kappa)  $ for each $\kappa>0$ considered
 were obtained by minimizing over a uniform grid on $[3,45]$ with increment $2\times10^{-3}$.
 Each simulation used   $\kmax =4000$.

\whamit{Test 2d.}
The optimal threshold $\barthresh^*(\kappa)  $  were estimated for each
 $\del\in\{2,6,11\}$ using the same grid and $\kmax $ as Test~2c.  
 Based on 
Monte-Carlo estimates of the expectations in \eqref{e:DriftAssumptionsGen}  we find that 
the sign constraints $m_0<0<m_1$ hold for each value of $\del$ considered:
\[
\begin{array}{c|cc}
	\del & m_0 & m_1 \\ \hline
	\rule{0pt}{3.2ex}
	2  & -0.773 & 0.700\\
	6  & -0.779 & 0.288\\
	11 & -0.780  & 0.148
\end{array}
\]

\whamit{Test 2e.}
For   $\theta\in\Re^4$ consider the  positive matrix with entries in $\bigstate$ defined 
for $\bst = (z,y)$,  $\bst'=(z',y')$ via
\[
\begin{aligned}
K_{\theta} (\bst,\bst') 
&=
\oP(\bst,\bst')  \exp \big( F_\theta (y,y') ) \big)
\\
&=
\oQ(z,z')g(y'\mid z') \exp \big( F_\theta (y,y') ) \big) 
\end{aligned}
\]
We then have   $\Lambda_0(F_\theta)=\log\rho_\theta $,  where  $\rho_\theta>0$ denotes the P-F eigenvalue of $K_{\theta}$, allowing us to compute 	$ \Gamma_0(\theta)   \eqdef  \Lambda_0(F_\theta)-\marg_1(F_\theta)$ for any $\theta$.



The vector $\theta^*=   (-0.236 ; 0.245 ;  0.013;  0.123)$,  solves    the primal  \eqref{e:jStarConvexRelaxationFD}, and 
\[
\begin{array}{c|c|c|c|c|c}
m_0^{\theta^*} & m_1^{\theta^*} & \thexp_+^{\theta^*} & \gamma^2 & \bFine & \oshoot_\infty	 \\
	\hline
	\rule{0pt}{3.2ex}
   -0.099& 0.1107& 1.000& 16.08& 0.107 &-5.18
\end{array}
\]
from which we obtain $\barthreshFineStar(\kappa) $ and $ \barJFineStar(\kappa) $ defined in \eqref{e:barFineBody}.

As in Model~1 the   estimates \eqref{e:oshootEst} were obtained
 based on $5\times10^5$ independent runs for a range of threshold values.     
 The results are nearly constant over the range of thresholds considered:  
\[
\begin{array}{c|rrrrrr}
\thresh
    & 11  & 15 & 18  & 24\\
\hline
\rule{0pt}{3.2ex}
\widehat{\oshoot}(\thresh)
    & -5.183  & -5.163 & -5.197  & -5.168
\end{array}
\]

The twisted distribution $  \twmarg_0^{[2]} $ associated with the function $F_{\theta^*}\colon\ystate^2\to\Re$ was computed and   we obtained    $\|  \twmarg_0^{[2]}  -  \marg_1^{[2]} \|_\infty
\approx  10^{-8}$, numerically verifying \eqref{e:model2bivariateMatching}  (a special case of \eqref{e:calcOptLin}).

\notes{  kmax etc is now summarized at the top:
\\
CUSUM* for $F_{\theta^*}$ used a uniform threshold grid on $[0.5,9.0]$ with increment $5\times10^{-4}$ and $\kmax =2\times 10^3$. 
}

\wham{Proof of \Cref{t:Model2clG}:}
We have  $
F_w(y,y') = w^\transpose\psi(y,y') = y'-y$,  for any $y,y'\in\ystate$,
by construction of the vector $w$ and the basis $\psi$.
\Cref{t:Null2} gives $\Lambda_0(F_{\theta+sw}) = \Lambda_0(F_\theta)$ for any $s$,  and   stationarity gives  $\marg_i^{[2]} (F_{\theta+sw}) =\marg_i^{[2]} (F_{\theta})$ for each $i$.  From the definitions we conclude that  
\[
\Gamma_0(\theta+sw)=\Gamma_0(\theta)\, \quad s\in\Re
\]
The proof is completed on recalling that   $F_{\theta+rv} = F_{\theta} + r$ and 
hence
 $\Gamma_0(\theta+rv)= \Gamma_0(\theta) $ for $r\in\Re$.    
\qed

\notes{sm2ac:  Here is an example where you are bringing in concepts outside of the current context:
\\
$\theta\mapsto\theta+sw$ leaves the objective and both constraints
in \eqref{e:jStarConvexRelaxationFD} unchanged.  
\\
The lemma statement is designed to avoid reference  \eqref{e:jStarConvexRelaxationFD}.   It doesn't matter that the lemma will eventually be applied to  \eqref{e:jStarConvexRelaxationFD}.    
 }

 \end{document}

%% file: bookmacrosRL_IEEE.tex
\def\urls#1{{\footnotesize\url{#1}}}

\def\mindex#1{\index{#1}}

\DeclareFontFamily{U}{mathx}{\hyphenchar\font45}
\DeclareFontShape{U}{mathx}{m}{n}{<-> mathx10}{}
\DeclareSymbolFont{mathx}{U}{mathx}{m}{n}
\DeclareMathAccent{\widebar}{0}{mathx}{"73}

\makeatletter
\newcommand\gobblepars{%
    \@ifnextchar\par%
        {\expandafter\gobblepars\@gobble}%
{}}
\makeatother

\def\whamrm#1{\smallbreak\pagebreak[3]%
	\noindent\text{\rm#1}\ \ \gobblepars}
	
\def\whamit#1{\smallbreak\pagebreak[3]%
	\noindent\textit{#1}\ \ \gobblepars}

\def\wham#1{\smallbreak\pagebreak[3]%
	\noindent\textup{\textbf{#1}}\ \ \gobblepars}

\def\fee{\upphi}

\newcommand{\bbblot}{\raise1pt\hbox{\vrule height .4ex width .4ex depth .05ex}}

\long\def\defbox#1{\framebox[.9\hsize][c]{\parbox{.85\hsize}{%
\parindent=0pt
\baselineskip=12pt plus .1pt      
\parskip=6pt plus 1.5pt minus 1pt 
 #1}}}

\long\def\beginbox#1\endbox{\subsection*{}%
\hbox{\hspace{.05\hsize}\defbox{\medskip#1\bigskip}}%
\subsection*{}}

\def\endbox{}

 \def\archival#1{} 

\def\FRAC#1#2#3{\genfrac{}{}{}{#1}{#2}{#3}}

\def\ddt{{\mathchoice{\FRAC{1}{d}{dt}}%
{\FRAC{1}{d}{dt}}%
{\FRAC{3}{d}{dt}}%
{\FRAC{3}{d}{dt}}}}

\def\ddtp{{\mathchoice{\FRAC{1}{d^{\hbox to 2pt{\rm\tiny +\hss}}}{dt}}%
{\FRAC{1}{d^{\hbox to 2pt{\rm\tiny +\hss}}}{dt}}%
{\FRAC{3}{d^{\hbox to 2pt{\rm\tiny +\hss}}}{dt}}%
{\FRAC{3}{d^{\hbox to 2pt{\rm\tiny +\hss}}}{dt}}}}

\def\ddyp{{\mathchoice{\FRAC{1}{d^{\hbox to 2pt{\rm\tiny +\hss}}}{dy}}%
{\FRAC{1}{d^{\hbox to 2pt{\rm\tiny +\hss}}}{dy}}%
{\FRAC{3}{d^{\hbox to 2pt{\rm\tiny +\hss}}}{dy}}%
{\FRAC{3}{d^{\hbox to 2pt{\rm\tiny +\hss}}}{dy}}}}

\def\half{{\mathchoice{\FRAC{1}{1}{2}}%
{\FRAC{1}{1}{2}}%
{\FRAC{3}{1}{2}}%
{\FRAC{3}{1}{2}}}}

\def\limsup{\mathop{\rm lim{\,}sup}}

\def\argmin{\mathop{\rm arg{\,}min}}
\def\argmax{\mathop{\rm arg{\,}max}}

\def\ustate{{\sf U}} 
\def\ystate{{\sf Y}} 
\def\zstate{{\sf Z}}

\def\ystate{{\sf Y}}

\def\by{\clB(\ystate)}

\def\cP{{\check{P}}}

\def\bfmath#1{{\mathchoice{\mbox{\boldmath$#1$}}%
{\mbox{\boldmath$#1$}}%
{\mbox{\boldmath$\scriptstyle#1$}}%
{\mbox{\boldmath$\scriptscriptstyle#1$}}}}

\def\bfPhi{\bfmath{\Phi}}

\def\bfmB{\bfmath{B}}

\def\bfmX{\bfmath{X}}

\def\bfmY{\bfmath{Y}}

\def\bfmhhaY{\bfmath{\hhaY}} 
\def\bfmhhaY{\hbox to 0pt{$\widehat{\bfmY}$\hss}\widehat{\phantom{\raise 1.25pt\hbox{$\bfmY$}}}}

\def\bfmZ{\bfmath{Z}}

\def\haP{{\widehat P}}

\def\clB{{\cal B}}
\def\clC{{\cal C}}

\def\clG{{\cal G}}

\def\clI{{\cal I}}
\def\clK{{\cal K}}
\def\clL{{\cal L}}

\def\clR{{\cal R}}
\def\clS{{\cal S}}

\def\clV{{\cal V}}

\def\clX{{\cal X}}
\def\clY{{\cal Y}}

\def\clL{{\cal L}}
\def\clI{{\cal I}}
\def\clC{{\cal C}}

\def\eqdef{\mathbin{:=}}

\def\Prob{{\sf P}}

\def\Expect{{\sf E}}

\def\ind{\bbbone}
  
 \def\epsy{\varepsilon}

\def\varble{\,\cdot\,}

\def\formtmp#1#2{{\vskip12pt\noindent\fboxsep=0pt\colorbox{#1}{\vbox{\vskip3pt\hbox to \textwidth{\hskip3pt\vbox{\raggedright\noindent\textbf{#2\vphantom{Qy}}}\hfill}\vspace*{3pt}}}\par\vskip2pt%
\noindent\kern0pt}}

\def\barJ{{\bar{J}}}

\def\barpsi{{\bar{\psi}}}

{\end{list}}

\def\ass(#1:#2){(#1\ref{#1:#2})}

\def\ritem#1{
\item[{\sf \ass(\current_model:#1)}]
}

\newenvironment{recall-ass}[1]{%
\begin{description}
\def\current_model{#1}}{
\end{description}
}

 \newcommand{\blot}{\vrule height 1.1ex width .9ex depth -.1ex }
\def\qedb{\ifmmode\blot\else{\vspace{-.2cm}\unskip\nobreak\hfil
\penalty50\hskip1em\null\nobreak\hfil\blot
\parfillskip=0pt\finalhyphendemerits=0\endgraf}\fi}

  \makeatletter
\DeclareRobustCommand{\sqcdot}{\mathbin{\mathpalette\morphic@sqcdot\relax}}
\newcommand{\morphic@sqcdot}[2]{%
  \sbox\z@{$\m@th#1\centerdot$}%
  \ht\z@=.33333\ht\z@
  \vcenter{\box\z@}%
}
\makeatother

\def\qedsymbol{\hbox{\tiny$\blacksquare$}}  %

\def\qed{\ifmmode\qedsymbol\else{\unskip\nobreak\hfil
\penalty50\hskip1em\null\nobreak\hfil\qedsymbol
\parfillskip=0pt\finalhyphendemerits=0\endgraf}\fi}

\newcounter{rmnum}

\newcounter{anum}

\newcommand{\field}[1]{\mathbb{#1}}

\def\Re{\field{R}}

\def\Prob{{\sf P}}

\def\Expect{{\sf E}}

\def\transpose{{\intercal}}

\def\argmin{\mathop{\rm arg\, min}}
\def\ind{\hbox{\large \bf 1}}

\def\epsy{\varepsilon}
\def\varble{\,\cdot\,} 
\def\dint{\mathop{\int\!\!\!\int}} 

\def\haJ{\widehat J}

\def\haY{\widehat{Y}}

\def\hhaY{\hbox to 0pt{$\haY$\hss}\widehat{\phantom{\raise 1.25pt\hbox{Y}}}}

\def\haY{\widehat Y}

\def\bfPhi{\bfmath{\Phi}}

\newlength{\dhatheight}
